\documentclass{article}
\usepackage{geometry}
\usepackage{amssymb}
\usepackage{amsmath,cases}
\usepackage{amsthm}
\usepackage{enumerate}
\usepackage{xcolor}
\usepackage{graphicx}
\usepackage{textcomp}
\usepackage{mathrsfs}
\usepackage{bbm}
\usepackage{bm}
\usepackage{caption}
\usepackage{chemarrow}
\usepackage{stmaryrd}
\usepackage{url}
\usepackage{psfrag}
\usepackage[active]{srcltx}   
\usepackage[all,cmtip]{xy}
\usepackage{pdfsync}																				
\usepackage{yhmath}
\usepackage[bookmarksnumbered=true, pdfauthor={Jérémie Bettinelli}]{hyperref}									

\newcommand{\hr}{\texorpdfstring}

\newcommand{\8}{\infty}
\newcommand{\choo}[2]{\begin{pmatrix}#1 \\ #2 \end{pmatrix}}

\newcommand{\bs}{\backslash}
\newcommand{\eps}{\varepsilon}
\newcommand{\ent}[2]{\llbracket #1,#2 \rrbracket}
\newcommand{\de}{\mathrel{\mathop:}\hspace*{-.6pt}=}

\newcommand{\ooo}{\bullet}
\newcommand{\ul}[1]{\underline #1}
\newcommand{\ol}[1]{\overline #1}

\newcommand{\ton}{\xrightarrow[n\to\infty]{}}
\newcommand{\tok}{\xrightarrow[k\to\infty]{}}
\newcommand{\tol}{\xrightarrow[n \to \infty]{(d)}}
\newcommand{\tolk}{\xrightarrow[k \to \infty]{(d)}}

\newcommand{\tode}[3]{\xrightarrow[#2 \to #3]{#1}}

\newcommand{\bbB}{\mathbb{B}}
\newcommand{\Ddd}{\mathbb{D}}
\newcommand{\ee}{\mathbbm{e}}
\newcommand{\E}[1]{\mathbb{E}\left[ #1 \right]}
\newcommand{\MM}{\mathbb{M}}
\newcommand{\N}{\mathbb{N}}
\newcommand{\Pb}{\mathbb{P}}
\newcommand{\bbQ}{\mathbb{Q}}
\newcommand{\R}{\mathbb{R}}
\newcommand{\Sss}{\mathbb{S}}

\newcommand{\Z}{\mathbb{Z}}

\newcommand{\A}{\mathcal{A}}
\newcommand{\cB}{\mathcal{B}}
\newcommand{\C}{\mathcal{C}}

\newcommand{\EE}{\mathcal{E}}

\newcommand{\calG}{\mathcal{G}}

\newcommand{\K}{\mathcal{K}}

\newcommand{\cN}{\mathcal{N}}

\newcommand{\Q}{\mathcal{Q}}

\newcommand{\Qns}{\mathcal{Q}_{n,\sigma}}
\newcommand{\Qnsn}{\mathcal{Q}_{n,\sigma_n}}
\newcommand{\Qnsb}{\mathcal{Q}_{n,\sigma}^\bullet}
\newcommand{\Qnsnb}{\mathcal{Q}_{n,\sigma_n}^\bullet}
\newcommand{\rR}{\mathcal{R}}
\newcommand{\cS}{\mathcal{S}}

\newcommand{\V}{\mathcal{V}}

\newcommand{\X}{\mathcal{X}}
\newcommand{\ZZ}{\mathcal{Z}}

\newcommand{\B}{\mathscr{B}}
\newcommand{\Bs}{\B_\sigma}
\newcommand{\Bsn}{\B_{\sigma_n}}

\newcommand{\FF}{\mathscr{F}}
\newcommand{\Fi}{\mathscr{F}_\infty}

\newcommand{\Leb}{\mathscr{L}}

\newcommand{\loo}{\mathscr{L}}
\newcommand{\PP}{\mathscr{P}}
\newcommand{\QQ}{\mathscr{Q}}

\newcommand{\TT}{\mathscr{T}}

\newcommand{\aaa}{\mathfrak{a}}
\newcommand{\bb}{\mathfrak{b}}
\newcommand{\fB}{\mathfrak{B}}

\newcommand{\e}{\mathfrak{e}}
\newcommand{\f}{\mathfrak{f}}
\newcommand{\fF}{\mathfrak{F}}

\newcommand{\lab}{\mathfrak{l}}
\newcommand{\Lab}{\mathfrak{L}}
\newcommand{\Li}{\mathfrak{L}_\infty}
\newcommand{\m}{\mathfrak{m}}

\newcommand{\q}{\mathfrak{q}}

\newcommand{\qis}{\mathfrak{q}_\infty^\sigma}
\newcommand{\fR}{\mathfrak{R}}

\newcommand{\tr}{\mathfrak{t}}

\newcommand{\oo}{\bm{0}}

\newcommand{\bmp}{\bm{p}}
\newcommand{\bmT}{\bm{T}}

\newcommand{\cov}{\operatorname{cov}}
\newcommand{\law}{\stackrel{(d)}{=}}
\newcommand{\un}[1]{\mathbbm{1}_{#1}}

\newcommand{\eqt}{\simeq_\infty}
\newcommand{\ro}{\partial}
\newcommand{\suc}{\operatorname{succ}}

\newcommand{\IP}{\textup{IP}}
\newcommand{\IPl}{\textup{IP}_{left}}
\newcommand{\IPr}{\textup{IP}_{right}}
\newcommand{\fl}{f\hspace{-0.6mm}l}
\newcommand{\fli}{\fl_\infty}

\newcommand{\dish}{\delta_\mathcal{H}}
\newcommand{\dH}{\operatorname{dim}_\mathcal{H}}
\newcommand{\diam}{\operatorname{diam}}
\newcommand{\dGH}{d_{GH}}
\newcommand{\dis}{\operatorname{dis}}

\newcommand{\disig}{d_\infty^\sigma}
\newcommand{\pii}{\pi_\infty}
\newcommand{\PM}{\textup{P}\mathbb{M}}

\newcommand{\BDG}{Bouttier--Di~Francesco--Guitter }
\newcommand{\g}{\gamma \, n^{1/4}}

\newcommand{\arc}{\smallfrown}
\newcommand{\arcr}{\curvearrowright}
\newcommand{\arcl}{\curvearrowleft}

\newcommand{\sulW}{{\rlap{\underline{\hspace{2mm}}}\widehat W\!}}
\newcommand{\ulz}{{\rlap{\underline{\hspace{1.7mm}}}\zeta}}

\newcommand{\sand}{\qquad\text{ and }\qquad}

\newcommand{\lhb}{[[}

\newcommand{\rhb}{]]}

\newcommand{\lp}{\left(}
\newcommand{\rp}{\right)}
\newcommand{\lb}{\left\{}
\newcommand{\rb}{\right\}}
\newcommand{\lf}{\left\lfloor}
\newcommand{\rf}{\right\rfloor}
\newcommand{\lc}{\left\lceil}
\newcommand{\rc}{\right\rceil}
\newcommand{\lt}{\left|}
\newcommand{\rt}{\right|}

\newcommand{\lbr}{\left[}
\newcommand{\rbr}{\right]}

\newcommand{\rno}{\right.}
\newcommand{\lh}{\|}
\newcommand{\rh}{\|}

\definecolor{gris}{gray}{0.7}
\definecolor{grisf}{gray}{0.4}
\definecolor{vert}{rgb}{0,.5547,0}

\theoremstyle{plain}
\newtheorem{thm}{Theorem}
\newtheorem{lem}[thm]{Lemma}
\newtheorem{prop}[thm]{Proposition}

\newtheorem{corol}[thm]{Corollary}
\newtheorem{defi}{Definition}

\theoremstyle{definition}
\newtheorem*{rem}{Remark}
\newtheorem*{nota}{Note}

\theoremstyle{remark}

\newenvironment{pre}[1][\proofname]{%
  \proof[\bfseries #1]%
}{\endproof}

\title{Scaling Limit of Random Planar Quadrangulations with a Boundary\thanks{This work is partially supported by ANR-08-BLAN-0190}}
\author{J\'er\'emie \textsc{Bettinelli}\\
	\small Institut \'Elie Cartan de Lorraine, B.P.\ 239\\
	\small F-54506 Vand\oe uvre-l\`es-Nancy Cedex\\
	\small\url{jeremie.bettinelli@normalesup.org}\\
	\small\url{http://www.normalesup.org/~bettinel}}

\begin{document}
\maketitle

\begin{abstract}

We discuss the scaling limit of large planar quadrangulations with a boundary whose length is of order the square root of the number of faces. We consider a sequence $(\sigma_n)$ of integers such that $\sigma_n/\sqrt{2n}$ tends to some $\sigma\in[0,\infty]$. For every $n \ge 1$, we call $\q_n$ a random map uniformly distributed over the set of all rooted planar quadrangulations with a boundary having $n$ faces and $2\sigma_n$ half-edges on the boundary. For $\sigma\in (0,\infty)$, we view $\q_n$ as a metric space by endowing its set of vertices with the graph metric, rescaled by $n^{-1/4}$. We show that this metric space converges in distribution, at least along some subsequence, toward a limiting random metric space, in the sense of the Gromov--Hausdorff topology. We show that the limiting metric space is almost surely a space of Hausdorff dimension $4$ with a boundary of Hausdorff dimension $2$ that is homeomorphic to the two-dimensional disc. For $\sigma=0$, the same convergence holds without extraction and the limit is the so-called Brownian map. For $\sigma=\infty$, the proper scaling becomes $\sigma_n^{-1/2}$ and we obtain a convergence toward Aldous's CRT.

\end{abstract}

\section{Introduction}

\subsection{Motivations}

In the present work, we investigate the scaling limit of random (planar) quadrangulations with a boundary. Recall that a planar map is an embedding of a finite connected graph (possibly with loops and multiple edges) into the two-dimensional sphere, considered up to direct homeomorphisms of the sphere. The faces of the map are the connected components of the complement of edges. A quadrangulation with a boundary is a particular instance of planar map whose faces are all \textit{quadrangles}, that is, faces incident to exactly~$4$ half-edges (or oriented edges), with the exception of one face of arbitrary even degree. The quadrangles will be called \textit{internal faces} and the other face will be referred to as the \textit{external face}. The half-edges incident to the external face will constitute the \textit{boundary} of the map. Beware that we do not require the boundary to be a simple curve. We will implicitly consider our maps to be rooted, which means that one of the half-edges is distinguished. In the case of quadrangulations with a boundary, the root will always lie on the boundary, with the external face to its left. 

In recent years, scaling limits of random maps have been the subject of many studies. The most natural setting is the following. We consider maps as metric spaces, endowed with their natural graph metric. We choose uniformly at random a map of ``size''~$n$ in some class, rescale the metric by the proper factor, and look at the limit in the sense of the Gromov--Hausdorff topology~\cite{gromov99msr}. The size considered is usually the number of faces. From this point of view, the most studied class is the class of planar quadrangulations. The pioneering work of Chassaing and Schaeffer~\cite{chassaing04rpl} revealed that the proper rescaling factor in this case is~$n^{-1/4}$. The problem was first addressed by Marckert and Mokkadem~\cite{marckert06limit}, who constructed a candidate limiting space called the \textit{Brownian map}, and showed the convergence toward it in another sense. Le~Gall~\cite{legall07tss} then showed the relative compactness of this sequence of metric spaces and that any of its accumulation points was almost surely of Hausdorff dimension~$4$. It is only recently that the problem was completed independently by Miermont~\cite{miermont11bms} and Le~Gall~\cite{legall11ubm}, who showed that the scaling limit is indeed the Brownian map. This last step, however, is not mandatory in order to identify the topology of the limit: Le~Gall and Paulin~\cite{legall08slb}, and later Miermont~\cite{miermont08sphericity}, showed that any possible limit is homeomorphic to the two-dimensional sphere. 

To be a little more accurate, Le~Gall considered in~\cite{legall07tss} the classes of $2p$-angulations, for $p\ge2$ fixed, and, in~\cite{legall11ubm}, the same classes to which he added the class of triangulations, so that the result about quadrangulations is in fact a particular case. We also generalized the study of~\cite{legall07tss,miermont08sphericity} to the case of bipartite quadrangulations in fixed positive genus $g \ge 1$ in~\cite{bettinelli10slr,bettinelli10tsl}, where we showed the convergence up to extraction of a subsequence and identified the topology of any possible limit as that of the surface of genus~$g$. In the present work, we adopt a similar point of view and consider the class of quadrangulations with a boundary, where the length of the boundary grows as the square root of the number of internal faces. We show the convergence up to extraction, and show that any possible limiting space is almost surely a space of Hausdorff dimension~$4$ with a boundary of Hausdorff dimension~$2$ that is homeomorphic to the two-dimensional disc. We also show that, if the length of the boundary is small compared to the square root of the number of internal faces, then the convergence holds (without extraction) and the limit is the Brownian map. When the length of the boundary is large with respect to the square root of the number of internal faces, then the proper scaling becomes the length of the boundary raised to the power $-1/2$, and we obtain a convergence toward the so-called Continuum Random Tree (CRT).

The study of these problems often starts with a bijection between the class considered and a class of simpler objects. In the case of planar quadrangulations, the bijection in question is the so-called Cori--Vauquelin--Schaeffer bijection~\cite{cori81planar,schaeffer98cac,chassaing04rpl} between planar quadrangulations and so-called well-labeled trees. This bijection has then been generalized in several ways. Bouttier, Di Francesco, and Guitter~\cite{bouttier04pml} extended it into a bijection coding all planar maps (and even more). Later, Chapuy, Marcus, and Schaeffer~\cite{chapuy07brm} considered bipartite quadrangulations of positive fixed genus. As quadrangulations with a boundary are a particular case of planar maps, we will use in this work a slightly amended instance of the \BDG bijection. Let us also mention that Bouttier and Guitter studied in~\cite{bouttier09dsqb} the distance statistics of quadrangulations with a boundary. In particular, their study showed the existence of the three different regimes we consider in this work. Additionally, Curien and Miermont~\cite{curien12uipqb} studied in a recent work the local limit of quadrangulations with a boundary. 

From now on, when we speak of quadrangulations, we always mean rooted planar quadrangulations with a boundary, and, by convention, we always draw the external face as the infinite component of the plane.

\subsection{Main results}

\subsubsection{Generic case}

Let~$\m$ be a map. We call $V(\m)$ its sets of vertices, $E(\m)$ its sets of edges, and $\vec E(\m)$ its set of half-edges. We say that a face~$f$ is \textbf{incident} to a half-edge~$\e$ (or that~$\e$ is incident to~$f$) if~$\e$ belongs to the boundary of~$f$ and is oriented in such a way that~$f$ lies to its left. We write~$\e_*$ the root of~$\m$, and, for any half-edge~$\e$, we call~$\bar\e$ its reverse, as well as~$\e^-$ and~$\e^+$ its origin and end. We denote by~$d_\m$ the graph metric on~$\m$ defined as follows: for any $a,b\in V(\m)$, the distance $d_\m(a,b)$ is the number of edges of any shortest path in~$E(\m)$ linking~$a$ to~$b$. Finally, we call $\Qns$ the set of all quadrangulations with a boundary having~$n$ internal faces and~$2\sigma$ half-edges on the boundary.

The Gromov--Hausdorff distance between two compact metric spaces $(\X,\delta)$ and $(\X',\delta')$ is defined by
$$\dGH\lp(\X,\delta),(\X',\delta')\rp \de \inf \big\{ \dish\big(\varphi(\X),\varphi'(\X')\big)\!\big\},$$
where the infimum is taken over all isometric embeddings $\varphi : \X \to \X''$ and $\varphi':\X'\to \X''$ of~$\X$ and~$\X'$ into the same metric space $(\X'', \delta'')$, and $\dish$ stands for the usual Hausdorff distance between compact subsets of $\X''$. This defines a metric on the set~$\MM$ of isometry classes of compact metric spaces (\cite[Theorem~7.3.30]{burago01cmg}), making it a Polish space\footnote{This is a simple consequence of Gromov's compactness theorem \cite[Theorem~7.4.15]{burago01cmg}.}.

Our main results for quadrangulations with a boundary are the following.

\begin{thm}\label{cvqb}
Let $\sigma \in (0,\8)$ and $(\sigma_n)_{n\ge 1}$ be a sequence of positive integers such that $\sigma_n / \sqrt{2n} \to \sigma$ as $n\to \infty$. Let $\q_n$ be uniformly distributed over the set $\Qnsn$ of all planar quadrangulations with a boundary having~$n$ internal faces and $2\sigma_n$ half-edges on the boundary. Then, from any increasing sequence of integers, we may extract a subsequence $(n_k)_{k\ge 0}$ such that there exists a random metric space $(\qis,\disig)$ satisfying
$$\Big( V(\q_{n_k}),\frac 1 {\gamma n_k^{1/4}}\, d_{\q_{n_k}} \Big) \tolk (\qis,\disig)$$
in the sense of the Gromov--Hausdorff topology, where
$$\gamma \de \lp\frac {8}9 \rp^{1/4}.$$

Moreover, the Hausdorff dimension of the limit space $(\qis,\disig)$ is almost surely equal to~$4$, regardless of the choice of the sequence of integers.
\end{thm}

Remark that the constant~$\gamma$ is not necessary in this statement (simply change $\disig$ into $\gamma \disig$). We
made it figure at this point for consistency with the other works on the subject and because of our definitions later in the
paper. Note also that, a priori, the metric space $(\qis,\disig)$ depends on the subsequence $(n_k)_{k\ge 0}$. In view of the recent developments made by Miermont~\cite{miermont11bms} and Le~Gall~\cite{legall11ubm} in the case without boundary, we conjecture that the extraction in Theorem~\ref{cvqb} is not necessary and that~$\disig$ can be explicitly expressed in a way similar to their expression. We also believe that the space $(\qis,\disig)$ only depends on~$\sigma$, and arises as some universal scaling limit for more general classes of random maps with a boundary. In particular, our approach should be generalizable to the case of $2p$-angulations, $p\ge 2$, by using the same kind of arguments as Le Gall in~\cite{legall07tss}.

As in the case without boundary, Theorem~\ref{cvqb} is nonetheless sufficient to identify the topology of the limit, regardless of the subsequence $(n_k)_{k\ge 0}$.

\begin{thm}\label{cvqtop}
For $\sigma >0$, any possible metric space $(\qis,\disig)$ from Theorem~\ref{cvqb} is a.s.\ homeomorphic to the $2$-dimensional disc~$\Ddd_2$.
\end{thm}

We may also compute the Hausdorff dimension of the boundary of the limiting space: we define $\partial \qis \subseteq \qis$ as the set of points having no neighborhood homeomorphic to a disc.

\begin{thm}\label{thmdimh}
For any $\sigma>0$, the boundary $\partial \qis$ is a subset of $(\qis, \disig)$ whose Hausdorff dimension is almost surely equal to~$2$.
\end{thm}

\subsubsection{Case \hr{$\sigma=0$}{s=0}}

In the case where $\sigma=0$, we may actually be a little more precise than in the previous theorems. In particular, we have a whole convergence, instead of just a convergence along subsequences. We find that, in the limit, the boundary ``vanishes'' in the sense that we obtain the same limit as in the case without boundary: the Brownian map~\cite{legall11ubm,miermont11bms}.

\begin{thm}\label{cvqbs0}
Let $(\sigma_n)_{n\ge 1}$ be a sequence of positive integers such that $\sigma_n / \sqrt{2n} \to 0$ as $n\to \infty$. Let $\q_n$ be uniformly distributed over the set $\Qnsn$ of all planar quadrangulations with a boundary having~$n$ internal faces and $2\sigma_n$ half-edges on the boundary. Then, 
$$\Big( V(\q_{n}),\frac 1 {\gamma n^{1/4}}\, d_{\q_{n}} \Big) \tol (\m_\8,D^*)$$
in the sense of the Gromov--Hausdorff topology, where $(\m_\8,D^*)$ is the Brownian map.
\end{thm}

As a consequence, we retrieve immediately the classical properties of the Brownian map, from which the results of the previous section are inspired. For instance, it is known that the Hausdorff dimension of $(\m_\8,D^*)$ is almost surely equal to~$4$ (\cite{legall07tss}), and that the metric space $(\m_\8,D^*)$ is a.s.\ homeomorphic to the $2$-dimensional sphere $\Sss_2$ (\cite{legall08slb,miermont08sphericity}).

\subsubsection{Case \hr{$\sigma=\infty$}{sigma=infinity}}

In the case $\sigma=\8$, the proper scaling factor is no longer~$n^{-1/4}$, but the length of the boundary raised to the power~$-1/2$. We find Aldous's so-called CRT~\cite{aldous91crt,aldous93crt} defined as follows. We denote the normalized Brownian excursion by~$\ee$, and we define the pseudo-metric
$$\delta_\ee(s,t)\de \ee(s)+\ee(t)-2 \min_{[s\wedge t, s \vee t]} \ee,\qquad 0 \le s,t \le 1.$$
It defines a metric on the quotient $\TT_\ee \de [0,1]_{/\{\delta_\ee=0\}}$, which, by a slight abuse of notation, we still write~$\delta_\ee$. The \textbf{Continuum Random Tree} is the random metric space $(\TT_\ee,\delta_\ee)$. Moreover, we also have a whole convergence in this case.

\begin{thm}\label{cvqbs8}
Let $(\sigma_n)_{n\ge 1}$ be a sequence of positive integers such that $\sigma_n / \sqrt{2n} \to \8$ as $n\to \infty$. Let $\q_n$ be uniformly distributed over the set $\Qnsn$ of all planar quadrangulations with a boundary having~$n$ internal faces and $2\sigma_n$ half-edges on the boundary. Then, 
$$\Big( V(\q_{n}),\frac 1 {(2\sigma_n)^{1/2}}\, d_{\q_{n}} \Big) \tol (\TT_\ee,\delta_\ee)$$
in the sense of the Gromov--Hausdorff topology.
\end{thm}

Let us try to give an intuition of what happens here. Roughly speaking, the boundary takes so much space that we need to rescale by a factor that suits its length. The faces, which should be in the scale~$n^{1/4}$, are then too much rescaled and disappear in the limit, leaving only the boundary visible. As a result, for~$n$ large enough and in the proper scale, the quadrangulation itself is not very far from its boundary, which in its turn is not very far from a random tree. This rough reasoning gives an intuition of why the CRT arises at the limit.

We also observe an interesting phenomenon if we take all these theorems into account. It can be expected that, if we take a uniform quadrangulation in~$\Qnsn$ with~$n$ large and~$\sigma_n$ large enough but not too large (probably in the scale $n^{1/2 + \eps}$ with $\eps >0$ small) then, in the scale~$n^{-1/4}$, it should locally resemble the Brownian map, whereas in the scale~$\sigma_n^{-1/2}$, it should look more like the CRT. We believe this picture could be turned into a rigorous statement but we choose not to pursue this route in the present paper.

\subsection{Organization of this paper and general strategy}

We begin by exposing in Section~\ref{secbdg} the version of the \BDG bijection that we will need. As we do not use it in its usual setting, we spend some time explaining it. In particular, we introduce a notion of bridge that is not totally standard. We then investigate in Section~\ref{secpt1} the scaling limit of the objects appearing in this bijection, and deduce Theorem~\ref{cvqb}.

Discrete forests play an important part in the coding of quadrangulations with a boundary through the \BDG bijection, and the analysis of Section~\ref{secpt1} leads to the construction of a continuum random forest, which may be seen as a generalization of Aldous's CRT~\cite{aldous91crt,aldous93crt}. We carry out the analysis of Le~Gall~\cite{legall07tss} to our case in Section~\ref{secrf} and see any limiting space of Theorem~\ref{cvqb} as a quotient of this continuum random forest via an equivalence relation defined in terms of Brownian labels on it.

Following Miermont~\cite{miermont08sphericity}, we then prove Theorem~\ref{cvqtop} in Section~\ref{secsurfb} thanks to the notion of regularity introduced by Whyburn~\cite{whyburn35sls,whyburn35rcm}. As we consider in this work surfaces with a boundary, the notion of 1-regularity used by Miermont in~\cite{miermont08sphericity} is no longer sufficient: we will also need here the notion of 0-regularity, which we will present in Section~\ref{secreg}.

Section~\ref{sig0} is devoted to the case $\sigma=0$ in which we use a totally different approach, consisting in comparing quadrangulations with a ``small'' boundary with quadrangulations without boundary. In Section~\ref{sig8}, we treat the case $\sigma=\8$ by a different method.

We will need to use the so-called Brownian snake to prove some remaining technical results. We prove these in Section~\ref{secsnake}. In particular, in Section~\ref{secippre}, we will look at the increase points of the Brownian snake we consider. From our approach, we can retrieve \cite[Lemma~3.2]{legall08slb}.

Finally, Section~\ref{secdev} is devoted to some developments and open questions.

\bigskip

Our general strategy is in many points similar to~\cite{bettinelli10slr,bettinelli10tsl}. Although we will try to make this work as self-contained as possible, we will often refer the reader to these papers when the proofs are readily adaptable, and will rather focus on the new ingredients. One of the main difficulties that was not present in~\cite{bettinelli10slr,bettinelli10tsl} arises from the fact that the Brownian labels on the continuum random forest we construct do not have the same diffusion factor on the floor than in the trees. To be a little more precise, the labels in the trees vary like standard Brownian motion, whereas on the floor they vary as a Brownian motion multiplied by the factor~$\sqrt 3$ (see Proposition~\ref{cvcont} for a rigorous statement). This factor comes from the fact that the bridge coding the external face in the \BDG bijection does not have the same variance as the Motzkin paths appearing everywhere else. Its presence generates new technical issues and forces us to find new proofs for some of Le~Gall's estimates.

A key point of our analysis is that, at the limit, the boundary does not have any pinch points (Lemma~\ref{lem0reg}). As the boundary of the map roughly corresponds to the floor of the forest (Proposition~\ref{propboundary}), it will be crucial to see that, in the quotient we define, the points of the floor are not identified with one another (Lemma~\ref{lemid}). We will see in Theorem~\ref{ipb} that two points are identified if they have the same labels and if the labels of the points ``between them'' are all greater. From the already known cases, we could think that everything will work similarly but, a priori, this factor $\sqrt 3$ could induce some identification of points on the floor of the forest. Fortunately, this does not happen. However, we can see from our proofs in Section~\ref{secippre} that this value is critical, in the sense that if it was strictly greater, then some of the points of the floor would be identified, so that the boundary would no longer be a simple curve. See the note page~\pageref{notroot3}.

The presence of this factor also suggests that the limiting spaces appearing in Theorem~\ref{cvqb} cannot easily be constructed from the Brownian map.

\bigskip

Except in Section~\ref{secsnake}, all the random variables considered in this work are taken on a common probability space $(\Omega, F, \Pb)$.

\paragraph{Acknowledgments}
The author wishes to heartily thank Gr\'egory Miermont for the precious advice and careful guidance he constantly offered during the realization of this work.

\section{The \BDG bijection}\label{secbdg}

As is often the case in this kind of problems, we start with a bijection allowing us to work with simpler objects. We use here a particular instance of the so-called \BDG bijection~\cite{bouttier04pml}, which has already been used in~\cite{bouttier09dsqb}. For more convenience, we modify it a little to better fit our purpose. This will allow us to code quadrangulations with a boundary by forests whose vertices carry integer labels.

\subsection{Forests}\label{secforb}

We use for forests the formalism of \cite{bettinelli10slr,bettinelli10tsl}, which we briefly recall here. We denote by $\N \de \{1,2,\dots\}$ the set of positive integers and for $i\le j$, $\ent i j \de [i,j] \cap \Z = \{ i, i+1, \dots, j \}$. For $u=(u_1,\dots,u_n)$, $v=(v_1,\dots,v_p) \in \bigcup_{n=1}^{\infty} \N^n$, we let $\lh u \rh \de n$ be the height of~$u$, and $uv \de (u_1,\dots,u_n,v_1,\dots,v_p)$ be the concatenation of~$u$ and~$v$. We say that~$u$ is an \textbf{ancestor} of~$uv$ and that~$uv$ is a \textbf{descendant} of~$u$. In the case where $v \in \N$, we use the terms \textbf{parent} and \textbf{child} instead.

\begin{defi}
A \textbf{forest} is a finite subset $\f \subset \bigcup_{n=1}^{\infty} \N^n$ satisfying:
\begin{enumerate}[($i$)]
 \item there is an integer $t(\f) \ge 1$, called the \textbf{number of trees} of~$\f$, such that $\f\cap \N = \ent 1 {t(\f)+1}$,
 \item if $u \in \f\bs\N$, then its parent belongs to $\f$,
 \item for every $u \in \f$, there is an integer $c_u(\f) \ge 0$ such that $ui \in \f$ if and only if $1 \le i \le c_u(\f)$,
 \item $c_{t(\f)+1}(\f) =0$.
\end{enumerate}
\end{defi}

\begin{figure}[ht]
	\begin{minipage}{0.55\linewidth}
		\psfrag{1}[][][.8]{$1$}
		\psfrag{2}[][][.8]{$2$}
		\psfrag{3}[][][.8]{$3$}
		\psfrag{4}[][][.8]{$4$}
		\psfrag{0}[][][.8]{$0$}
		\psfrag{9}[][][.8]{-$1$}
		\psfrag{8}[][][.8]{-$2$}
		\psfrag{a}[r][r][.8]{$\f(0)$, $\f(8)$}
		\psfrag{b}[r][r][.8]{$\f(1)$, $\f(5)$, $\f(7)$}
		\psfrag{c}[r][r][.8]{$\f(2)$, $\f(4)$}
		\psfrag{d}[r][r][.8]{$\f(3)$}
		\psfrag{e}[][][.8]{$\f(6)$}
		\psfrag{f}[][][.8]{$\f(9)$}
		\psfrag{g}[][][.8]{$\f(10)$}
	\centering\includegraphics[height=25mm]{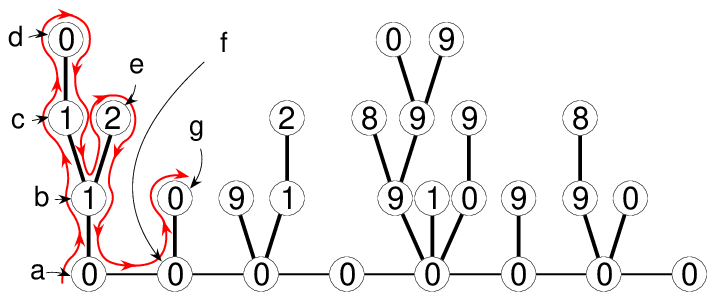}
	\end{minipage}
	\begin{minipage}{0.43\linewidth}
		\psfrag{C}[][][.8]{\textcolor{red}{$C_{\f}$}}
		\psfrag{L}[][][.8]{\textcolor{blue}{$L_{\f,\lab}$}}
	\centering\includegraphics[height=25mm]{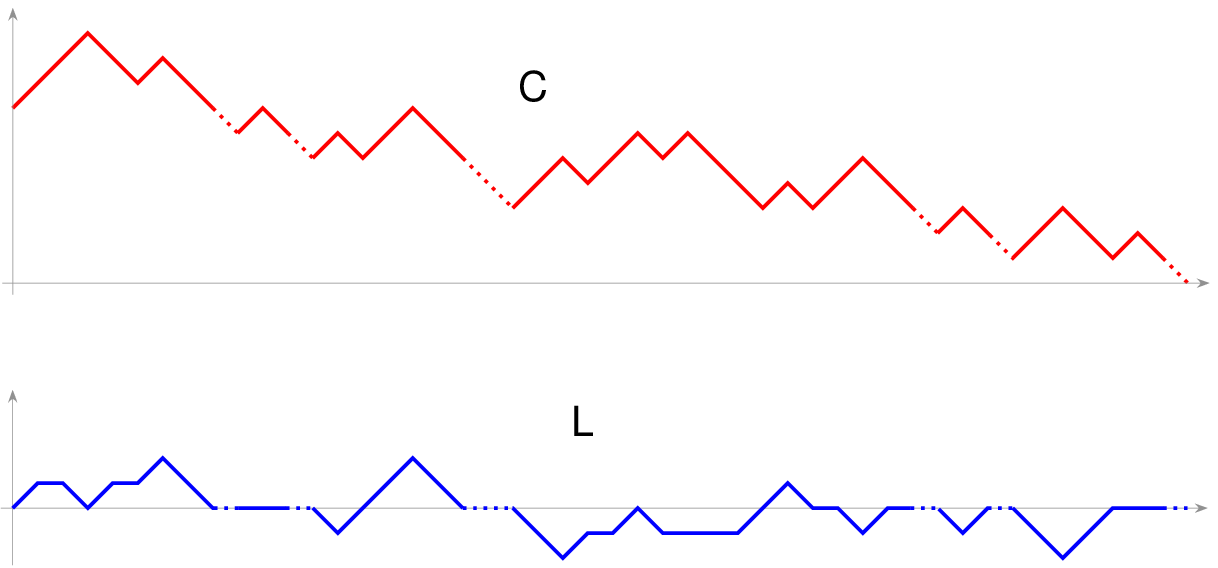}
	\end{minipage}
	\caption[Facial sequence and contour pair of a well-labeled forest.]{The facial sequence and contour pair of a well-labeled forest from $\fF_7^{20}$. The paths are dashed on the intervals corresponding to floor edges.}
	\label{forest}
\end{figure}

The set $\fl\de\f\cap \N$ is called the \textbf{floor} of the forest~$\f$. When $u\in \fl$, we sometime denote it by~$(u)$ to avoid confusion between the integer~$u$ and the point $(u)\in\f$. For $u=(u_1,\dots,u_p) \in \f$, we set $\aaa(u) \de u_1 \in \fl$ its oldest ancestor. For $1\le j \le t(\f)$, the set $\{u\in \f\,:\, \aaa(u)=j\}$ is called the tree of~$\f$ rooted at~$(j)$. Beware that the point $(t(\f)+1)\in \fl$ is not a tree. As we will see later, it is here for convenience. The points $u,v \in \f$ are called \textbf{neighbors}, and we write $u \sim v$, if either~$u$ is a parent or child of~$v$, or $u,v \in \fl$ and $|u-v| = 1$. On the figures, we always draw edges between neighbors (see Figure~\ref{forest}). We say that an edge drawn between a parent and its child is a \textbf{tree edge} whereas an edge drawn between two consecutive tree roots will be called a \textbf{floor edge}.

\begin{defi}
A \textbf{well-labeled forest} is a pair $(\f,\lab)$ where $\f$ is a forest and $\lab:\f \to \Z$ is a function satisfying:
\begin{enumerate}[($i$)]
 \item for all $u\in \fl$, $\lab(u) = 0$,
 \item if $u \sim v$, then $|\lab(u) - \lab(v)| \le 1$.
\end{enumerate}
\end{defi}

Let $\fF_\sigma^n \de \lb (\f,\lab):\, t(\f)=\sigma,\,|\f|= n+\sigma +1 \rb$ be the set of well-labeled forests with~$\sigma$ trees and~$n$ tree edges. By a simple application (see for example \cite[Lemma~3]{bettinelli10slr}) of the so-called cycle lemma \cite[Lemma~2]{bertoin03ptf}, and the fact that to every forest with~$n$ tree edges correspond exactly~$3^n$ labeling functions, we obtain that 
\begin{equation}\label{cardfns}
|\fF_\sigma^n|=3^n\,\frac \sigma {2n+\sigma}\choo{2n+\sigma}{n}.
\end{equation}

\bigskip

For a forest $\f$ with~$\sigma$ trees and~$n$ tree edges, we define its \textbf{facial sequence} $\f(0),\f(1),\dots,\f(2n+\sigma)$ as follows (see Figure~\ref{forest}): $\f(0) \de (1)$, and for $0\le i \le 2n+\sigma-1$,
\begin{itemize}
 \item if $\f(i)$ has children that do not appear in the sequence $\f(0),\f(1),\dots,\f(i)$, then $\f(i+1)$ is the first of these children, that is, $\f(i+1)\de \f(i)j_0$ where
 $$j_0 = \min\lb j\ge 1\,:\, \f(i)j \notin \lb \f(0),\f(1),\dots,\f(i) \rb \rb,$$
 \item otherwise, if $\f(i)\notin\fl$, then $\f(i+1)$ is the parent of $\f(i)$,
 \item if neither of these cases occur, which implies that $\f(i)\in\fl$, then $\f(i+1) \de \f(i)+1$. 
\end{itemize}
A well-labeled forest $(\f,\lab)$ is then entirely determined by its so-called \textbf{contour pair} $(C_\f,L_{\f,\lab})$ consisting in its \textbf{contour function} $C_\f:[0,2n+\sigma]\to \R_+$ and its \textbf{spatial contour function} $L_{\f,\lab}:[0,2n+\sigma]\to \R$ defined by
$$C_\f(i) \de \lh \f(i) \rh + t(\f) -\aaa\lp \f(i) \rp \sand L_{\f,\lab}(i) \de \lab(\f(i)),\qquad 0 \le i \le 2n+\sigma,$$
and linearly interpolated between integer values (see Figure~\ref{forest}).

\subsection{Bridges}

\begin{defi}
We say that a sequence of integers $(\bb(0), \bb(1), \dots, \bb(\sigma))$ for some $\sigma \ge 1$ is a \textbf{bridge} if $\bb(0)= 0$, $\bb(\sigma) \le 0$, and, for all $0\le i \le \sigma - 1$, we have $\bb(i+1)-\bb(i) \ge -1$. The integer~$\sigma$ will be called the \textbf{length} of the bridge.
\end{defi}

The somehow unusual condition $\bb(\sigma) \le 0$ will become clear in the following section: it will be used to keep track of the position of the root in the quadrangulation. We let~$\Bs$ be the set of all bridges of length~$\sigma$. In the following, when we consider a bridge $\bb\in \Bs$, we will always implicitly extend its definition to $[0,\sigma]$ by linear interpolation between integer values. 

\begin{lem}\label{cardbsi}
The cardinality of the set $\Bs$ is
$$\lt \Bs \rt = \choo{2\sigma}{\sigma}.$$
\end{lem}

\begin{pre}
With a bridge $(\bb(i))_{0\le i \le \sigma}\in \Bs$, we associate the following sequence
$$(\tilde\bb_j)_{1\le j \le 2\sigma} \de (\underbrace{+1,+1,\dots,+1}_{\bb(0)-\bb(\sigma) \textup{ times}},-1,\underbrace{+1,+1,\dots,+1}_{\bb(1)-\bb(0)+1 \textup{ times}},-1,\underbrace{+1,+1,\dots,+1}_{\bb(2)-\bb(1)+1 \textup{ times}},\dots, -1, \hspace*{-2mm}\underbrace{+1,+1,\dots,+1}_{\bb(\sigma)-\bb(\sigma-1)+1 \textup{ times}}\hspace*{-2mm}).$$
The set $\Bs$ is then in one-to-one correspondence with the set of sequences in $\{-1,+1\}^{2\sigma}$ counting exactly~$\sigma$ times the number $-1$. The number of bridges of length~$\sigma$ is then the number of choices we have to place these~$\sigma$ numbers among the $2\sigma$ spots.
\end{pre}

\subsection{The bijection}

A \textbf{pointed quadrangulation} (with a boundary) is a pair $(\q,v^\bullet)$ consisting in a quadrangulation (with a boundary) $\q$ together with a distinguished vertex $v^\bullet \in V(\q)$. We define
$$\Qnsb \de \lb (\q,v^\bullet):\, \q\in \Qns,\, v^\bullet \in V(\q) \rb$$
the set of all pointed quadrangulations with~$n$ internal faces and $2\sigma$ half-edges on the boundary. The \BDG bijection may easily be adapted into a bijection between the sets $\Qnsb$ and $\fF_\sigma^n\times \Bs$. We briefly describe it here, and refer the reader to \cite{bouttier04pml} for proofs and further details.

\subsubsection{From quadrangulations to forests and bridges}\label{bdgqfb}

\begin{figure}[ht]
		\psfrag{d}[][][.8]{$d$}
		\psfrag{e}[r][r][.8]{$d+1$}
		\psfrag{f}[l][l][.8]{$d+2$}
		\psfrag{g}[l][l][.8]{$d+1$}
		\psfrag{0}[][][.8]{$4$}
		\psfrag{1}[][][.8]{$5$}
		\psfrag{2}[][][.8]{$6$}
		\psfrag{3}[][][.8]{$7$}
		\psfrag{v}[][][.8]{$v^\circ$}
		\psfrag{a}[][][.6]{$v_0$}
		\psfrag{b}[l][l][.6]{$v_1=v_{i_1}$}
		\psfrag{c}[l][l][.6]{$v_2=v_{i_2}$}
		\psfrag{q}[][][.6]{$v_3$}
		\psfrag{s}[][][.6]{$v_{i_\sigma}$}
	\centering\includegraphics{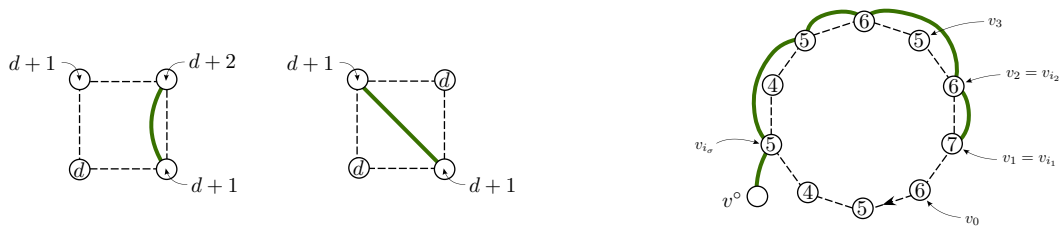}
	\caption[Adding new edges in the \BDG bijection.]{\textbf{Left.} Adding the new edge to an internal face. \textbf{Right.} Example of the operation on the external face. In this example, $\bb=(0,-1,-1,-2,-2,-1)$.}
	\label{faces}
\end{figure}

Let us start with the mapping from $\Qnsb$ onto $\fF_\sigma^n\times \Bs$. Let $(\q,v^\bullet)\in\Qnsb$. We label the vertices of~$\q$ as follows: for every vertex $v\in V(\q)$, we set $\tilde \lab(v)=d_\q(v^\bullet,v)$. Because~$\q$ is bipartite, the labels of both ends of any edge differ by exactly~$1$. As a result, the internal faces can be of two types: the labels around the face are either $d$, $d+1$, $d+2$, $d+1$, or $d$, $d+1$, $d$, $d+1$ for some~$d$. We add a new edge to every internal face as shown on the left part of Figure~\ref{faces}.

The operation regarding the external face is a little more intricate. We denote by~$v_0$, $v_1$, \dots, $v_{2\sigma-1}$ its vertices read in counterclockwise order\footnote{Recall that the external face is drawn as the unbounded face of the plane, so that the counterclockwise order on the plane is actually the clockwise order \emph{around the face}.}, starting at the origin of the root, $v_0=\e_*^-$ (and we use the convention $v_{2\sigma}=v_0$). We only consider the vertices $v_i$ such that $\tilde \lab(v_{i+1})=\tilde \lab(v_i)-1$. Note that, because $\tilde \lab(v_{i+1})-\tilde \lab(v_i) \in \{-1,+1\}$, there are exactly~$\sigma$ such vertices. We denote them by $v_{i_1}$, $v_{i_2}$, \dots, $v_{i_\sigma}$, with $0\le i_1< i_2<\dots<i_\sigma <2\sigma$. Finally, we add a new vertex $v^\circ$ inside the external face, and draw extra edges linking $v_{i_k}$ to $v_{i_{k+1}}$ for all $1\le k \le \sigma-1$, and $v_{i_\sigma}$ to $v^\circ$. See the right part of Figure~\ref{faces}.

We then only keep the new edges we added and the vertices in $(V(\q)\bs\{v^\bullet\})\cup\{v^\circ\}$. We obtain a forest~$\f$ whose floor is drawn in the external face: $(k)=v_{i_k}$ for $1\le k \le \sigma$, and $(\sigma+1)=v^\circ$. To obtain the labels of~$\f$, we shift the labels tree by tree, in such a way that the floor labels are~$0$: we define $\lab(u) \de \tilde \lab(u)-\tilde \lab(\aaa(u))$, and $\lab(v^\circ)=0$. Finally, the bridge~$\bb$ records the labels of the floor before the shifting operation: for $0\le k \le \sigma-1$, we let $\bb(k) \de \tilde \lab(v_{i_{k+1}}) -\tilde \lab(v_{i_1})$, and $\bb(\sigma)=\tilde \lab (v_0)-\tilde \lab(v_{i_1})$ (so that $\bb(\sigma)$ keeps track of the position of the root).

The pointed quadrangulation $(\q,v^\bullet)$ corresponds to the pair $((\f,\lab),\bb)$.

\begin{figure}[ht]
		\psfrag{v}[][][.8]{$v^\bullet$}
		\psfrag{w}[][][.8]{$v^\circ$}
		\psfrag{q}[][]{$\q$}
		\psfrag{f}[][]{$(\f,\lab)$}
		\psfrag{b}[][]{\textcolor{red}{$\bb$}}
		\psfrag{s}[][][.8]{$\sigma$}
		\psfrag{0}[][][.8]{$0$}
		\psfrag{1}[][][.8]{$1$}
		\psfrag{2}[][][.8]{$2$}
		\psfrag{3}[][][.8]{$3$}
		\psfrag{9}[][][.8]{-$1$}
		\psfrag{8}[][][.8]{-$2$}
	\centering\includegraphics{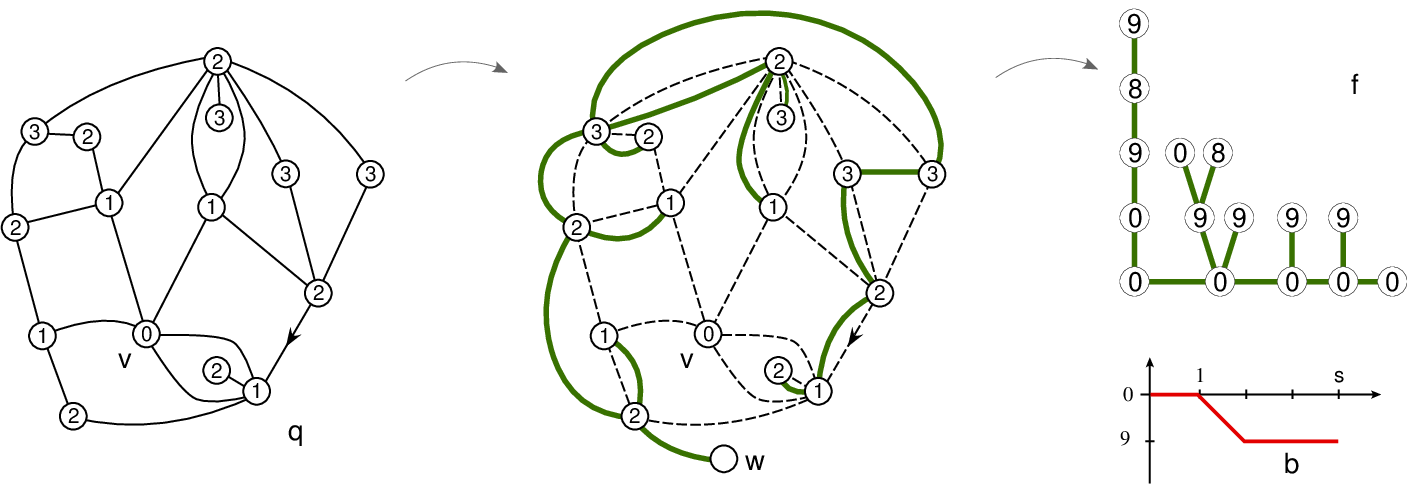}
	\caption[The \BDG bijection.]{The mapping from $\Qnsb$ onto $\fF_\sigma^n\times \Bs$. In this picture, $n= 10$ and $\sigma= 4$.}
\end{figure}

\subsubsection{From forests and bridges to quadrangulations}\label{bdgfbq}

Let us now describe the mapping from $\fF_\sigma^n\times \Bs$ onto $\Qnsb$. Let $(\f,\lab) \in \fF_\sigma^n$ be a well-labeled forest and $\bb \in \Bs$ be a bridge. As above, we write $\f(0)$, $\f(1)$, \dots, $\f(2n+\sigma)$ the facial sequence of~$\f$. The pointed quadrangulation $(\q,v^\bullet)$ corresponding to $((\f,\lab),\bb)$ is then constructed as follows. First, we shift all the labels of~$\f$ tree by tree according to the bridge $\bb$: precisely, we define $\hat \lab(u)\de\lab(u)+\bb(\aaa(u)-1)$. Then, we shift all the labels in such a way that the minimal label is equal to~$1$: let us set $\tilde \lab \de \hat\lab-\min\hat\lab +1$ as this shifted labeling function. We add an extra vertex~$v^\bullet$ carrying the label $\tilde \lab(v^\bullet)\de 0$ inside the only face of~$\f$. Finally, following the facial sequence, for every $0\le i \le 2n+\sigma-1$, we draw an arc---without intersecting any edge of~$\f$ or arc already drawn---between $\f(i)$ and $\f(\suc(i))$, where $\suc(i)$ is the \textbf{successor} of~$i$, defined by
\begin{equation}\label{succf}
\suc(i) \de \lb \begin{array}{cl}
		\!\inf S_\ge	&\text{ if } S_\ge \neq \varnothing\\
		\!\inf S_\le	&\text{ otherwise}
		\end{array}\rno\text{ with }
		\begin{array}{l}
		S_\ge\de\{k\in \ent i{2n\!+\!\sigma\!-\!1}\,:\, \tilde \lab(\f(k))= \tilde \lab(\f(i)) - 1\}\\
		S_\le \de \{k\in\ent0{i-1}\,:\, \tilde \lab(\f(k))= \tilde \lab(\f(i)) - 1\}
		\end{array}	
\end{equation}
with the conventions $\inf \varnothing = \infty$, and $\f(\infty) = v^\bullet$.

Because there may be more that one arc linking $\f(i)$ to $\f(\suc(i))$, we will speak of the arc linking~$i$ to $\suc(i)$ to avoid any confusion, and we will write it
$$i \arc \suc(i) \qquad \text{ or }\qquad \suc(i) \arc i.$$
When we need an orientation, we will write $i \arcr \suc(i)$ the arc oriented from~$i$ toward $\suc(i)$ and $i \arcl \suc(i)$ the arc oriented from $\suc(i)$ toward~$i$. The quadrangulation $\q$ is then defined as the map whose set of vertices is $(\f\bs\{(\sigma+1)\}) \cup \{v^\bullet\}$, whose edges are the arcs we drew, and whose root is either $\suc^{-\bb(\sigma)}(0) \arcr \suc^{-\bb(\sigma)+1}(0)$ if $\bb(\sigma) > \bb(\sigma-1) -1$, or $2n+\sigma - 1 \arcl \suc(2n+\sigma-1)$ if $\bb(\sigma) = \bb(\sigma-1) -1$.

\subsubsection{Some remarks}\label{secrem}

\noindent
\textbf{1.} Because of the way we drew the arcs of $\q$ in Section~\ref{bdgfbq}, it is easy to see that for any vertex $v\in V(\q)$, $\tilde \lab(v)=d_\q(v^\bullet,v)$, so that both functions $\tilde \lab$ of Sections~\ref{bdgqfb} and~\ref{bdgfbq} coincide.

\medskip
\noindent
\textbf{2.} Note that the sequence $\tilde \bb$ from the proof of Lemma~\ref{cardbsi} reads the increments of the labels around the boundary: $\tilde \bb_j=\tilde \lab(v_{j}) - \tilde \lab(v_{j-1})$ for $1 \le j \le 2\sigma$.

\medskip
\noindent
\textbf{3.} Using Lemma~\ref{cardbsi}, equation~\eqref{cardfns}, and the fact that every quadrangulation in $\Qns$ has exactly $n+\sigma +1$ vertices, we recover the following formula (see e.g.~\cite{bender94,schaeffer97,bouttier09dsqb} for other proofs)
$$|\Qns|=\frac{|\fF_\sigma^n|\ |\Bs|}{n+\sigma+1}=
	\frac{3^n (2\sigma)!\, (2n+\sigma -1)!}{\sigma! (\sigma-1)!\, n!\, (n+\sigma+1)!}.$$

\medskip
\noindent
\textbf{4.} If $(C,L)$ is the contour pair of $(\f,\lab)$, then we may retrieve the oldest ancestor of~$\f(i)$ thanks to~$C$ by the relation
$$\aaa\big(\f(i)\big) -1 = \sigma - \underline C(i),$$
where we use the notation
$$\underline{X}_s \de \inf_{[0, s]} X$$
for any process $(X_s)_{s\ge 0}$. The function
$$\Lab \de \Big( L(s) + \bb\big( \sigma - \underline C(s) \big) \Big)_{0\le s \le 2n +\sigma}$$
then records the labels of the forest, once shifted tree by tree according to the bridge~$\bb$. As a result, we see that $\Lab(i)-\min \Lab +1$ represents the distance in~$\q$ between~$v^\ooo$ and the point corresponding to~$\f(i)$.

\medskip
\noindent
\textbf{5.} This gives a natural way to explore the vertices of~$\q$: we denote by~$\q(i)$ the vertex corresponding to~$\f(i)$. In particular, $\{\q(i),\,0 \le i \le 2n+\sigma-1\} = V(\q)\bs\{v^\bullet\}$. We end this section by giving an upper bound for the distance between two vertices~$\q(i)$ and~$\q(j)$, in terms of the function $\Lab$:
\begin{equation}\label{dlemmeb}
d_\q(\q(i),\q(j)) \le \Lab(i) + \Lab(j) - 2 \max \lp \min_{ k \in \overrightarrow{\ent i j}} \Lab(k),\min_{ k \in \overrightarrow{\ent j i}} \Lab(k) \rp +2
\end{equation}
where we define
\begin{equation}\label{oraijb}
\overrightarrow{\ent i j} \de \lb \begin{array}{cll}
					\ent i j			&\text{ if } & i \le j,\\
					\ent i {2n+\sigma-1} \cup \ent 0 j	&\text{ if } & j < i.
				    \end{array}\rno
\end{equation}
This kind of bounds is often used in these problems (see e.g.~\cite{legall07tss,miermont09trm,bettinelli10tsl}). We refer the reader to \cite[Lemma~3]{miermont09trm} for a detailed proof.

\section{Proof of Theorem~\ref{cvqb}}\label{secpt1}

\subsection{Convergence of the coding functions}\label{secccf}

Let $(\sigma_n)_{n\ge 1}$ be a sequence of positive integers such that
$$\sigma_{(n)} \de \frac{\sigma_n}{\sqrt{2n}}\ton \sigma \in [0,\8].$$
Until further notice, we suppose that $\sigma \in (0,\8)$. The remaining cases $\sigma=0$ and $\sigma=\8$ will be treated separately in Section~\ref{singcases}. Let $\q_n$ be uniformly distributed over the set $\Qnsn$ of quadrangulation with $n$ internal faces and $2\sigma_n$ half-edges on the boundary. Conditionally given $\q_n$, we let $v_n^\bullet$ be uniformly distributed over the set $V(\q_n)$. Because every quadrangulation in $\Qnsn$ has exactly $n+\sigma_n +1$ vertices (by Euler characteristic formula), we see that $(\q_n,v_n^\bullet)$ is uniformly distributed over $\Qnsnb$, and thus corresponds through the \BDG bijection to a pair $\big((\f_n,\lab_n),\bb_n\big)$ uniformly distributed over the set $\fF_{\sigma_n}^n\times \Bsn$.

\subsubsection{Brownian bridges, first-passage Brownian bridges, and Brownian snake}\label{secbbfpbbs}

Let us define the space
$$\K \de \bigcup_{x \in \R_+} \C([0, x],\R)$$
of continuous real-valued functions on $\R_+$ killed after some time. For an element $f \in \K$, let~$\zeta(f)$ denote its lifetime, that is, the only~$x$ such that $f \in \C([0, x],\R)$. We endow this space with the following metric:
$$d_\K (f,g) \de |\zeta(f) - \zeta(g) | + \sup_{y \ge 0} \lt f\big(y\wedge \zeta(f)\big)-g\big(y\wedge\zeta(g)\big)\rt.$$

We write $B_{[0,\sigma]}^{0\to 0}$ a Brownian bridge on $[0,\sigma]$ from~$0$ to~$0$, defined as a standard Brownian motion on $[0,\sigma]$ started at~$0$, conditioned on being at~$0$ at time~$\sigma$ (see for example \cite{bertoin03ptf, bettinelli10slr, billingsley68cpm, revuz99cma}). We also denote by $F_{[0,1]}^{\sigma \to 0}$ a first-passage Brownian bridge on $[0,1]$ from~$\sigma$ to~$0$, defined as a standard Brownian motion on $[0,1]$ started at~$\sigma$, and conditioned on hitting~$0$ for the first time at time~$1$. We refer the reader to~\cite{bettinelli10slr} for a proper definition of this conditioning, as well as for some convergence results of the discrete analogs.

The so-called Brownian snake's head may then be defined as the process $\big( F_{[0,1]}^{\sigma \to 0}, Z_{[0,1]}\big)$, where, conditionally given $F_{[0,1]}^{\sigma \to 0}$, the process $\big(Z_{[0,1]}(s)\big)_{0\le s \le 1}$ is a centered Gaussian process with covariance function
\begin{equation}\label{defsnakeshead}
\cov\big(Z_{[0,1]}(s), Z_{[0,1]}(s')\big) = \inf_{[s\wedge s',s\vee s']} \big(F_{[0,1]}^{\sigma \to 0} - \underline F_{[0,1]}^{\sigma \to 0}\big).
\end{equation}
We refer to \cite{bettinelli10slr, duquesne02rtl, legall99sbp} for more details.

\subsubsection{Convergence of the bridge and the contour pair of the well-labeled forest}

We let $(C_n, L_n)$ be the contour pair of $(\f_n,\lab_n)$, and we define the scaled versions of~$C_n$, $L_n$, and~$\bb_n$ by
\begin{align*}
C_{(n)} &\de \lp \frac {C_n\big((2n +\sigma_n-1)s\big)} {\sqrt{2n}} \rp_{0\le s \le 1}&
L_{(n)} &\de \lp \frac {L_n\big((2n +\sigma_n-1)s\big)} {\g}\rp_{0\le s \le 1}
\end{align*}
$$\bb_{(n)} \de \lp \frac {\bb_n(\sqrt{2n}\,s)} {\g}\rp_{0\le s \le \sigma_{(n)}}$$
where the constant~$\gamma$ was defined during the statement of Theorem~\ref{cvqb}.

\begin{rem}
Following \cite{bettinelli10slr,bettinelli10tsl}, the notation with a parenthesized $n$ will always refer to suitably rescaled objects, as in the definitions above.
\end{rem}

The aim of this section is the following proposition.

\begin{prop}\label{cvcont}
The triple $(C_{(n)},L_{(n)},\bb_{(n)})$ converges in distribution in the space $(\K,d_\K)^3$ toward a triple $(C_\infty,L_\infty,\bb_\infty)$ whose law is defined as follows:
\begin{itemize}
	\item the processes $(C_\infty,L_\infty)$ and $\bb_\infty$ are independent,
	\item the process $(C_\infty,L_\infty)$ has the law of a Brownian snake's head on $[0,
1]$ going from $\sigma$ to $0$:
	$$(C_\infty,L_\infty) \law \lp F^{\sigma \to 0}_{[0,1]}, Z_{[0,1]}\rp,$$
	\item the process $\bb_\infty$ has the law of a Brownian bridge on $[0,\sigma]$ from~$0$ to~$0$, scaled by the factor~$\sqrt 3$:
	$$\bb_\infty \law \sqrt 3 \, B_{[0,\sigma]}^{0\to 0}.$$
\end{itemize}
\end{prop}

\begin{pre}
By \cite[Corollary 16]{bettinelli10slr}, the pair $(C_{(n)}, L_{(n)})$ converges in distribution\footnote{In \cite{bettinelli10slr}, the processes considered were the same except that the term $(2n +\sigma_n-1)$ was replaced with $2n$. The fact that $\sigma_n/2n \to 0$ and the uniform continuity of the process $\big( F^{\sigma \to 0}_{[0,1]}, Z_{[0,1]}\big)$ yield the result as stated here.} toward the pair $\big( F^{\sigma \to 0}_{[0,1]}, Z_{[0,1]}\big)$, in the space $\lp \K, \, d_\K \rp^2$. Moreover, $(C_n,L_n)$ and $\bb_n$ are independent, so that it only remains to show that $\bb_{(n)}$ converges in distribution toward $\sqrt 3 \, B_{[0,\sigma]}^{0\to 0}$. To that end, we will use \cite[Lemma~10]{bettinelli10slr}.

Let $(X_i)_{i \ge 1}$ be a sequence of i.i.d. random variables with distribution given by 
$$\Pb(X_i=p)=2^{-p-2},\qquad p\ge -1.$$
We set $\Sigma_0\de 0$ and, for $j \ge 1$, $\Sigma_j \de \sum_{i=1}^j X_i$. For $k\ge 0$ fixed, and $n$ such that $\sigma_n \ge k$, we also define a process $(S_n^k(i))_{0\le i \le \sigma_n}$ distributed as $(\Sigma_i)_{0\le i \le \sigma_n}$ conditioned on the event $\{\Sigma_{\sigma_n}=-k\}$. We extend its definition to $[0,\sigma_n]$ by linear interpolation between integer values. Because~$X_1$ is centered, has moments of any order, and has variance~$2$, we may apply \cite[Lemma~10]{bettinelli10slr} and we see that the process
\begin{equation}\label{cvbb}
\lp \frac{S_n^k(\sqrt{2n}\,s)}{\g} \rp_{0\le s \le \sigma_{(n)}} \tol \sqrt 3\, B_{[0,\sigma]}^{0\to 0}.
\end{equation}

Moreover, it is easy to see that the bridge $S_n^k$ is uniform over the set $\{ \bb \in \Bsn \, :\, \bb(\sigma_n)=-k\}$. Indeed, for any $\bb \in \Bsn$ such that $\bb(\sigma_n)=-k$, we have
$$\Pb(S_n^k=\bb)=\frac{\Pb\big(\forall i \in \ent 1 {\sigma_n},\ X_i=\bb(i)-\bb(i-1)\big)}{\Pb(\Sigma_{\sigma_n}=-k)}=\frac{2^{-2\sigma_n+k}}{\Pb(\Sigma_{\sigma_n}=-k)},$$
which does not depend on~$\bb$ but only on~$n$ and~$k$. For such a~$\bb$, we set
$$c_{n,k} \de \frac{\Pb(\bb_n=\bb)}{\Pb(S_n^k=\bb)} = \choo {2\sigma_n}{\sigma_n}^{-1} \choo{2\sigma_n - k-1}{\sigma_n-1}.$$
(We may use the bijection of Lemma~\ref{cardbsi} to compute the denominator.) We have that
$$c_{n,k} = \frac{1}{2} \frac{(2\sigma_n-k-1)!}{(2\sigma_n-1)!} \frac{\sigma_n!}{(\sigma_n-k)!}
    \le \frac{1}{2} \prod_{i=0}^{k-1} \frac{\sigma_n-i}{\sigma_n-i+\sigma_n-1}
    \le 2^{-k} ,$$
and that $c_{n,k} \to 2^{-k-1}$ as $n \to \infty$. Now, let $\varphi : \K \to \R$ be a bounded measurable function. Using~\eqref{cvbb} and the fact that $c_{n,k}=\Pb(\bb_n(\sigma_n)=-k)$, we obtain by dominated convergence that
\begin{align*}
\E{\varphi(\bb_{(n)})} &= \sum_{k=0}^{\sigma_n} c_{n,k} \ \E{\varphi \lp \lp \frac {S_n^k(\sqrt{2n}\,s)} {\g}\rp_{0\le s \le \sigma_{(n)}}\rp} \ton \E{\varphi\big(\sqrt 3 \, B_{[0,\sigma]}^{0\to 0}\big)}.
\end{align*}
This completes the proof.
\end{pre}

\bigskip

Recall the notation $\q_n(i)$ introduced at the end of Section~\ref{secbdg} for the vertex corresponding to $\f_n(i)$ through the \BDG bijection. Remember that $d_{\q_n}(v_n^\bullet, \q_n(i))=\Lab_n(i)-\min \Lab_n +1$, where
\begin{equation}\label{Labn}
\Lab_n \de \Big( L_n(s) + \bb_n\big( \sigma_n - \underline C_n(s) \big) \Big)_{0\le s \le 2n +\sigma_n}.
\end{equation}
The rescaled version of $\Lab_n$ is then given by
$$\Lab_{(n)} \de \lp \frac {\Lab_n\big((2n +\sigma_n-1)s\big)} {\g}\rp_{0\le s \le 1}=
\Big( L_{(n)}(s) + \bb_{(n)}\big( \sigma_{(n)} - \underline C_{(n)}(s) \big) \Big)_{0\le s \le 1}.$$

\begin{corol}\label{cvln}
The process $(C_{(n)},\Lab_{(n)})$ converges in distribution in the space $(\K,d_\K)^2$ toward the process $(C_\infty,\Lab_\infty)$, where
\begin{equation}\label{defli}
\Lab_\infty \de \Big( L_\infty(s) + \bb_\infty\big( \sigma - \underline C_\infty(s) \big) \Big)_{0\le s \le 1}.
\end{equation}
\end{corol}

\subsection{Proof of Theorem~\ref{cvqb}}\label{secprob}

The proof of  Theorem~\ref{cvqb} is very similar to~\cite[Section~6]{bettinelli10slr}, so that we only sketch it. Our approach is adapted from Le~Gall \cite{legall07tss} for the first assertion, and from Le~Gall and Miermont \cite{legall09scaling} for the Hausdorff dimension. In addition, we use this occasion to introduce some notation that will be useful later.

We define on $\ent 0 {2n+\sigma_n-1}$ the pseudo-metric~$d_n$ by
$$d_n(i,j) \de d_{\q_n}\lp \q_n(i),\q_n(j)\rp,$$
we extend its definition to non integer values by linear interpolation: for $s,t$ in $[0,2n + \sigma_n -1]$,
\begin{equation*}
d_n(s,t)\de \ul s\, \ul t\, d_n(\lf s\rf+1,\lf t\rf+1) +
	      \ul s\, \ol t\, d_n(\lf s\rf+1,\lf t\rf) +
	      \ol s\, \ul t\, d_n(\lf s\rf,\lf t\rf+1) +
	      \ol s\, \ol t\, d_n(\lf s\rf,\lf t\rf),
\end{equation*}
where $\lf s\rf \de \sup\{k \in \Z,\, k\le s\}$, $\ul s \de s - \lf s\rf$ and $\ol s \de \lf s\rf+1 - s$, and we define its rescaled version: for $s,t \in [0,1]$, we let
\begin{equation*}
d_{(n)}(s,t) \de \frac 1 {\g}\, d_n\big((2n+\sigma_n-1)s,(2n+\sigma_n-1)t\big).
\end{equation*}
We also define the equivalence relation~$\sim_n$ on $\ent 0 {2n+\sigma_n-1}$ by declaring that $i \sim_n j$ when $\q_n(i)=\q_n(j)$, which is equivalent to $d_n(i,j) =0$. The function~$d_{(n)}$ may then be seen as a metric on
\begin{equation*}
\QQ_n \de (2n+\sigma_n-1)^{-1} \, \ent 0 {2n+\sigma_n-1}_{/\sim_n},
\end{equation*}
and, as~$v_n^\ooo$ is the only point of~$\q_n$ that does not lie in $\{\q_n(i)\, : \, 0 \le i \le {2n+\sigma_n-1}\}$, we have
\begin{equation}\label{ghb}
\dGH \lp \bigg( \QQ_n ,d_{(n)} \bigg), \bigg( V(\q_n),\frac 1 {\g}\, d_{\q_n} \bigg) \rp \le \frac 1 {\g}.
\end{equation}

The bound~\eqref{dlemmeb} gives us a control on the metric~$d_{(n)}$, from which we can derive the following lemma (see \cite[Lemma~19]{bettinelli10slr}).

\begin{lem}\label{tightdb}
The distributions of the quadruples of processes
$$\lp C_{(n)},L_{(n)},\bb_{(n)}, \lp d_{(n)}(s,t) \rp_{0 \le s,t \le 1}\rp, \qquad n\ge 1$$
form a relatively compact family of probability distributions.
\end{lem}

As a result of Lemma~\ref{tightdb}, from any increasing sequence of integers, we may extract a (deterministic) subsequence $(n_k)_{k\ge 0}$ such that there exists a random function $\disig \in \C([0,1]^2,\R)$ satisfying
\begin{equation}\label{dinftyb}
\lp d_{(n_k)}(s,t) \rp_{0 \le s,t \le 1} \tode{(d)}{k}{\infty} \lp \disig(s,t) \rp_{0 \le s,t \le 1},
\end{equation}
and such that this convergence holds jointly with the convergence of Proposition~\ref{cvcont} and Corollary~\ref{cvln}. From now on, we fix such a subsequence $(n_k)_{k\ge 0}$. We will generally focus on this particular subsequence in the following, and we will often assume convergences when $n\to \8$ to hold along this particular subsequence. By Skorokhod's representation theorem, we may and will moreover assume that this joint convergence holds almost surely. In the limit, the bound~\eqref{dlemmeb} becomes
\begin{equation}\label{boundsig0}
\disig(s,t) \le d^\circ_\infty(s,t) \de \Lab_\infty(s) + \Lab_\infty(t) - 2 \max \lp \min_{ x \in \overrightarrow{[s,t]}} \Lab_\infty(x),\min_{ x \in \overrightarrow{[t,s]}} \Lab_\infty(x) \rp, \quad 0 \le s,t \le 1,
\end{equation}
where
\begin{equation}\label{orastb}
\overrightarrow{[s,t]} \de \lb \begin{array}{cll}
					[s,t]			&\text{ if } & s \le t, \\
					\left[s,1\right] \cup [0,t]	&\text{ if } & t < s.
				    \end{array}\rno
\end{equation}

\begin{nota}
Beware not to confuse $d^\circ_\infty$ with $d^0_\infty$. In fact, we will never use the latter symbol so there should not be any confusion. 
\end{nota}

Adding to this the fact that the functions~$d_{(n)}$ obey the triangle inequality, we see that the function~$\disig$ is a pseudo-metric. We define the equivalence relation associated with it by saying that $s \sim_\infty t$ if $\disig(s,t)=0$, and we set $\qis \de [0,1]_{/\sim_\infty}$. The convergence claimed in Theorem~\ref{cvqb} holds along the same subsequence $(n_k)_{k\ge 0}$.

To see this, we use the characterization of the Gromov--Hausdorff distance via correspondences. Recall that a correspondence between two metric spaces $(\X,\delta)$ and $(\X',\delta')$ is a subset $\rR\subseteq \X \times \X'$ such that for all $x\in \X$, there is at least one $x'\in \X'$ for which $(x,x')\in \rR$ and vice versa. The distortion of the correspondence $\rR$ is defined by
$$\dis(\rR) \de \sup \lb |\delta(x,y) - \delta(x',y')|\,:\, (x,x'),(y,y')\in \rR \rb.$$
Then we have \cite[Theorem~7.3.25]{burago01cmg}
\begin{equation}\label{dghcorres}
\dGH(\X,\X') = \frac 12 \inf_{\rR} \dis(\rR),
\end{equation}
where the infimum is taken over all correspondences between $\X$ and $\X'$.

We denote by~$\bmp_n$ the canonical projection from $\ent 0 {2n+\sigma_n-1}$ to $\ent 0 {2n+\sigma_n-1}_{/\sim_n}$. For $t\in [0,1]$, we define $\bmp_{(n)}(t) \de (2n+\sigma_n-1)^{-1}\, \bmp_{n}(\lf (2n+\sigma_n-1) t \rf)$, and we denote by~$\qis(t)$ the equivalence class of~$t$ in~$\qis$. We then define the correspondence~$\rR_n$ between the spaces $(\QQ_n ,d_{(n)} )$ and $\lp \qis,\disig \rp$ as the set
\begin{equation*}
\rR_n \de \lb \lp \bmp_{(n)}(t),\qis(t)\rp, \, t\in [0,1] \rb.
\end{equation*}
Its distortion is
$$\dis(\rR_n) = \sup_{0 \le s,t \le 1} \Big| d_{(n)}\lp \frac{\lf (2n+\sigma_n-1) s \rf}{2n+\sigma_n-1},\frac{\lf (2n+\sigma_n-1) t \rf}{2n+\sigma_n-1}\rp - \disig(s,t) \Big|,$$
and, thanks to~\eqref{dinftyb},
$$\dGH \lp \lp \QQ_{n_k},d_{(n_k)} \rp, \lp \qis,\disig \rp \rp \le \frac 12\, \dis \lp \rR_{n_k} \rp \tok 0.$$
Combining this with~\eqref{ghb}, we obtain the first assertion of Theorem~\ref{cvqb}.

\bigskip

The Hausdorff dimension of the limit may be computed by the technique we used in \cite{bettinelli10slr}. Because the proof is very similar, and is not really related to our purpose here, we leave it to the reader. The idea is roughly the following. To prove that the Hausdorff dimension is less than $4$, we use the fact that $\Lab_\infty$ is almost surely $\alpha$-H\"older for all $\alpha \in (0 , 1/4)$, yielding that the canonical projection from $([0,1],|\cdot|)$ to $(\qis,\disig)$ is also $\alpha$-H\"older for the same values of~$\alpha$. To prove that it is greater than~$4$, we show that the size of the balls of diameter~$\delta$ is of order~$\delta^4$. To see this, we first bound from below the distances in terms of label variation along the branches of the forest, and then use twice the law of the iterated logarithm: this tells us that, for a fixed $s\in [0,1]$, the points outside of the set $[s-\delta^4,s+\delta^4]$ code points that are at distance at least~$\delta^2$ from~$\qis(s)$ in the forest, so that their distance from~$\qis(s)$ is at least~$\delta$ in the map.
See \cite[Section~6.4]{bettinelli10slr} for a complete proof. We will also use a similar approach to show Theorem~\ref{thmdimh} in Section~\ref{secboundary}.

\section{Maps seen as quotients of real forests}\label{secrf}

In the discrete setting, the metric space $(V(\q_n),d_{\q_n})$ may either be seen as a quotient of $\ent 0 {2n+\sigma_n-1}$, as in last section, or directly as the space~$\f_n$ endowed with the proper metric. In the continuous setting, we defined~$\qis$ as a quotient of $[0,1]$, but it will also be useful to see it as a quotient of a continuous analog to~$\f_n$. We obtain a quotient, because some points may be very close in the discrete forest, and become identified in the limit. Finding a criterion telling which points are identified in the limit will be the object of Section~\ref{secipb}. In a first time, we define the continuous analog to forests.

\subsection{Real forests}

\begin{figure}[ht]
		\psfrag{a}[][][.8]{$\alpha$}
		\psfrag{b}[][][.8]{$\beta$}
		\psfrag{A}[][][.8]{\textcolor{blue}{$\lhb \alpha, \beta \rhb$}}
		\psfrag{B}[][][.8]{\textcolor{vert}{$\lhb b, c \rhb$}}
		\psfrag{C}[][]{$[\alpha,\beta]$}
		\psfrag{c}[][][.8]{$a$}
		\psfrag{d}[][][.8]{$b$}
		\psfrag{g}[][][.8]{$c$}
		\psfrag{l}[][][.8]{$l$}
		\psfrag{r}[][][.8]{$r$}
		\psfrag{s}[][][.8]{$l'$}
		\psfrag{t}[][][.8]{$r'$}
		\psfrag{p}[][][.8]{$\rho_a$}
		\psfrag{z}[][][.8]{$\rho_b$}
		\psfrag{q}[][][.8]{$\rho$}
		\psfrag{w}[][][.8]{$\ro$}
		\psfrag{u}[][]{$\tau$}
		\psfrag{v}[][]{$\tau_a$}
		\psfrag{y}[][]{\textcolor{red}{$\fl$}}
	\centering\includegraphics{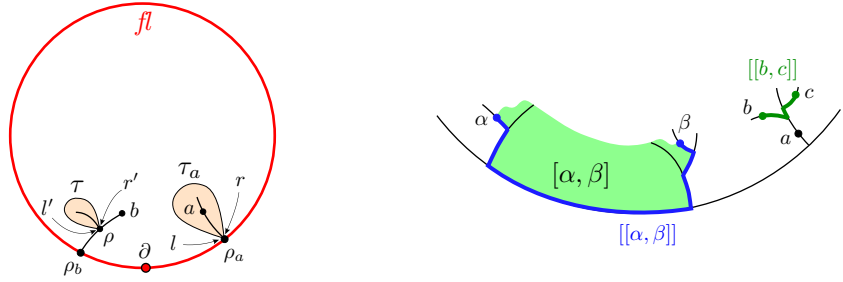}
	\caption[Notation for real forests.]{\textbf{Left.} On this picture, we can see the root~$\ro$, the floor~$\fl$, an example of tree~$\tau_a$ (coded by $[l,r]$), and an example of tree~$\tau$ to the left of $\lhb \rho_b,b \rhb$ rooted at~$\rho$ (coded by $[l',r']$). \textbf{Right.} On this picture, $a$ is an ancestor of~$b$ and $c$, and we can see the sets $\lhb b,c \rhb$, $[\alpha, \beta]$, and $\lhb \alpha,\beta \rhb$.}
	\label{notqb}
\end{figure}

We define here real forests in a way convenient to our purpose, by adapting the notions used in~\cite[Section~3]{bettinelli10tsl}. We will also need basic facts on real trees (see for example~\cite{legall05rta}). We consider a continuous function $h:[0,1] \to \R_+$ such that $h(1)=0$, and we define on $[0,1]$ the relation~$\simeq$ as the coarsest equivalence relation such that $0 \simeq 1$, and $s \simeq t$ if 
\begin{equation}\label{seqt}
h(s)=h(t)=\inf_{[s\wedge t, s\vee t]} h.
\end{equation}
In other words, the second relation identifies the points ``facing each other under the graph of~$h$.'' We call \textbf{real forest} any set $\FF \de [0,1]_{/\simeq}$ obtained by such a construction. It is possible to endow it with a natural metric, but we will not use it in this work. We now define the notions we will use throughout this work (see Figure~\ref{notqb}). For $s \in [0,1]$, we write $\FF(s)$ its equivalence class in the quotient $\FF = [0,1]_{/\simeq}$. In a way, we see $(\FF(s))_{0\le s \le 1}$ as the continuous facial sequence of~$\FF$. We call \textbf{root} of $\FF$ the point $\ro\de \FF(0)=\FF(1)$.

\begin{defi}
The \textbf{floor} of $\FF$ is the set $\fl \de \FF\lp \lb s \, :\, h(s) = \underline h(s) \rb \rp$.
\end{defi}

For $a=\FF(s)\in\FF\bs\fl$, let $l \de \inf\{t \le s \, :\, \underline h(t)=\underline{h}(s) \}$ and $r \de \sup\{t \ge s \, :\, \underline h(t)=\underline{h}(s) \}$. Note that, once endowed with the natural metric, the set $\tau_a \de \FF([l,r])$ is a real tree rooted at $\rho_a \de \FF(l)=\FF(r)\in\fl$. In the following, we will not need metric properties about real trees, we will only see them as topological spaces.

\begin{defi}
We call \textbf{tree} of~$\FF$ a set of the form $\tau_a$ for any $a\in\FF\bs\fl$.
\end{defi}

If $a\in\fl$, we simply set $\rho_a \de a$. Let~$\tau$ be a tree of~$\FF$ rooted at~$\rho$, and $a,b \in \tau$. We let $\lhb a,b \rhb$ bet the range of the unique injective path linking~$a$ to~$b$. In particular, the set $\lhb \rho,a\rhb$ represents the ancestral lineage of~$a$ in the tree $\tau$. We say that~$a$ is an \textbf{ancestor} of~$b$, and we write $a \preceq b$, if $a \in \lhb \rho,b \rhb$. We write $a \prec b$ if $a \preceq b$ and $a \neq b$.

Let $a,b\in \FF$ be two points. There is a natural way to explore the forest $\FF$ from~$a$ to~$b$. If $\inf \FF^{-1}(a) \le \sup \FF^{-1}(b)$, then let $t\de \inf \{ r \ge \inf \FF^{-1}(a)\, :\, b=\FF(r) \}$ and $s\de \sup \{ r \le t\, :\, a=\FF(r) \}$. If $\sup \FF^{-1}(b) < \inf \FF^{-1}(a)$, then let $t \de \inf \FF^{-1} (b)$ and $s\de \sup \FF^{-1}(a)$. We define
\begin{equation}\label{abb}
{[a,b]} \de \FF \lp \overrightarrow{[s,t]} \rp,\vspace{-2mm}
\end{equation}
where $\overrightarrow{[s,t]}$ is defined by~\eqref{orastb}. We may now extend the definition of $\lhb a,b \rhb$ to any two points in $\FF$. First, for $a,b \in \fl$, we let $\lhb a,b \rhb \de [a,b]\cap \fl$. Then, for any points $a$, $b\in\FF$ such that $\rho_a\neq \rho_b$, we define
$$\lhb a, b \rhb \de \lhb a, \rho_a \rhb \cup \lhb \rho_a, \rho_b \rhb \cup \lhb \rho_b, b \rhb,$$
so that it is the range of the unique injective path from~$a$ to~$b$ that stays inside $[a,b]$.

\begin{defi}
Let $b=\FF(t) \in \FF \bs \fl$ and $\rho \in \lhb \rho_b,b\rhb \bs \{\rho_b,b\}$. Let $l' \de \inf\{s\le t \, :\, \FF(s)=\rho\}$ and $r' \de \sup\{s\le t \, :\, \FF(s)=\rho\}$. Then, provided $l' \neq r'$, we call \textbf{tree to the left} of $\lhb \rho_b,b\rhb$ rooted at $\rho$ the set $\FF([l',r'])$.

We define the \textbf{tree to the right} of $\lhb \rho_b,b\rhb$ rooted at $\rho$ in a similar way, by replacing ``$\le$'' with ``$\ge$'' in the definitions of $l'$ and $r'$. 
\end{defi}

\begin{defi}
We call \textbf{subtree} of $\FF$ any tree of $\FF$, or any tree to the left or right of $\lhb \rho_b,b\rhb$ for some $b \in \FF \bs \fl$.
\end{defi}

Note that subtrees of $\FF$ are real trees, and that trees of~$\FF$ are also subtrees of~$\FF$. The maximal interval $[s,t]$ such that $\tau = \FF([s,t])$ is called the \textbf{interval coding} the subtree~$\tau$.

\bigskip

We denote by~$\FF_n$ the real forest obtained from the function $s\in[0,1] \mapsto C_n\big((2n +\sigma_n)s\big)$, as well as~$\Fi$ the real forest obtained from the function $C_\infty$. We also denote by $\simeq_{(n)}$ and $\simeq_\infty$ the corresponding equivalence relations. We write $\ro_\infty$ the root of $\Fi$, and $\fli$ its floor. It is more natural to use~$\f_n$ rather than $\FF_n$ in the discrete setting. As~$\f_n$ may be viewed as a subset of $\FF_n$ (when identifying $(\sigma_n +1)$ with $(1)$), we will use for~$\f_n$ the formalism we defined above simply by restriction. Note that the notions of floor and trees are consistent with the definitions we gave in Section~\ref{secforb} in this case.

Remark that, because the function $C_\infty$ is a first-passage Brownian bridge, there are almost surely no trees rooted at the root $\ro_\infty$ of $\Fi$, and all the points of~$\Fi$ are of order less than~$3$, in the sense that for all $a\in \Fi$ and every connected subset $C \subseteq \Fi$, the number of connected components of $C\bs\{a\}$ is less than~$3$. We will not use this remark in the following, so that we do not go into further details.

\subsection{Quotient of real forests}

Similarly to the notation $\f_n(i)$ and $\q_n(i)$ in the discrete setting, we denote by $\Fi(s)$ (resp.~$\qis(s)$) the equivalence class of $s\in [0,1]$ in $\Fi = [0,1]_{/\simeq_\infty}$ (resp.\ in $\qis = [0,1]_{/\sim_\infty}$).

\begin{lem}\label{coarseb}
The equivalence relation $\simeq_\infty$ is coarser than~$\sim_\infty$.
\end{lem}

\begin{pre}
First, notice that, by~\eqref{boundsig0}, we have $\disig(0,1) \le d^\circ_\infty(0,1)=0$, so that $0 \sim_\infty 1$. The remaining is then identical to the first part of the proof of~\cite[Lemma~6]{bettinelli10tsl}.
\end{pre}

This allows us to define a pseudo-metric and an equivalence relation on~$\Fi$, still denoted by~$\disig$ and~$\sim_\infty$, by setting $\disig\big(\Fi(s),\Fi(t)\big) \de \disig(s,t)$ and declaring $\Fi(s) \sim_\infty \Fi(t)$ if $s \sim_\infty t$. The metric space $(\qis,\disig)$ is thus isometric to $\big({\Fi}_{/\sim_\infty},\disig\big)$. We also define $d^\circ_\infty$ on~$\Fi$ by letting
$$d^\circ_\infty(a,b) \de \inf \lb d^\circ_\infty(s,t) \,:\, a=\Fi(s),\, b=\Fi(t) \rb.$$
We will see in Lemma~\ref{minb} that there is a.s.\ only one point where the function~$\Lab_\infty$ reaches its minimum. If  $s^\bullet\in[0,1]$ denotes this point, then it is not hard (see \cite[Lemma~7]{bettinelli10tsl}) to see from the fourth remark of Section~\ref{secrem} that
$$\disig(s,s^\bullet) = \Lab_\infty(s) - \Lab_\infty(s^\bullet).$$
By the triangle inequality, we obtain that $s\sim_\infty t$ implies $\Lab_\infty(s)=\Lab_\infty(t)$, so that, in particular, $s\simeq_\infty t$ implies $\Lab_\infty(s)=\Lab_\infty(t)$, by Lemma~\ref{coarseb}. It is then licit to see~$\Lab_\infty$ as a function on~$\Fi$ by setting $\Lab_\infty\big(\Fi(s)\big) \de \Lab_\infty(s)$. This yields a more explicit expression for~$d^\circ_\infty$:
\begin{equation}\label{dzerob}
d^\circ_\infty(a,b) = \Lab_\infty(a) + \Lab_\infty(b) - 2 \max \lp \min_{ x \in [a,b]} \Lab_\infty(x),\min_{ x \in [b,a]} \Lab_\infty(x) \rp,
\end{equation}
where $[a,b]$ was defined by~\eqref{abb}. Similarly, for $a\in \f_n$, we set $\Lab_n(a) \de \lab_n(a)+\bb_n(\aaa(a)-1)$, so that $\Lab_n(\f_n(i))=\Lab_n(i)$ for all $0 \le i \le 2n+\sigma_n -1$.

\subsection{Point identifications}\label{secipb}

\subsubsection{Criterion telling which points are identified}

Our analysis starts with the following two observations on the process $(C_\infty(s), \Lab_\infty(s))_{0\le s\le 1}$.

\begin{lem}\label{minb}
The set of points where $\Lab_\infty$ reaches its minimum is a.s.\ a singleton.
\end{lem}

Let $f:[0,\ell] \to \R$ be a continuous function. We say that $s\in [0,\ell)$ is a \textbf{right-increase point} of~$f$ if there exists $t\in (s,\ell]$ such that $f(r) \ge f(s)$ for all $s\le r \le t$. A \textbf{left-increase point} is defined in a symmetric way. We denote by $\IP(f)$ the set of all (left or right) increase points of~$f$.

\begin{lem}\label{pcb}
Almost surely, $\IP(C_\infty)$ and $\IP(\Lab_\infty)$ are disjoint sets.
\end{lem}

The proofs of these lemmas make intensive use of the so-called Brownian snake, so that we postpone them to Section~\ref{secsnake}. We have the following criterion:

\begin{thm}\label{ipb}
Almost surely, for every $a,b \in \Fi$, $a\sim_\infty b$ is equivalent to $d^\circ_\infty(a,b)=0$. In other words,
$$\disig(a,b)=0\ \Leftrightarrow\ d^\circ_\infty(a,b)=0.$$
\end{thm}

We call \textbf{leaves} the points of~$\Fi$ whose equivalence class for~$\simeq_\infty$ is trivial. It will be important in what follows to observe that, by Lemma~\ref{pcb} and Theorem~\ref{ipb}, only leaves of~$\Fi$ can be identified by~$\sim_\infty$.

The proof of Theorem~\ref{ipb} is based on Lemma~\ref{minb}, Lemma~\ref{pcb}, and Lemma~\ref{legalllemb} below, which we will prove in Section~\ref{secsnake}. Once we have these lemmas, the arguments of the proof of \cite[Theorem~8]{bettinelli10tsl} (which uses the ideas of~\cite{legall07tss}) may readily be adapted to our case. For the sake of self-containment, we give here the main ingredients. By the bound~\eqref{boundsig0}, we already have one implication:
$$d^\circ_\infty(a,b)=0 \ \Rightarrow \ \disig(a,b)=0.$$
The converse is shown in three steps. First, we show that the floor points are not identified (by~$\sim_\8$) with any other points, then that points are not identified with their strict ancestors, and finally the general case. 
As an example, we will treat here the first step mentioned above. As we will see, the adaptation is almost verbatim, and is a little easier. The other steps use the same ideas and are even more straightforwardly adaptable, so that we leave them to the reader. Precisely, we are going to show the following lemma:

\begin{lem}\label{lemid}
Almost surely, for every $b\in \Fi$ and every $a\in\fli\bs\{\rho_b\}$, we have $a\not\sim_\infty b$.
\end{lem}

\subsubsection{Set overflown by a path and paths passing through subtrees}\label{secof}

We give in this section the two notions we will need for discrete paths. In the following, we will never consider paths using the edges of the forest, but always paths using the edges of the map, and we will always use the letter ``$\wp$'' to denote these paths.

The first notion is the notion of a set overflown by a path: roughly speaking, imagine a squirrel jumping from tree to tree in the forest along the edges of a path~$\wp$ in the map. Then the set overflown by~$\wp$ is the ground covered by the squirrel along its journey. Let us denote by~$\fl_n$ the floor of~$\f_n$. Let $i\in \ent 0{2n+\sigma_n-1}$, and let $\suc(i)$ be its successor in $(\f_n,\lab_n)$, defined by~\eqref{succf}. We moreover suppose that $\suc(i) \neq \infty$. We say that the arc $i \arc \suc(i)$ linking $\f_n(i)$ to $\f_n(\suc(i))$ overflies the set
$$\f_n \lp \overrightarrow{\ent i {\suc(i)}} \rp \cap \fl_n,$$
where $\overrightarrow{\ent i {\suc(i)}}$ was defined by~\eqref{oraijb}. We define the set overflown by a path~$\wp$ in~$\q_n$ that avoids the base point $v_n^\bullet$ as the union of the sets its arcs overfly. Beware that, in this definition, the orientations of the edges is not taken into account. In particular, the reverse of a path~$\wp$ overflies the same set as the path~$\wp$ itself.

\begin{figure}[ht]
		\psfrag{o}[][]{}
		\psfrag{g}[][][1.2]{\textcolor{red}{$\wp$}}
	\centering\includegraphics{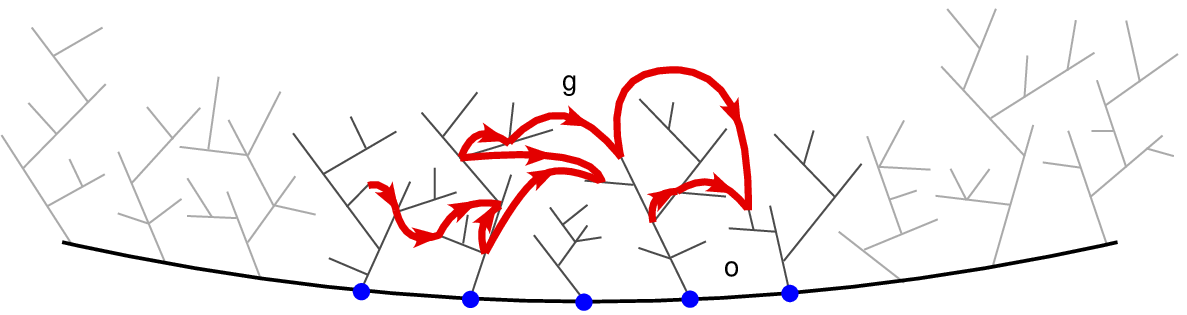}
	\caption[Set overflown by a path.]{The set overflown by the path~$\wp$ is the set of (blue) large dots. Note that the middle tree is overflown although it is not visited. The arrows on the arcs of~$\wp$ point from corners to their successors in the \BDG bijection.}
\end{figure}

\begin{rem}
Note that, by the \BDG construction, all the labels of the set overflown by a path are larger than or equal to the minimum label on the path. Note also that the set overflown by a path is a connected subset of~$\fl_n$.
\end{rem}


\bigskip

The second notion is the notion of path passing through a subtree: here again, imagine a squirrel moving along the path~$\wp$. The path~$\wp$ passes through a subtree~$\tau$ if the squirrel visits~$\tau$, and moreover enters it when going in one direction (from left to right or from right to left) and exits it while going in the same direction. Let~$\tau$ be a subtree of~$\f_n$ and $\wp=(\wp(0),\wp(1),\dots,\wp(r))$ a path in~$\q_n$ that avoids the base point~$v_n^\bullet$.
We say that the path~$\wp$ passes through the subtree~$\tau$ between times~$i$ and~$j$, where $0 < i \le j < r$, if
\begin{itemize}
	\item $\wp(i-1) \notin \tau$; $\wp(\ent i j) \subseteq \tau$; $\wp(j+1) \notin \tau$,
	\item $\Lab_n(\wp(i)) - \Lab_n(\wp(i-1)) = \Lab_n(\wp(j+1)) - \Lab_n(\wp(j))$.  
\end{itemize}

We say that a sequence of vertices $a_n\in \f_n$ converges toward a point $a\in \Fi$ if there exists a sequence of integers $s_n \in \ent 0 {2n+\sigma_n-1}$ coding~$a_n$ (i.e.\ $a_n=\f_n(s_n)$) such that $s_n/(2n+\sigma_n-1)$ admits a limit $s$ coding~$a$, i.e.\ such that $a=\Fi(s)$. Let $\ent {l_n}{r_n}$ be the intervals coding subtrees $\tau_n \subseteq \f_n$. We say that the subtree $\tau_n$ converges toward a subtree $\tau \subseteq \Fi$ if the sequences $l_n/(2n+\sigma_n-1)$ and $r_n/(2n+\sigma_n-1)$ admit limits $l$ and $r$ such that the interval coding $\tau$ is $[l, r]$. The key lemma of our approach is the following. It is adapted from Le~Gall \cite[End of Proposition~4.2]{legall07tss}, and will be proved in Section~\ref{secsnake}.

\begin{lem}\label{legalllemb}
With full probability, the following occurs. Let $a, b \in \Fi$ be such that $\Lab_\infty(a)=\Lab_\infty(b)$. We suppose that there exists a subtree~$\tau$ rooted at~$\rho$ such that $\inf_{\tau} \Lab_\infty < \Lab_\infty(a) < \Lab_\infty(\rho)$. We further suppose that we can find vertices~$a_n$, $b_n \in \f_n$ and subtrees~$\tau_n$ in~$\f_n$ converging respectively toward~$a$, $b$, $\tau$ and satisfying the following property: for infinitely many~$n$'s, there exists a geodesic path~$\wp_n$ in~$\q_n$ from~$a_n$ to~$b_n$ that avoids the base point~$v_n^\bullet$ and passes through the subtree~$\tau_n$.

Then, $a\not\sim_\infty b$.
\end{lem}

\subsubsection{Proof of Lemma~\ref{lemid}}

\begin{pre}[Proof of Lemma~\ref{lemid}]
We argue by contradiction and suppose that we can find $b\in \Fi$ and $a\in\fli\bs\{\rho_b\}$ such that $a\sim_\infty b$. It is easy to find $a_n \in \fl_n$ and $b_n \in \f_n$ converging respectively toward~$a$ and~$b$. Let~$\wp_n$ be a geodesic path (in~$\q_n$, for $d_{\q_n}$) from~$a_n$ to~$b_n$. For~$n$ large, $\wp_n$ avoids the base-point, because otherwise, $a$ and~$b$ would have the minimal label and this would contradict Lemma~\ref{minb}. For such an~$n$, $\wp_n$ has to overfly at least $\lhb \rho_{b_n}, a_n \rhb$ or $\lhb a_n, \rho_{b_n} \rhb$. To see this, let $(x,y) \in \lhb \rho_{b_n}, a_n \rhb \times \lhb a_n, \rho_{b_n} \rhb$. When we remove from~$\f_n$ all the edges incident to~$x$ and all the edges incident to~$y$, we obtain several connected components, and the points~$a_n$ and~$b_n$ do not belong to the same of these components. There has to be an arc of~$\wp_n$ that links a point belonging to the component containing~$a_n$ to one of the other components. Such an arc overflies~$x$ or~$y$.

Let us suppose that, for infinitely many $n$'s, $\wp_n$ overflies $\lhb \rho_{b_n}, a_n \rhb$. By the remark concerning the labels on the set overflown by a path in the previous section, a simple argument (see \cite[Lemma~14]{bettinelli10tsl}) shows that $\Lab_\infty(c) \ge \Lab_\infty(a) = \Lab_\infty(b)$ for all $c \in \lhb \rho_{b}, a \rhb$. The labels on~$\fli$ are given by the process~$\bb_\infty$, defined during Proposition~\ref{cvcont}: for $x \in [0,\sigma]$, we define $T_x \de \inf \{ r \ge 0 \, :\, C_\infty(r)= \sigma - x \}$, so that $\fli = \Fi(\{T_x,\, 0 \le x \le \sigma \})$, and
$$\big(\Lab_\infty(T_x)\big)_{0 \le x \le \sigma}=\big(\bb_\infty(x)\big)_{0 \le x \le \sigma}.$$
As~$\bb_\infty$ has the law of a certain Brownian bridge (scaled by~$\sqrt 3$), and as local minimums of Brownian motion are distinct, we can find $d\in \lhb \rho_{b}, a \rhb \bs \{a,\rho_b\}$ such that $\Lab_\infty(c) > \Lab_\infty(a)$ for all $c \in \lhb d,a \rhb\bs \{a\}$.

Because $a\in \fli$, every number coding it is an increase point of $C_\infty$ and thus is not an increase point of $\Lab_\infty$, by Lemma~\ref{pcb}. As a result, there exists a tree $\tau^1$ rooted at $\rho^1 \in \lhb d,a \rhb\bs \{a\}$ satisfying $\inf_{\tau^1} \Lab_\infty < \Lab_\infty(a) < \Lab_\infty(\rho^1)$ (see Figure~\ref{tau1}).

\begin{figure}[ht]
		\psfrag{a}[][]{$a$}
		\psfrag{b}[][]{$b$}
		\psfrag{r}[][]{$\rho_b$}
		\psfrag{d}[][]{$d$}
		\psfrag{1}[][][.8]{$\rho^1$}
		\psfrag{i}[][]{$\tau^1$}
		\psfrag{m}[r][r][.8]{$\inf \Lab_\infty < \Lab_\infty(a)$}
		\psfrag{l}[r][r][.8]{\textcolor{red}{$\Lab_\infty > \Lab_\infty(a)$}}
		\psfrag{A}[][][.8]{\textcolor{blue}{$\lhb \rho_b, a \rhb$}}
	\centering\includegraphics{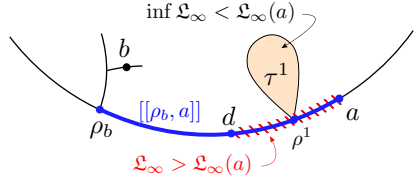}
	\caption[Visual aid for the proof of Lemma~\ref{lemid}.]{The tree $\tau^1$.}
	\label{tau1}
\end{figure}

Similarly, if for infinitely many $n$'s, $\wp_n$ overflies $\lhb a_n, \rho_{b_n} \rhb$, then we can find a tree $\tau^2$ rooted at $\rho^2 \in \lhb a, \rho_b \rhb\bs \{a,\rho_b\}$ satisfying $\inf_{\tau^2} \Lab_\infty < \Lab_\infty(a) < \Lab_\infty(\rho^2)$. Three cases may occur:
\begin{enumerate}[($i$)]
	\item for $n$ large enough, $\wp_n$ does not overfly $\lhb a_n, \rho_{b_n} \rhb$ (and therefore overflies $\lhb \rho_{b_n}, a_n \rhb$),
	\item for $n$ large enough, $\wp_n$ does not overfly $\lhb \rho_{b_n}, a_n \rhb$ (and therefore overflies $\lhb a_n, \rho_{b_n} \rhb$),
	\item $\wp_n$ overflies $\lhb \rho_{b_n}, a_n \rhb$ for infinitely many $n$'s, and $\lhb a_n, \rho_{b_n} \rhb$ also for infinitely many $n$'s.
\end{enumerate}

In case ($i$), the tree $\tau^1$ is well defined. Let $\tau_n^1 \subseteq \f_n$ be a tree rooted at $\rho_n^1 \in \lhb \rho_{b_n}, a_n \rhb$ converging to~$\tau^1$. We claim that, for $n$ sufficiently large, $\wp_n$ passes through $\tau_n^1$. First, notice that by continuity, for $n$ large enough, $\inf_{\tau_n^1} \Lab_n < \inf_{\wp_n} \Lab_n$. The idea is that, at some point, $\wp_n$ has to go from a tree located at the right of $\tau_n^1$ to a tree located at its left, and, because it does not overfly $\lhb a_n, \rho_{b_n} \rhb$, it has no other choice than passing through~$\tau^1_n$ (see Figure~\ref{taun1}).

More precisely, we denote by $\ent{s_n^1}{t_n^1}$ the set coding the subtree $\tau_n^1$, and we let $\omega_n = \f_n(p_n) \in \lhb a_n, \rho_{b_n} \rhb$ be a point that is not overflown by~$\wp_n$. Then, we define
$$A_n \de \f_n\lp \overrightarrow{\ent{t_n^1+1}{p_n}} \rp.$$
We denote by $\wp_n(i-1)$ the last point of~$\wp_n$ belonging to~$A_n$. Such a point exists because $a_n \in A_n$ and $b_n \notin A_n$. For $n$ large, because~$\wp_n$ does not overfly $\omega_n$, and because $\inf_{\tau_n^1} \Lab_n < \inf_{\wp_n} \Lab_n$, we see that $\wp_n(i) \in \tau_n^1$. Let $\wp_n(j+1)$ be the first point after $\wp_n(i)$ not belonging to $\tau_n^1$. It exists because $b_n \notin \tau_n^1$. Using the facts that~$\wp_n$ does not overfly~$\omega_n$, and that $\wp_n(j+1)\notin A_n$, we see that~$\wp_n$ passes through~$\tau_n^1$ between times~$i$ and~$j$.

\begin{figure}[ht]
		\psfrag{c}[][][1.2]{\textcolor{red}{$\wp_n$}}
		\psfrag{A}[][][.9]{$A_n$}
		\psfrag{a}[][]{$a_n$}
		\psfrag{o}[l][l]{$\omega_n$}
		\psfrag{s}[][][.8]{$s_n^1$}
		\psfrag{t}[][][.8]{$t_n^1$}
		\psfrag{p}[][][.8]{$p_n$}
		\psfrag{i}[][]{$\tau_n^1$}
	\centering\includegraphics{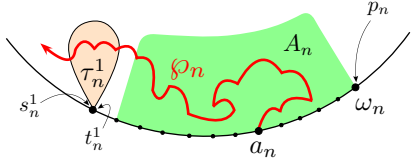}
	\caption[Visual aid for the proof of Lemma~\ref{lemid}.]{The path~$\wp_n$ passing through the tree~$\tau_n^1$.}
	\label{taun1}
\end{figure}

In case ($ii$), we apply the same reasoning with $\tau^2$ instead of $\tau^1$. In case ($iii$), both trees~$\tau^1$ and~$\tau^2$ are well defined and we obtain that~$\wp_n$ has to pass through one of their discrete approximations. We then conclude by Lemma~\ref{legalllemb} that $a\not\sim_\infty b$, which contradicts our hypothesis.
\end{pre}

\section{Regularity of quadrangulations}\label{secsurfb}

Recently, the notion of regularity has been used to identify the topology of the scaling limit of random uniform planar quadrangulations in~\cite{miermont08sphericity}, and then positive genus quadrangulations in~\cite{bettinelli10tsl}. In both these references, it is the notion of 1-regularity that is used, roughly stating that there are no small loops separating the surface into large components. In the case of surfaces with a boundary, a new problem arises, and we also need the notion of 0-regularity for the boundary. In this section, we present both these notions, which were introduced in a slightly different context (see the discussion in \cite[Section~2]{miermont08sphericity}) by Whyburn \cite{whyburn35sls, whyburn35rcm}, and then use them to prove Theorem~\ref{cvqtop}.

\subsection{0-regularity and 1-regularity}\label{secreg}

Recall that we wrote $(\MM,\dGH)$ for the set of isometry classes of compact metric spaces, endowed with the Gromov--Hausdorff metric. A metric space $(\X,\delta)$ is called a \textbf{path metric space} if any two points $x,y \in \X$ can be joined by a path isometric to the segment $[0,\delta(x,y)]$. We let~$\PM$ be the set of isometry classes of path metric spaces, which is a closed subset of~$\MM$, by \cite[Theorem~7.5.1]{burago01cmg}.

\begin{defi}
We say that a sequence $(\X_n)_{n \ge 1}$ of compact metric spaces is \textbf{1-regular} if for every $\eps > 0$, there exists $\eta >0$ such that for $n$ large enough, every loop of diameter less than~$\eta$ in~$\X_n$ is homotopic to~$0$ in its $\eps$-neighborhood.
\end{defi}

The theorem (derived from \cite[theorem~7]{begle44rc}) that was used in \cite{miermont08sphericity,bettinelli10tsl} states that the limit of a converging 1-regular sequence of path metric spaces all homeomorphic to the $g$-torus is either reduced to a singleton (this can only happen when $g=0$), or homeomorphic to the $g$-torus as well. In other words, this gives a sufficient condition for the limit to be homeomorphic to the surface we started with. In the case of the $2$-dimensional disc~$\Ddd_2$, this condition is no longer sufficient. For example, take for the space~$\X_n$ the union of two unit discs whose centers are at distance $2-1/n$. This peanut-shaped space is homeomorphic to~$\Ddd_2$, and it is easy to see that~$(\X_n)_n$ is 1-regular and converges to the wedge sum (or bouquet) of two discs. The following definition discards this kind of degeneracy.

\begin{defi}
We say that a sequence $(\X_n)_{n \ge 1}$ of compact metric spaces is \textbf{0-regular} if for every $\eps > 0$, there exists $\eta >0$ such that for~$n$ large enough, every pair of points in~$\X_n$ lying at a distance less than~$\eta$ from each other belong to a connected subset of~$\X_n$ of diameter less than~$\eps$.
\end{defi}

We will rely on the following theorem, which is a simple consequence of \cite[Theorem~6.4]{whyburn35sls}. Recall that the boundary of a path metric space $(\X,\delta)$ is the set $\partial \X \subseteq \X$ of points having no neighborhood homeomorphic to a disc, equipped with the restriction of the metric $\delta$.

\begin{prop}[Whyburn]\label{whyburn}
Let $(\X_n)_{n \ge 1}$ be a sequence of path metric spaces all homeomorphic to the $2$-dimensional disc~$\Ddd_2$, converging for the Gromov--Hausdorff topology toward a metric space~$\X$ not reduced to a single point. Suppose that the sequence $(\X_n)_{n \ge 1}$ is 1-regular, and that the sequence $(\partial \X_n)_n$ is 0-regular. 

Then~$\X$ is homeomorphic to~$\Ddd_2$ as well.
\end{prop}

In~\cite{whyburn35sls}, Whyburn actually considered convergence in the sense of the Hausdorff topology, and made the extra hypothesis that~$\partial \X_n$ converges to a set~$B$. He concluded that~$\X$ was homeomorphic to $\Ddd_2$ and that $\partial\X$ was equal to~$B$. To derive the version that we state here, we proceed as follows. First, by \cite[Lemma~A.1]{greven09cdr}, we can find a compact metric space $\ZZ$, and isometric embeddings $\varphi$, $\varphi_1$, $\varphi_2$, \dots~of $\X$, $\X_1$, $\X_2$, \dots~into~$\ZZ$ such that $\varphi_n(\X_n)$ converges toward $\varphi(\X)$ for the Hausdorff topology in~$\ZZ$. Then, by \cite[Theorem~7.3.8]{burago01cmg}, the family $\{\partial(\varphi_n(\X_n))\}$ is relatively compact for the Hausdorff topology. Let us consider a subsequence along which $\partial(\varphi_n(\X_n))$ converges to a set~$B$. Applying Whyburn's original theorem along this subsequence, we obtain that $\varphi(\X)$ is homeomorphic to $\Ddd_2$, so that $\X$ is homeomorphic to $\Ddd_2$ as well. We moreover obtain that $\partial(\varphi(\X))=B$, and, using the same argument, we see that any accumulation point of the sequence $(\partial(\varphi_n(\X_n)))_n$ has to be $\partial(\varphi(\X))$, so that $\partial(\varphi_n(\X_n))=\varphi_n(\partial\X_n)$ actually converges toward $\partial(\varphi(\X))=\varphi(\partial\X)$ for the Hausdorff topology. This last observation will be used in Section~\ref{secboundary} to identify the boundary of~$\qis$.

\subsection{Representation as metric surfaces}\label{secrepsurf}

As the space $(V(\q_n),d_{\q_n})$ is not a surface, we cannot directly apply Proposition~\ref{whyburn}. In a first time, we will construct a path metric space $(\cS_n,\delta_n)$ homeomorphic to~$\Ddd_2$, and an embedded graph that is a representative of the map~$\q_n$, such that the restriction of $(\cS_n,\delta_n)$ to the embedded graph is isometric to $( V(\q_n), d_{\q_n})$. We use the same method as Miermont in \cite[Section~3.1]{miermont08sphericity} (see also \cite[Section~5.2]{bettinelli10tsl}), roughly consisting in gluing hollow boxes together according to the structure of~$\q_n$.

Let us be a little more specific. Let~$f_*$ be the external face of~$\q_n$, $F(\q_n)$ its set of internal faces, and $F_*(\q_n) \de F(\q_n)\cup\{f_*\}$ the set of all its faces. Let also~$\calG$ be a regular $2\sigma_n$-gon with unit length edges embedded in~$\R^2$, and let us denote by~$z_k$, $0\le k \le 2\sigma_n$, its vertices (with $z_0=z_{2\sigma_n}$). With every quadrangle $f \in F(\q_n)$, we associate a copy of the ``hollow bottomless unit cube,'' and with~$f_*$, we associate a ``hollow bottomless $2\sigma_n$-sided prism'': we define
$$X_f \de [0,1]^3\bs \big((0,1)^2 \times [0,1) \big), \quad f \in F(\q_n),  \sand X_{f_*} \de \big( \mathcal G \times [0,1] \big) \bs \big( \mathring{\mathcal G} \times [0,1) \big),$$
where $\mathring{\mathcal G}$ denotes the interior of $\mathcal G$, and we endow these spaces with the intrinsic metric $D_{f}$ inherited from the Euclidean metric. This means that the distance between two points~$x$ and~$y$ is the Euclidean length of a minimal path in~$X_f$ linking~$x$ to~$y$. Note in particular that if~$x$ and~$y$ are on the boundary, this path is entirely contained in the boundary. This will ensure that, when we will glue these spaces together, we will not alter the graph metric. Note also that, so far, the external face is not really treated differently from the other faces (except for the fact that it has a different number of edges). In the end, we will remove the ``top'' $\mathring{\calG} \times \{1\}$ from~$X_{f_*}$.

Now, we associate with every half-edge $e\in \vec E(\q_n)$ a path~$c_e$ parameterizing the corresponding edge of the polygon~$\partial X_f$, where~$f$ is the face incident to~$e$. We denote by $e_1$, $e_2$, \dots, $e_{2\sigma_n}$ the half-edges bordering~$f_*$ ordered in the clockwise order (recall that, by convention, $f_*$ is the infinite face of $\q_n$, so that the order is reversed), and define
$$c_{e_k}(t)	\de\big((1-t)z_{k-1} + t z_k,0\big) \in X_{f_*}, \qquad t\in[0,1],\quad 1 \le k \le 2\sigma_n.$$
In a similar way, for every internal face $f \in F(\q_n)$, and every half-edge~$e$ incident to it, we define a function $c_e:[0,1]\to \partial X_f$ parameterizing an edge of~$\partial X_f$. We do this in such a way that the parameterization of~$\partial X_f$ is coherent with the counterclockwise order around~$f$ (see \cite[Section~3.1]{miermont08sphericity} or \cite[Section~5.2]{bettinelli10tsl}).

We may now glue these spaces together along their boundaries: we define the relation~$\approx$ as the coarsest equivalence relation for which $c_e(t) \approx c_{\bar e}(1-t)$ for all $e\in \vec E(\q_n)$ and $t \in [0,1]$, where~$\bar e$ denotes the reverse of~$e$. The topological quotient $\hat\cS_n \de (\bigsqcup_{f\in F_*(\q_n)} X_f)_{/\approx}$ is then a $2$-dimensional CW-complex satisfying the following properties. Its $1$-skeleton $\EE_n = (\bigsqcup_{f\in F_*(\q_n)} \partial X_f)_{/\approx}$ is an embedding of~$\q_n$ with faces $X_f \bs \partial X_f$. The edge $\{e,\bar e\} \in E(\q_n)$ corresponds to the edge of~$\hat\cS_n$ made of the equivalence classes of the points in $c_e([0,1])$. Its $0$-skeleton $\V_n$ is in one-to-one correspondence with $V(\q_n)$, and its vertices are the equivalence classes of the vertices of the polygons~$\partial X_f$'s.

We endow the space $\bigsqcup_{f\in F_*(\q_n)} X_f$ with the largest pseudo-metric $\delta_n$ compatible with $D_f$, $f\in F_*(\q_n)$ and $\approx$, in the sense that $\delta_n(x,y) \le D_f(x,y)$ for $x,y\in X_f$ and $\delta_n(x,y)=0$ whenever $x \approx y$. Its quotient, which we still denote by~$\delta_n$, then defines a pseudo-metric on $\hat\cS_n$ (which is actually a true metric, as we will see in Proposition~\ref{surfp}). We also define $\delta_{(n)} \de \delta_n /(\g)$ its rescaled version. Finally, we set $\cS_n \de (\bigsqcup_{f\in F_*(\q_n)} Y_f)_{/\approx} \subseteq \hat\cS_n$, where $Y_{f_*}\de X_{f_*}\bs (\mathring{\mathcal G} \times \{1\})$ and $Y_f \de X_f$ when $f \neq f_*$.

\begin{prop}[{\cite[Proposition~1]{miermont08sphericity}}]\label{surfp}
The space $(\hat\cS_n,\delta_n)$ is a path metric space homeomorphic to~$\Sss_2$. Moreover, the metric space $(\V_n,\delta_n)$ is isometric to $(V(\q_n), d_{\q_n})$, and any geodesic path in $\hat\cS_n$ between points in $\V_n$ is a concatenation of edges of $\hat\cS_n$.
\end{prop}

We readily obtain the following corollary.

\begin{corol}\label{surfb}
The space $(\cS_n,\delta_n)$ is a path metric space homeomorphic to $\Ddd_2$. Moreover, the metric space $(\V_n,\delta_n)$ is isometric to $(V(\q_n), d_{\q_n})$, and any geodesic path in~$\cS_n$ between points in~$\V_n$ is a concatenation of edges of $\cS_n$. Finally, $\dGH \big( (V(\q_n), d_{\q_n}), (\cS_n,\delta_n) \big) \le 3$, so that, by Theorem~\ref{cvqb},
$$\big(\cS_{n_k},\delta_{(n_k)} \big) \tolk (\qis, \disig)$$
in the sense of the Gromov--Hausdorff topology.
\end{corol}

Note that, although the boundary of $\q_n$ is not topologically a circle in general, $\partial\cS_n$ (which corresponds to $\partial\mathcal G \times \{1\}$ in $Y_{f_*}$) always is. In what follows, we will see $V(\q_n)$ as a subset of~$\cS_n$. In other words, we identify~$\V_n$ with $V(\q_n)$.

\subsection{Proof of Theorem~\ref{cvqtop}}

We now prove that $(\qis,\disig)$ is a.s.\ homeomorphic to~$\Ddd_2$ thanks to Proposition~\ref{whyburn} and Corollary~\ref{surfb}. As $(\qis,\disig)$ is a.s.\ not reduced to a point\footnote{It is for example a.s.\ of Hausdorff dimension~4 by Theorem~\ref{cvqb}.}, it is enough to show that the sequence $(\partial \cS_{n_k})_k$ is 0-regular, and that the sequence $(\cS_{n_k})_k$ is 1-regular. The 1-regularity of $(\cS_{n_k})_k$ is readily adaptable from \cite[Section~5.3]{bettinelli10tsl} so that we begin with the 0-regularity of the boundary. We denote by $\pii : \Fi \to \qis$ the canonical projection.

\subsubsection{0-regularity of the boundary}

\begin{lem}\label{lem0reg}
The sequence $(\partial \cS_{n_k})_k$ is 0-regular.
\end{lem}

\begin{pre}
The idea is that $\fli$ has no cut points in $\Fi$, and because the points in $\fli$ are not identified with any other points, $\pii(\fli)$ does not have any cut points either.

We argue by contradiction and assume that, with positive probability, along some (random) subsequence of the sequence $(n_k)_{k\ge 0}$, there exist $\eps >0$, $x_n$, $y_n \in \partial\cS_n$ such that $\delta_{(n)}(x_n,y_n) \to 0$, and~$x_n$ and~$y_n$ do not belong to the same connected component of $B_{(n)}(x_n,\eps) \cap \partial\cS_n$, where $B_{(n)}(x_n,\eps)$ denotes the open ball of radius~$\eps$ centered at~$x_n$ for the metric~$\delta_{(n)}$. We reason on this event.

As~$x_n$ and~$y_n$ do not belong to the same connected component of $B_{(n)}(x_n,\eps) \cap \partial\cS_n$, we can find~$x_n'$, $y_n'\in \partial\cS_n\bs B_{(n)}(x_n,\eps)$ such that~$x_n'$ belongs to one of the two arcs joining~$x_n$ to~$y_n$ in $\partial\cS_n$, and such that~$y_n'$ belongs to the other one. We are going to approach these four points with points of~$\fl_n$.

We denote by $\partial\q_n \subseteq \vec E(\q_n)$ the set of half-edges incident to the external face of~$\q_n$. With every point $x \in\partial \cS_n$ naturally corresponds a half-edge $e(x)\in \partial\q_n$: if $x$ corresponds to $((1-t)z_{k-1} + t z_k,1) \in X_{f_*}$ for some $t\in [0,1)$, then~$e(x)$ is the half-edge~$e_k$. We consider the first half-edge $e\in \partial\q_n$ after $e(x_n)$ ($e(x_n)$ included) in the clockwise order such that $\Lab_n(e^+)=\Lab_n(e^-)+1$, and we set $a_n \de e^+$. By definition of the \BDG bijection, $a_n \in \fl_n$. Moreover, $a_n$ is ``close'' to~$x_n$, in the sense that $\delta_n(a_n,x_n) \le 1 + \sup_{0 \le i < 2\sigma_n} |\bb_n(i+1) - \bb_n(i)+2|$, so that
$$\delta_{(n)}(a_n,x_n) \le \frac{3}{\g} + \sup_{x} \lt \bb_{(n)}\lp x+ \frac{1}{\sqrt{2n}} \rp - \bb_{(n)}(x)\rt \le \frac{3}{\g} + \omega_{\bb_{(n)}}(\eta),$$
as soon as $n \ge 1/2\eta^2$. Here, $\omega_{\bb_{(n)}}$ denotes the modulus of continuity of $\bb_{(n)}$. Hence, we obtain that $\limsup \delta_{(n)}(a_n,x_n) \le \omega_{\bb_{\infty}}(\eta)$, for all $\eta >0$, so that $\delta_{(n)}(a_n,x_n) \to 0$.

\begin{figure}[ht]
		\psfrag{+}[][][.6]{$+$}
		\psfrag{-}[][][.6]{$-$}
		\psfrag{0}[][][.8]{$6$}
		\psfrag{1}[][][.8]{$7$}
		\psfrag{2}[][][.8]{$8$}
		\psfrag{q}[][]{$\partial\q_n$}
		\psfrag{s}[][]{$\partial\cS_n$}
		\psfrag{a}[][]{$a_n$}
		\psfrag{x}[][]{$x_n$}
		\psfrag{e}[][]{\textcolor{red}{$e(x_n)$}}
		\psfrag{f}[][]{\textcolor{vert}{$\fl_n$}}
	\centering\includegraphics{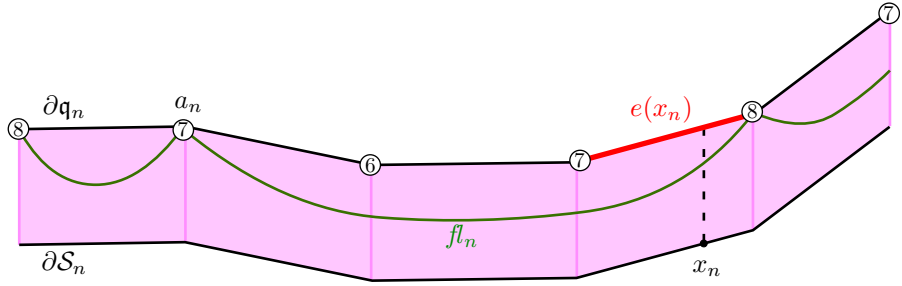}
	\caption[Approaching a point in $\partial\cS_n$ with a point in $\fl_n$.]{Approaching a point $x_n\in\partial\cS_n$ with a point $a_n\in\fl_n$.}
	\label{xnan}
\end{figure}

We define in a similar way points~$b_n$, $a_n'$, and~$b_n'$ in~$\fl_n$ corresponding to~$y_n$, $x_n'$, and~$y_n'$. Exchanging~$x_n'$ and~$y_n'$ if necessary, we may suppose that the points~$a_n$, $a_n'$, $b_n$, $b_n'$ are encountered in this order when traveling in the counterclockwise order around $\partial\q_n$. Up to further extraction, we may suppose that $(a_n,a_n',b_n,b_n') \to (a,a',b,b')\in \fli^4$, so that $a' \in \lhb a,b \rhb$ and $b' \in \lhb b,a \rhb$. Moreover, because $\delta_{(n)}(x_n,x_n') \ge \eps$, we see that $\disig(a,a') \ge \eps$. Similarly, we obtain that $\disig(b,a') \ge \eps$, $\disig(a,b') \ge \eps$, and $\disig(b,b') \ge \eps$, so that $a\neq b$. Finally, the fact that $\delta_{(n)}(x_n,y_n) \to 0$ implies that $\disig(a,b)=0$, so that $a \sim_\infty b$. This contradicts Lemma~\ref{lemid}.
\end{pre}

\subsubsection{1-regularity of \hr{$\cS_n$}{Sn}}

In order to show that the sequence $(\cS_{n_k})_k$ is 1-regular, we first only consider simple loops made of edges in~$\cS_n$. A simple loop~$\wp$ splits~$\cS_n$ into two domains. By the Jordan curve theorem, one of these is homeomorphic to a disc. We call it the \textbf{inner domain} of~$\wp$. The other domain contains~$\partial\cS_n$ in its closure, and we call it the \textbf{outer domain} of~$\wp$.

\begin{lem}\label{1regb}
A.s., for all $\eps >0$, there exists $0 < \eta < \eps$ such that for all~$k$ sufficiently large, the inner domain of any simple loop made of edges in $\cS_{n_k}$ with diameter less than~$\eta$ has diameter less than~$\eps$.
\end{lem}

The proof of this Lemma is readily adaptable from the proof of \cite[Lemma~22]{bettinelli10tsl}, which uses the method employed by Miermont in \cite{miermont08sphericity}. The general idea is that a loop separates some part of the map from the base point. As a result, the labels in one of the two domains are larger than the labels on the loop. In the forest, this corresponds to having a part with labels larger than the labels on the ``border.'' In the continuous limit, this creates an increase point for both~$C_\infty$ and~$\Lab_\infty$. We recall now the main steps.

\begin{pre}
We argue by contradiction and suppose that, with positive probability, there exists $\eps > 0$ for which, along some (random) subsequence of the sequence $(n_k)_{k\ge 0}$, there exist simple loops~$\wp_n$ made of edges in~$\cS_n$ with diameter tending to~$0$ (with respect to the rescaled metric $\delta_{(n)}$) and whose inner domains are of diameter larger than~$\eps$. We reason on this event. We will show in the proof of Proposition~\ref{propboundary} that $\partial\cS_{n_k}$ tends, for the Gromov--Hausdorff topology, toward~$\pii(\fli)$. Because~$\fli$ is not reduced to a singleton, we see by Lemma~\ref{lemid} that~$\pii(\fli)$ is not a singleton either, so that $\diam(\pii(\fli))>0$. To avoid trivialities, we moreover suppose that $\eps < \diam(\pii(\fli))$. Because $\partial\cS_n$ is entirely contained in the outer domain of~$\wp_n$, we obtain that, for~$n$ large enough, the outer domain of~$\wp_n$ is also of diameter larger than~$\eps$.

Let~$s_n^\bullet$ be an integer where~$\Lab_n$ reaches its minimum, and $w^\bullet_n \de \f_n(s_n^\bullet)$ the corresponding point in the forest. Let us suppose for the moment that $w_n^\bullet \notin \wp_n$. We take~$x_n$ as far as possible from~$\wp_n$ in the connected component of the complement of~$\wp_n$ that does not contain~$w_n^\bullet$, and we denote by~$y_n$ the first vertex of the path $\lhb x_n, w^\bullet_n \rhb$ that belongs to~$\wp_n$. Up to further extraction, we suppose that $s_n^\bullet / (2n+\sigma_n-1) \to s^\bullet\de\operatorname{argmin} \Lab_\infty$, $x_n \to x$, and $y_n \to y$. Because of the way~$x_n$ and~$y_n$ were chosen, it is not hard to see that $x \neq y$.

Let us first suppose that $y \neq w^\bullet \de \Fi(s^\bullet)$. In particular, $w_n^\bullet \notin \wp_n$ for~$n$ large, so that~$x_n$ and~$y_n$ are well defined. In this case, $y \in \lhb x, w^\bullet \rhb \bs \{x,w^\bullet\}$, so that the points in $\Fi^{-1}(y)$ are increase points of $C_\infty$. By Lemma~\ref{pcb}, we can find a subtree~$\tau$, not containing~$y$, satisfying $\inf_{\tau} \Lab_\infty < \Lab_\infty(y)$ and rooted on $\lhb x, y\rhb$. We consider a discrete approximation~$\tau_n$ of this subtree, rooted on $\lhb x_n, y_n \rhb$. As the labels on~$\wp_n$ differ by $o(n^{1/4})$, when~$n$ is sufficiently large, we thus have $\inf_{\tau_n} \Lab_n < \inf_{\wp_n} \Lab_n $.

As the labels of the forest represent the distances in~$\q_n$ to the base point (up to some additive constant), we see that all the labels of the points in the same domain as~$x_n$ are larger than $\inf_{\wp_n} \Lab_n$.  As a consequence, $\tau_n$ cannot be entirely included in this domain, so that the set $\wp_n\cap \tau_n$ is not empty. We take $z_n \in \wp_n\cap \tau_n$, and, up to further extraction, we suppose that $z_n \to z$. On the one hand, $\delta_{(n)}(y_n,z_n) \le \diam(\wp_n)$, so that $y \sim_\infty z$. On the other hand, $z \in \tau$ and $y \notin \tau$, so that $y \neq z$. Because~$y$ is not a leaf, this contradicts Theorem~\ref{ipb}.

The case $y = w^\bullet$ is treated with a slightly different argument. As the argument is exactly the same as in~\cite{bettinelli10tsl}, we do not treat it here.
\end{pre}

We now turn to general loops that are not necessarily made of edges. Here again, we use an argument similar to the one used in~\cite{miermont08sphericity,bettinelli10tsl}, with some minor changes. We fix $\eps>0$, and we let~$\eta$ be as in Lemma~\ref{1regb}. For~$k$ sufficiently large, the conclusion of Lemma~\ref{1regb} holds, together with the inequality $\eta\, \gamma \, n_k^{1/4} \ge 12$.

We call \textbf{pane} of $\cS_n$ the projection in $\cS_n$ of a $[z_{j-1},z_j]\times [0,1] \subseteq X_{f_*}$ for some $1 \le j \le 2\sigma_n$, with the notation of Section~\ref{secrepsurf}. We also call \textbf{semi-edge} the projection in $\cS_n$ of either $\{z_j\}\times [0,1] \subseteq X_{f_*}$ or $[z_{j-1},z_j]\times\{1\} \subseteq X_{f_*}$ for some $1 \le j \le 2\sigma_n$. These correspond to the edges of the prism $X_{f_*}$ that are not already edges in $\cS_n$. Let us consider a loop $\loo$ drawn in $\cS_{n_k}$ with diameter less than $\eta/2$. Consider the union of the closed internal faces\footnote{We call \textbf{closed face} the closure of a face.} and panes visited by $\loo$. The boundary of this union consists in simple loops made of edges and semi-edges in $\cS_{n_k}$. It should be clear that one of these loops entirely contains~$\loo$ in the closure of its inner domain. Let us denote this loop by~$\lambda$.

Let $\tilde\lambda$ be the largest (in the sense of the inclusion of the inner domains) simple loop made of edges contained in the closure of the inner domain of~$\lambda$ (that is, the loop obtained by removing the semi-edges of the form $\{z_j\}\times [0,1]$ and changing the ones of the form $[z_{j-1},z_j]\times\{1\}$ by $[z_{j-1},z_j]\times\{0\}$). Because every internal face and every pane of $\cS_{n_k}$ has diameter less than $3/(\gamma \, n_k^{1/4})$, we see that $\diam(\tilde\lambda) \le \diam(\loo)+ 6/(\gamma \, n_k^{1/4}) \le \eta$. Then, by Lemma~\ref{1regb}, the diameter of the inner domain of~$\tilde\lambda$ is less than~$\eps$. As a result, the diameter of the inner domain of~$\lambda$ is less than~$2\eps$, so that~$\loo$ is homotopic to~$0$ in its $2\eps$-neighborhood.

\subsection{Boundary of \hr{$\qis$}{q}}\label{secboundary}

We use the observation following Proposition~\ref{whyburn} to show that the boundary of~$\qis$ is (the image in~$\qis$ of) the floor~$\fli$ of~$\Fi$, and then give a lower bound on its Hausdorff dimension. We postpone the proof of the upper bound to Section~\ref{secboundaryub}, because we will need the notation of Section~\ref{secsnake}.

\begin{prop}\label{propboundary}
The boundary of~$\qis$ is given by $\partial\qis=\pii(\fli)$.
\end{prop}

\begin{pre}
We define a pseudo-metric $\tilde d_{GH}$ on the set of triples $(\X,\delta,A)$ where $(\X,\delta)$ is a compact metric space and $A\subseteq \X$ is a closed subset of~$\X$ by
$$\tilde d_{GH}\big((\X,\delta,A),(\X',\delta',A')\big) \de \inf \big\{ \dish\big(\varphi(\X),\varphi'(\X')\big) \vee \dish\big(\varphi(A),\varphi'(A')\big)\big\},$$
where the infimum is taken over all isometric embeddings $\varphi : \X \to \X''$ and $\varphi':\X'\to \X''$ of~$\X$ and~$\X'$ into the same metric space $(\X'', \delta'')$. By slightly adapting the proof of \cite[Theorem~7.3.30]{burago01cmg}, we can show that $\tilde d_{GH}((\X,\delta,A),(\X',\delta',A'))=0$ if and only if there is an isometry from $(\X,\delta)$ onto $(\X',\delta')$ whose restriction to~$A$ maps~$A$ onto~$A'$.

\bigskip
	
We proceed in three steps. First, note that the observation following Proposition~\ref{whyburn} implies that
\begin{equation}\label{eq1}
\tilde d_{GH}\big( (\cS_{n_k},\delta_{(n_k)},\partial\cS_{n_k}),(\qis,\disig,\partial\qis)\big) \tok 0.
\end{equation}

Secondly, we show that
\begin{equation}\label{eq2}
\tilde d_{GH}\big( (\cS_{n},\delta_{(n)},\partial\cS_{n}),(V(\q_n)\bs\{v_n^\ooo\},\delta_{(n)},\fl_n\bs\{v_n^\circ\})\big) \ton 0,
\end{equation}
where~$v_n^\circ$ is the extra vertex of the floor added when performing the \BDG bijection. We work here in~$\cS_n$ and see $V(\q_n)\bs\{v_n^\ooo\}$ as one of its subsets. Because of the way~$\cS_n$ is constructed, we see that $\dish(\cS_n,V(\q_n)\bs\{v_n^\ooo\})\le{3}/{\g}$. Using the technique we used in the proof of Lemma~\ref{lem0reg} to approach the points of~$\partial\cS_n$ by points lying in~$\fl_n\bs\{v_n^\circ\}$, and the fact that every point in~$\fl_n\bs\{v_n^\circ\}$ is at distance at most $1/(\g)$ from~$\partial\cS_n$, we obtain that
$$\dish(\partial\cS_n,\fl_n\bs\{v_n^\circ\}) \le \frac{3}{\g} + \omega_{\bb_{(n)}}(\eta),$$
as soon as $n \ge 1/2\eta^2$. As a result, $\limsup \tilde d_{GH}( (\cS_{n},\delta_{(n)},\partial\cS_{n}),(V(\q_n)\bs\{v_n^\ooo\},\delta_{(n)},\fl_n\bs\{v_n^\circ\})) \le \omega_{\bb_{\infty}}(\eta)$ for all $\eta >0$, and~\eqref{eq2} follows by letting $\eta \to 0$.

Finally, we see that
\begin{equation}\label{eq3}
\tilde d_{GH}\big( (V(\q_{n_k})\bs\{v_{n_k}^\ooo\},\delta_{(n_k)},\fl_{n_k}\bs\{v_{n_k}^\circ\}),(\qis,\disig,\pii(\fli))\big) \tok 0.
\end{equation}
Recall that $(V(\q_{n})\bs\{v_{n}^\ooo\},\delta_{(n)})$ is isometric to the space $(\QQ_n,d_{(n)})$ defined in Section~\ref{secprob}. We slightly abuse notation and view $\fl_n\bs\{v_n^\circ\}$ as a subset of~$\QQ_n$. We set $r_n\de \dis(\rR_n)/2$, where~$\rR_n$ is the correspondence between~$\QQ_n$ and~$\qis$ defined during Section~\ref{secprob}, and we define the pseudo-metric~$\Delta_n$ on the disjoint union $\QQ_n \sqcup \qis$ by $\Delta_n(x,y)\de d_{(n)}(x,y)$ if $x$, $y \in \QQ_n$, $\Delta_n(x,y)\de \disig(x,y)$ if $x$, $y \in \qis$, 
$$\Delta_n(x,y)\de \inf\{d_{(n)}(x,x') + r_n + \disig(y',y)\, :\, (x',y')\in \rR_n\}$$
if $x \in \QQ_n$ and $y \in \qis$, and $\Delta_n(x,y)\de \Delta_n(y,x)$ if $x \in \qis$ and $y \in \QQ_n$. It is a simple exercise to verify that~$\Delta_n$ is indeed a pseudo-metric and that $\dish(\QQ_n,\qis) \le r_n$. We showed in Section~\ref{secprob} that~$r_{n_k} \to 0$ as $k \to 0$, so that it is sufficient to prove that $\dish(\fl_{n_k},\pii(\fli)) \to 0$ as well. Let us argue by contradiction and suppose that this is not the case. There exists $\eps >0$ such that one of the following occurs:
\begin{enumerate}[($i$)]
	\item for infinitely many $n$'s, we can find a point~$t_n$ in the set $(2n+\sigma_n-1)^{-1} \, \ent 0 {2n+\sigma_n-1}$ such that $\bmp_{(n)}(t_n)\in \fl_n\bs\{v_n^\circ\}$, and $\Delta_n(\bmp_{(n)}(t_n),\pii(\fli)) \ge\eps$,
	\item for infinitely many $n$'s, there is $s_n\in[0,1]$ such that $\Fi(s_n) \in \fli$, and $\Delta_n(\qis(s_n),\fl_n\bs\{v_n^\circ\}) \ge\eps$.
\end{enumerate}
In the first case, up to extraction, we may suppose that $t_n \to t$. The fact that $\bmp_{(n)}(t_n)\in \fl_n\bs\{v_n^\circ\}$ yields that $C_{(n)}(t_n)= \underline C_{(n)}(t_n)$, so that $C_\infty(t)= \underline C_\infty(t)$ by continuity, and $\Fi(t) \in \fli$. We then have
$$\eps \le \Delta_n\big(\bmp_{(n)}(t_n),\pii(\fli)\big) \le \Delta_n\big(\bmp_{(n)}(t_n),\qis(t)\big) \le \disig(t_n,t) + r_n \to 0$$
along some subsequence. This is a contradiction. In the second case, we may also suppose that $s_n \to s$, and we have $\Fi(s)\in \fli$. We set $t_n \de \inf\{ t\, : \, C_{(n)}(t)=\underline C_{(n)}(s) \}$, so that $\bmp_{(n)}(t_n) \in \fl_n\bs\{v_n^\circ\}$. Up to further extraction, we have that $t_n \to t$, and because $C_\infty(t)=\ul C_\infty(s) =C_\infty(s)$, we see that $s \eqt t$, which yields $\disig(s,t)=0$. Finally,
$$\eps \le \Delta_n(\qis(s_n),\fl_n\bs\{v_n^\circ\}) \le \Delta_n\big(\qis(s_n),\bmp_{(n)}(t_n)\big) \le \disig(s_n,t_n) + r_n\to 0$$
along some subsequence.

\bigskip

Now, \eqref{eq1}, \eqref{eq2}, and~\eqref{eq3} yield that $\tilde d_{GH}\big((\qis,\disig,\partial\qis),(\qis,\disig,\pii(\fli))\big) = 0$, so that there exists an isometry $\varphi : \qis \to \qis$ such that $\pii(\fli) = \varphi(\partial\qis) = \partial(\varphi(\qis)) = \partial\qis$.
\end{pre}

We are now able to bound from below the Hausdorff dimension of~$\partial\qis$. We start with a lemma.

\begin{lem}\label{lemborninfl}
For~$a$, $b\in \fli$, we have
$$\disig(a,b) \ge \Li(a) - \max \lp \min_{\lhb a, b \rhb} \Li , \min_{\lhb b, a \rhb} \Li \rp.$$
\end{lem}

\begin{pre}
Let~$a_n$, $b_n\in \fl_n$ be points converging to~$a$ and~$b$, and let~$\wp_n$ be a geodesic from~$a_n$ to~$b_n$. Reasoning as in the beginning of the proof of Lemma~\ref{lemid}, we see that~$\wp_n$ either overflies $\lhb a_n, b_n \rhb$ for infinitely many~$n$'s, or it overflies $\lhb b_n, a_n \rhb$ for infinitely many~$n$'s. 

In the first case, let $c \in \lhb a, b \rhb$, and let $c_n \in \lhb a_n, b_n \rhb$ be a point converging toward~$c$. For the values of~$n$ for which~$\wp_n$ overflies $\lhb a_n, b_n \rhb$, we obtain by the remark of Section~\ref{secof}, and the triangle inequality, that
$$\Lab_n(c_n) \ge \Lab_n(a_n) - d_{\q_n}(a_n,b_n).$$
Taking the limit after renormalization along these values of~$n$, we obtain that $\Li(c) \ge \Li(a) - \disig(a,b)$. Taking the infimum for~$c$ over $\lhb a, b\rhb$, we find $\disig(a,b) \ge \Li(a) - \min_{\lhb a, b \rhb} \Li$. In the second case, a similar reasoning yields that $\disig(a,b) \ge \Li(a) - \min_{\lhb b,a \rhb} \Li$.
\end{pre}

\begin{pre}[Proof of Theorem~\ref{thmdimh} (lower bound)]
Recall that, for $x \in [0,\sigma]$, we defined $T_x \de \inf \{ r \ge 0 \, :\, C_\infty(r)= \sigma - x \}$. We also set $\fl(x) \de \qis(T_x)$, so that $\pii(\fli) = \{\fl(x),\, 0 \le x \le \sigma \}$.

To obtain the lower bound, we proceed as follows. We define the measure $\Lambda_{\fl}$ on $\qis$ supported by $\pii(\fli)$ as the image of the Lebesgue measure on $[0,\sigma]$ by the map $y \in [0, \sigma ] \mapsto \fl(y)$. Let us fix $x \in [0,\sigma]$. Because the process $y \in [0,\sigma] \mapsto \Li(T_y)=\bb_\infty(y)$ has the law of a Brownian bridge (up to a factor $\sqrt 3$), the law of the iterated logarithm ensures us that, a.s., for $\eta>0$, and~$\delta$ small enough,
\begin{equation}\label{logit}
\Li(T_x) - \min_{y\in \overrightarrow{[x-\delta^{2-\eta},x]}} \Li(T_y) > \delta \sand \Li(T_x) - \min_{y\in \overrightarrow{[x,x+\delta^{2-\eta}]}} \Li(T_y) > \delta.
\end{equation}
For $a \in \qis$ and $r>0$, we denote by $B_\infty(a,r)\subseteq \qis$ the open ball centered at~$a$ with radius~$r$ for the metric~$\disig$. Using Lemma~\ref{lemborninfl}, we see that, whenever~\eqref{logit} holds,
$$B_\infty(\fl(x),\delta)\cap \pii( \fli) \subseteq \fl\big((x-\delta^{2-\eta}, x+\delta^{2-\eta})\big),$$
so that $\Lambda_{\fl}(B_\infty(\fl(x),\delta)) \le 2\delta^{2-\eta}$. Finally, we obtain that, a.s., for all $a \in \pii(\fli)$,
$$\limsup_{\delta \to 0} \frac{\Lambda_{\fl}(B_\infty(a,\delta))}{\delta^{2-\eta}} \le 2.$$
We then conclude that $\dH(\qis, \disig) \ge 2 - \eta$ for all $\eta >0$ by standard density theorems for Hausdorff measures (\cite[Theorem 2.10.19]{federer69gmt}).
\end{pre}

\section{Singular cases}\label{singcases}

\subsection{Case \hr{$\sigma=0$}{sigma=0}}\label{sig0}

In the case $\sigma=0$, we could apply a reasoning similar to the one we used in Sections~\ref{secpt1} through~\ref{secsurfb}. We would obtain for the law of $(C_\infty,L_\infty)$ the law of a Brownian snake driven by a normalized Brownian excursion, and we would use a result of Whyburn \cite[Corollary~5.21 and Theorem~6.3]{whyburn35sls} treating the case where $\diam(\partial\X_n)\to 0$. Instead, we use a more direct approach, roughly consisting in saying that a uniform quadrangulation with ``small'' boundary is close to a uniform quadrangulation without boundary. A non-negligible advantage of this method is that it gives a more precise statement, Theorem~\ref{cvqbs0}, and completely identifies the limiting space as the Brownian map.

Let us begin with a lemma giving an upper bound on the Gromov--Hausdorff distance between a quadrangulation with a boundary and the quadrangulation obtained by applying Schaeffer's bijection to one of the trees of the forest that corresponds through the \BDG bijection.

\begin{lem}\label{lemqfqt}
Let $(\f,\lab) \in \fF_\sigma^n$ be a well-labeled forest, $\bb \in \Bs$ a bridge, $\tr$ a tree of~$\f$ rooted at~$\rho$, and $\bm b\in \{-1,0\}$. Then $(0,\bm b)\in \B_1$, and, up to a trivial transformation, $(\tr,\lab_{|\tr})$ may be seen as an element of~$\fF_1^{|\tr|-1}$. We denote by $\q_\f \in \Qns$ (resp.\ $\q_\tr \in \Q_{|\tr|-1,1}$) the quadrangulation corresponding to $((\f,\lab),\bb)$ [resp.\ to $((\tr,\lab_{|\tr}),(0,\bm b))$] through the \BDG bijection (we omit here the distinguished vertices). Then
$$\dGH\lp \big(\q_\f,d_{\q_\f}\big),\big(\q_\tr,d_{\q_\tr}\big)\rp \le 2 \lp \max_{\f\bs \ring\tr} \hat \lab - \min_{\f\bs \ring\tr} \hat \lab +1 \rp,$$
where $\ring\tr \de \tr\bs\{\rho\}$, and
$$\hat \lab(u)\de\lab(u)+\bb(\aaa(u)-1),\qquad u \in \f$$ 
is the labeling function of~$\f$, shifted tree by tree according to the bridge, as in Section~\ref{bdgfbq}.
\end{lem}

\begin{pre}
Before we begin, let us introduce some useful notation. For arcs $i_1 \arc i_2$, $i_2 \arc i_3$, \dots, $i_{r-1} \arc i_r$, we write
$$i_1 \arc i_2 \arc \dots \arc i_r$$
the path obtained by concatenating them. Consistently with Section~\ref{bdgqfb}, let~$v^\ooo$ be the extra vertex we add when performing the \BDG bijection and let $v^\circ\de (\sigma+1) \in \f$ be the last vertex of~$\f$. We will identify the sets $\tr\cup\{v^\ooo \}$ with $V(\q_\tr)$, as well as $(\f\bs\{v^\circ\}) \cup \{v^\ooo \}$ with $V(\q_\f)$. Then the set
$$\rR\de \big\{ (a,a)\, : \, a \in \tr\cup\{v^\ooo \} \big\} \cup \big\{ (a,\rho)\, :\, a\in \f\bs(\tr\cup\{v^\circ\}) \big\}$$
is a correspondence between~$\q_\f$ and~$\q_\tr$. Without loss of generality, we may suppose that~$\tr$ is the first tree of~$\f$. This yields in particular that an integer $i\in \ent 0 {2|\tr|-2}$ codes the same vertex in~$\tr$ and in~$\f$, namely $\tr(i)=\f(i)$. Because we will apply the \BDG bijection at the same time to both $((\f,\lab),\bb)$ and $((\tr,\lab_{|\tr}),(0,\bm b))$, we will write $\suc_\f(i)$ the successor of $i\in \ent 0 {2n+\sigma-1}$ in the forest~$\f$, and $\suc_\tr(i)$ the successor of $i\in \ent 0 {2|\tr|-2}$ in the tree~$\tr$, in order to avoid confusion. We also set $l_2 \de \max_{\f\bs \ring\tr} \hat \lab$ and $l_1 \de \min_{\f\bs \ring\tr} \hat \lab$ for more clarity. Using the characterization~\eqref{dghcorres} of the Gromov--Hausdorff distance via correspondences, we see that it suffices to show that, for all $(a,a')$, $(b,b') \in \rR$, we have
$$\lt d_{\q_\f}(a,b) - d_{\q_\tr}(a',b') \rt \le 4(l_2 - l_1 +1).$$

\paragraph{First case: $a$, $b \in \f\bs(\ring\tr\cup\{v^\circ\})$.}
In this case, either $[a,b]$ or $[b,a]$ entirely lies inside $\f\bs\ring\tr$. As a result, \eqref{dlemmeb} gives
$$\lt d_{\q_\f}(a,b) - d_{\q_\tr}(\rho,\rho) \rt \le \hat \lab (a) + \hat \lab (b) - 2 \min_{\f\bs\ring\tr} \hat \lab + 2 \le 2(l_2 - l_1 +1).$$

\paragraph{Second case: $a$, $b \in \tr\cup\{v^\ooo \}$.}
We may suppose $a\neq b$. We proceed in two steps. We first claim that
$$d_{\q_\tr}(a,b) \le d_{\q_\f}(a,b).$$
To see this, let $\wp = (\wp(0), \wp(1), \dots, \wp(k))$ be any path (not necessarily geodesic) between~$a$ and~$b$ in~$\q_\f$. We will construct a shorter path from~$a$ to~$b$ in~$\q_\tr$, and our claim will immediately follow. Our construction is based on the simple observation that, if an arc exists in~$\q_\f$ between two points of $\tr\cup\{v^\ooo \}$, then the same arc also exists in~$\q_\tr$. We then only have to replace the portions of~$\wp$ that ``exit'' $\tr\cup\{v^\ooo \}$ with shorter paths in~$\q_\tr$. Precisely, we can restrict ourselves to the case where $\wp(r) \in \f\bs(\tr\cup\{v^\ooo \})$ for $0 < r <k$, with $k \ge 2$. We will also need to observe that a path linking two vertices of label~$l$ and~$l'$ has length at least $|l-l'|$.

Let us denote by~$i$ the integer such that the arc $(\wp(0),\wp(1))$ is either $i \arcr \suc_\f(i)$ or $i \arcl \suc_\f(i)$. We will say that $(\wp(0),\wp(1))$ is oriented \textbf{to the right} in the first case, and \textbf{to the left} in the second case. We also define~$j$ in a similar way for the arc $(\wp(k),\wp(k-1))$. Four possibilities are then to be considered (see Figure~\ref{pathremp}):
\begin{itemize}
\item Both $(\wp(0),\wp(1))$ and $(\wp(k),\wp(k-1))$ are oriented to the right. Without loss of  generality, we may suppose $i < j$. Properties of the \BDG bijection then show that $\hat \lab(\f(j)) \ge \hat \lab(\f(i))$, and we have
$$k \ge 1 + |(\hat \lab(\f(j))-1) -(\hat \lab(\f(i))-1)|+1 = \hat \lab(\f(j)) - \hat \lab(\f(i)) +2.$$
The following path in~$\q_\tr$,
$$j \arc \suc_\tr(j) \arc \dots \arc \underbrace{\suc_\tr^{\hat \lab(\f(j)) -\hat \lab(\f(i)) +1} (j)}_{= \suc_\tr(i)} \arc i,$$
links~$a$ to~$b$ in~$\q_\tr$ and has length less than~$k$. The equality in the last line is an easy consequence of the \BDG construction.

\item Both $(\wp(0),\wp(1))$ and $(\wp(k),\wp(k-1))$ are oriented to the left. Here again, we may suppose $i <j$. In this case, $\hat \lab(\f(j)) > \hat \lab(\f(i))$, and
$$\suc_\f(j) \arc \suc_\tr(\suc_\f(j)) \arc \dots \arc \underbrace{\suc_\tr^{\hat \lab(\f(j)) -\hat \lab(\f(i))} (\suc_\f(j))}_{= \suc_\f(i)}$$
fulfills our requirements.

\item $(\wp(0),\wp(1))$ is oriented to the right, and $(\wp(k),\wp(k-1))$ is oriented to the left. Necessarily, we have $\suc_\f(j) < i$, or $\suc_\f(j)=\infty$. If $\hat \lab(\f(i)) \ge \hat\lab(\f(j))$, then we take
$$i \arc \suc_\tr(i) \arc \dots \arc \underbrace{\suc_\tr^{\hat \lab(\f(i)) -\hat \lab(\f(j)) +1} (i)}_{= \suc_\f(j)},$$
otherwise, we take
$$\suc_\f(j) \arc \suc_\tr(\suc_\f(j)) \arc \dots \arc \underbrace{\suc_\tr^{\hat \lab(\f(j)) -\hat \lab(\f(i))} (\suc_\f(j))}_{= \suc_\tr(i)} \arc i,$$

\item $(\wp(0),\wp(1))$ is oriented to the left, and $(\wp(k),\wp(k-1))$ is oriented to the right. By considering the path $\bar\wp\de (\wp(k),\wp(k-1),\dots,\wp(0))$ instead of~$\wp$, we are back to the previous case.
\end{itemize}

\begin{figure}[ht]
		\psfrag{i}[][]{$i$}
		\psfrag{a}[r][r]{$\suc_\f(i)$}
		\psfrag{j}[][]{$j$}
		\psfrag{z}[r][r]{$\suc_\f(j)$}
		\psfrag{t}[][]{$\tr$}
		\psfrag{p}[][]{\textcolor{red}{$\wp$}}
	\centering\includegraphics{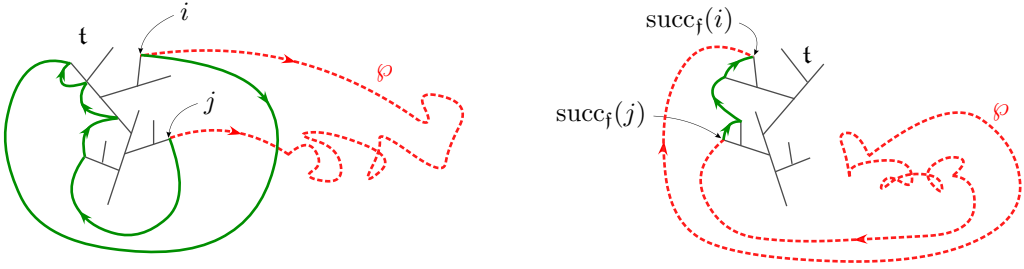}
	\caption[Gromov--Hausdorff distance between~$\q_\f$ and~$\q_\tr$.]{On this picture, $\tr$ is the only part of~$\f$ represented. The dashed (red) line represents the path~$\wp$ (in~$\q_\f$) and the (green) solid path is the path in~$\q_\tr$. Both first cases are represented.}
	\label{pathremp}
\end{figure}

\bigskip
We now show that
$$d_{\q_\f}(a,b) \le d_{\q_\tr}(a,b) + 2(l_2 - l_1 +1).$$
Let us consider a path~$\wp$ of length~$k$ in~$\q_\tr$ from~$a$ to~$b$. We are going to construct a path in~$\q_\f$ from~$a$ to~$b$, with length less than $k + 2(l_2 - l_1 +1)$. The only arcs present is~$\q_\tr$ but not in~$\q_\f$ are of the form $i \arc \suc_\tr(i)$ with $\suc_\tr(i) < i$ or $\suc_\tr(i)=\infty$, and $l_1 + 1 \le \hat \lab(\f(i)) \le l_2+1$. For convenience, let us call \textbf{pathological} such arcs. For all pathological arcs $i \arc \suc_\tr(i)$ and $j \arc \suc_\tr(j)$ with $i < j$, we can construct the path
\begin{equation}\label{path1}
j \arc \suc_\f(j) \arc \dots \arc \underbrace{\suc_\f^{\hat\lab (\f(j)) - \hat \lab (\f(i))+1}(j)}_{= \suc_f(i)} \arc i
\end{equation}
linking $\f(i)$ to $\f(j)$ in~$\q_\f$, its length being $\hat\lab (\f(j)) - \hat \lab (\f(i))+2$. We can also construct the path
\begin{equation}\label{path2}
j \arc \suc_\f(j) \arc \dots \arc \underbrace{\suc_\f^{\hat\lab(\f(j)) - l_1+1}(j)}_{=\suc_\f^{\hat\lab(\f(j)) - l_1}(\suc_\tr(j))} \arc \dots \arc \suc_\tr(j)
\end{equation}
linking $\f(j)$ to $\f(\suc_\tr(j))$ in~$\q_\f$, its length being $2(\hat\lab(\f(j)) - l_1)+1 \le 2(l_2-l_1+1)+1$. Using these paths, we construct our path in~$\q_\f$ as follows. If~$\wp$ does not use any pathological arcs, then~$\wp$ can be seen as a path in~$\q_\f$. If~$\wp$ uses exactly one pathological arc, we construct our new path by changing this arc into a path of the form~\eqref{path2}. By doing so, we obtain a path from~$a$ to~$b$ in~$\q_\f$ with length smaller than $k-1 + 2(l_2-l_1+1)+1$. Now, if~$\wp$ uses more than two pathological arcs, let $i \arc \suc_\tr(i)$ be the first one it uses, and $j \arc \suc_\tr(j)$ the last one. Let us denote by~$i_1$ and~$i_2$ the indices at which~$\wp$ uses them: $(\wp(i_1),\wp(i_1+1))=i \arcr \suc_\tr(i)$ or $i \arcl \suc_\tr(i)$ and $(\wp(i_2),\wp(i_2+1))=j \arcr \suc_\tr(j)$ or $j \arcl \suc_\tr(j)$. Changing~$\wp$ into its reverse~$\bar\wp$ if needed, we may suppose that $(\wp(i_1),\wp(i_1+1))=i \arcr \suc_\tr(i)$. If $(\wp(i_2),\wp(i_2+1))=j \arcl \suc_\tr(j)$, we change the portion $(\wp(i_1),\wp(i_1+1),\dots,\wp(i_2+1))$ into the path~\eqref{path1}, and obtain a new path shorter than~$\wp$. Finally, if $(\wp(i_2),\wp(i_2+1))=j \arcr \suc_\tr(j)$, we change the portion $(\wp(i_1),\wp(i_1+1),\dots,\wp(i_2+1))$ into the path~\eqref{path1} concatenated with the path~\eqref{path2}, and obtain a new path satisfying our requirements.

\paragraph{Third case: $a \in \tr\cup\{v^\ooo \}$, $b \in \f\bs(\tr\cup\{v^\circ\})$.}
We can write
\begin{align*}
\lt d_{\q_\f}(a,b) - d_{\q_\tr}(a,\rho) \rt &\le \lt d_{\q_\f}(a,b) - d_{\q_\f}(a,\rho) \rt + \lt d_{\q_\f}(a,\rho) - d_{\q_\tr}(a,\rho) \rt\\
									&\le d_{\q_\f}(b,\rho) + \lt d_{\q_\f}(a,\rho) - d_{\q_\tr}(a,\rho) \rt\\
									&\le 4(l_2 - l_1 +1),
\end{align*}
by applying the first case to $(b,\rho)$ and the second case to $(a,\rho)$. This ends the proof.
\end{pre}

We may now proceed to the proof of Theorem~\ref{cvqbs0}. We use the same notation as in Section~\ref{secccf}, and Corollary~\ref{cvln} remains true, if the process $(C_\infty,\Lab_\infty)$ has the law of a Brownian snake driven by a normalized Brownian excursion. As we will not need the explicit law of the process $(C_\infty,\Lab_\infty)$ in what follows, we do not prove this, and refer the reader to~\cite{bettinelli10slr}, in particular to Proposition~15 for similar results. By Skorokhod's representation theorem, we still assume that this convergence holds almost surely.

\begin{pre}[Proof of Theorem~\ref{cvqbs0}]
We define~$\tr_n$ as the largest tree of~$\f_n$ (if there are more than one largest tree, we take~$\tr_n$ according to some convention, for example the first one), and we consider a random variable $\bm b_n$ uniformly distributed over $\{-1,0\}$, independent of~$\q_n$. We denote by~$\hat \q_n$ the quadrangulation corresponding, as in the statement of Lemma~\ref{lemqfqt}, to $((\tr_n,\lab_{n\,|\tr_n}),(0,\bm b_n))$ through the \BDG bijection. Then, conditionally given $|\tr_n|=k+1$, the quadrangulation $\hat \q_n$ is uniformly distributed over the set~$\Q_{k,1}$.

From now on, we work on the set of full probability where the convergence $C_{(n)} \to C_\8$ holds. Let $\eps \in (0,1/4)$, and $2\eta \de \min_{[\eps, 1-\eps]} C_\8 >0$. As $C_{(n)}$ tends to $C_\8$, for~$n$ large enough, we have $\min_{[\eps, 1-\eps]} C_{(n)} \ge \eta$ and $\sigma_{(n)} < \eta$. As a result,
\begin{align*}
s_n&\de \inf\lb r \le \frac 1 2 \, : \, C_{(n)}(r)=\ul C_{(n)}\lp \frac 12 \rp \rb \le \eps,\\
t_n&\de \sup\lb r \ge \frac 1 2 \, : \, C_{(n)}(r)=\ul C_{(n)}\lp \frac 12 \rp \rb \ge1- \eps,
\end{align*}
so that $\tr_n$ is coded by $[(2n + \sigma_n -1) \, s_n, (2n + \sigma_n -1) \, t_n]$. Note that, in particular, this implies that $|\tr_n| \ge n \, (1-2\eps)$. This fact will be used later. By Lemma~\ref{lemqfqt},
\begin{multline*}
\limsup_{n\to\8} \dGH \lp \Big( V(\q_{n}),d_{\q_{n}}/(\gamma n^{1/4}) \Big),\Big( V(\hat\q_{n}), d_{\hat\q_{n}}/(\gamma n^{1/4}) \Big) \rp \hspace{-20mm} \\
	\le 2 \limsup_{n\to \8} \Bigg( \sup_{\overrightarrow{[1-\eps,\eps]}} \Lab_{(n)} - \inf_{\overrightarrow{[1-\eps,\eps]}} \Lab_{(n)}  + \frac{1}{\g}  \Bigg)= 2 \Bigg( \sup_{\overrightarrow{[1-\eps,\eps]}} \Lab_{\8} - \inf_{\overrightarrow{[1-\eps,\eps]}} \Lab_{\8} \Bigg) \tode{}{\eps}{0} 0.
\end{multline*}
Let us set $\hat\delta_{(n,k)} \de d_{\hat\q_{n}}/(\gamma k^{1/4})$. We then have to see that $(V(\hat\q_n),\hat\delta_{(n,n)})$ converges toward the Brownian map $(\m_\8,D^*)$. Let $f:\MM \to \R$ be uniformly continuous and bounded.
By the Portmanteau theorem \cite[Theorem 2.1]{billingsley68cpm}, we only need to show that
$$\E{f\big( V(\hat\q_n),\hat\delta_{(n,n)} \big)} \ton \E{f\big( \m_\8,D^* \big)}.$$

Let $\eps >0$. If we delete from~$\hat \q_n$ the only edge on the boundary that is not the root, we obtain a quadrangulation without boundary, which, conditionally given $|\tr_n|=k+1$, is uniformly distributed over the set of planar quadrangulations with~$k$ faces. As this operation does not affect the underlying metric space, by \cite[Theorem~1]{miermont11bms} or \cite[Theorem~1.1]{legall11ubm}, we obtain that the distribution of $(V(\hat\q_n),\hat\delta_{(n,k)})$ conditioned on $|\tr_n|=k+1$ converges toward the distribution of $(\m_\8,D^*)$ as $k\to \8$. As $(\m_\8,D^*)$ is a compact metric space, we can find large~$n_0$ and~$M$ such that, for all $k \ge n_0/2$ and~$n$ for which $\Pb(|\tr_n|=k+1) >0$,
\begin{equation}\label{sig0eq1}
\Pb\lp \diam\big( V(\hat\q_n),\hat\delta_{(n,k)} \big) \ge M \ \big|\ |\tr_n|=k+1 \rp <\frac{\eps}{2\, \sup f},
\end{equation}
and
\begin{equation}\label{sig0eq2}
\lt \,\E{f \big( V(\hat\q_n),\hat\delta_{(n,k)} \big)\ \big|\ |\tr_n|=k+1} - \E{f\big( \m_\8,D^* \big)} \rt < \eps.
\end{equation}
We then choose~$\eta\in (0,1/2)$ such that, for all $(\X,\delta)$, $(\X',\delta')$, 
\begin{equation}\label{sig0eq3}
\dGH\big((\X,\delta),(\X',\delta')\big) \le \frac12\, M\,\lp 1-(1-\eta)^{1/4}\rp \quad \Rightarrow \quad
\lt f\big((\X,\delta)\big) - f\big((\X',\delta')\big) \rt < \eps.
\end{equation}
For $n\ge n_0$, we then have
\begin{multline*}
\lt\, \E{f\big( V(\hat\q_n),\hat\delta_{(n,n)} \big)} - \E{f\big( \m_\8,D^* \big)} \rt \le 2\sup f\ \Pb\big(|\tr_n| \le n\, (1-\eta)\big) \\
+ \sum_{k=\lc n\,(1-\eta)\rc}^n \Pb(|\tr_n|=k+1)\ \lt \,\E{f \big( V(\hat\q_n), \hat\delta_{(n,n)} \big)\ \big|\ |\tr_n|=k+1} - \E{f\big( \m_\8,D^* \big)} \rt.
\end{multline*}
By the observation we previously made, we see that the first term in the right-hand side tends to~$0$ as $n\to \8$. To conclude, it will be sufficient to show that the term between vertical bars in the sum is smaller than~$3\eps$. Using~\eqref{sig0eq1}, \eqref{sig0eq2}, and the fact that $n\, (1-\eta) \ge n_0/2$, we obtain that it is smaller than
$$2\eps + \E{\lp f \big( V(\hat\q_n), \hat\delta_{(n,n)} \big) - f \big( V(\hat\q_n), \hat\delta_{(n,k)} \big)\rp\un{\{\diam( V(\hat\q_n),\hat\delta_{(n,k)}) < M \}} \ \big|\ |\tr_n|=k+1}.$$
By taking a trivial correspondence between $( V(\hat\q_n), \hat\delta_{(n,n)})$ and $( V(\hat\q_n), \hat\delta_{(n,k)})$, it is not hard to see that the Gromov--Hausdorff distance between these two spaces is smaller than
$$\frac{1}{2}\,\diam\big( V(\hat\q_n),\hat\delta_{(n,k)}\big) \lp 1 - \lp {k}/{n}\rp^{1/4}\rp.$$
We finally obtain the desired bound thanks to~\eqref{sig0eq3}.
\end{pre}

\subsection{Case \hr{$\sigma=\infty$}{sigma=infinity}}\label{sig8}

In this case, the scaling factor changes. We use the same formalism as in the beginning of Section~\ref{secccf}, except that we now suppose that the sequence $(\sigma_n)_{n\ge 1}$ satisfies $\sigma_n / \sqrt{2n} \to \8$ as $n\to \infty$. By \cite[Lemma~10]{bettinelli10slr}, the process
\begin{equation}\label{bnsn}
\lp \frac{\bb_n(\sigma_n\, s)}{(2\sigma_n)^{1/2}} \rp_{0\le s \le 1}
\end{equation}
converges in distribution toward a standard Brownian bridge $\fB=B_{[0,1]}^{0\to 0}$. By Skorokhod's representation theorem, we will assume that this convergence holds almost surely. We define on $[0,1]$ the pseudo-metric
$$\delta_\fB(s,t)\de \fB(s) + \fB(t) - 2 \max \lp \min_{ r \in \overrightarrow{[s,t]}} \fB(r),\min_{r \in \overrightarrow{[t,s]}} \fB(r) \rp, \quad 0 \le s,t \le 1.$$
By Vervaat's transformation \cite[Theorem~1]{vervaat79rbb}, the metric space $(\TT_\fB\de [0,1]_{/\{\delta_\fB=0\}},\delta_\fB)$ is isometric to the CRT $(\TT_\ee,\delta_\ee)$. We will show the convergence toward this space, by using correspondences.

\begin{pre}[Proof of Theorem~\ref{cvqbs8}]
We denote by $\bmp_\fB : [0,1]\to \TT_\fB$ the canonical projection, and we define the correspondence~$\fR_n$ between $(V(\q_n)\bs\{v_n^\ooo\},(2\sigma_n)^{-1/2} d_{\q_n})$ and $(\TT_\fB,\delta_\fB)$ by
$$\fR_n \de \big\{ \big(\q_n(i),\bmp_\fB(s)\big)\, :\, i \in \ent 0 {2n+\sigma_n-1},\, s\in [0,1],\, \sigma_n-\ul C_n(i)=\lf \sigma_n \, s\rf \big\}.$$
In terms of forests, this roughly consists in saying that all the vertices of a tree are in correspondence with a small segment corresponding to the edge of the floor following the root of the tree. It is sufficient to show that the distortion of~$\fR_n$ tends to~$0$ as $n\to\8$. For any two points~$a$ and $b\in \f_n$, we have the following bounds:
$$\Lab_n(a) + \Lab_n(b) - 2 \max \lp \min_{ \lhb a,b \rhb} \Lab_n,\min_{\lhb b,a \rhb} \Lab_n \rp \le d_{\q_n}(a,b) \le \Lab_n(a) + \Lab_n(b) - 2 \max \lp \min_{ [a,b]} \Lab_n,\min_{[b,a]} \Lab_n \rp + 2.$$
The second inequality is merely the bound~\eqref{dlemmeb}, and the first one is easily obtained by a technique similar to the one we used in the proof of Lemma~\ref{lemborninfl} (see also \cite[Lemma~20]{bettinelli10slr} or~\cite{curien11cactus}). It is thus easy to see (recall the definition~\eqref{Labn} of~$\Lab_n$) that, for~$i$, $j\in \ent 0 {2n+\sigma_n-1}$,
$$\lt \ d_{\q_n}\big(\q_n(i),\q_n(j)\big) -\lp  \bb_n(u) + \bb_n(v) - 2 \max \lp \min_{ \overrightarrow{\ent u v}} \bb_n,\min_{\overrightarrow{\ent v u}} \bb_n \rp \rp \rt \le 3\,(\sup \lab_n - \inf \lab_n) + 2,$$
where we wrote $u \de \sigma_n-\ul C_n(i)$ and $v\de\sigma_n-\ul C_n(j)$. Using the convergence of the process~\eqref{bnsn} stated before, we obtain
$$\limsup_{n\to\8} \dis (\fR_n) \le \limsup_{n\to\8} \frac{3\,(\sup \lab_n - \inf \lab_n)}{(2\sigma_n)^{1/2}}.$$
It remains to show that the latter quantity is equal to~$0$ in probability. This is a consequence of Lemma~\ref{cnln0}, which follows.
\end{pre}

We still denote by $(C_n, L_n)$ the contour pair of $(\f_n,\lab_n)$, but we now define the scaled versions of~$C_n$ and $L_n$ by
$$C_{[n]} \de \lp \frac {C_n(k_n \, s)} {\sigma_n} \rp_{0\le s \le 1}
\sand L_{[n]} \de \lp \frac {L_n(k_n\,s)} {\sqrt{\sigma_n}}\rp_{0\le s \le 1},$$
where we wrote $k_n \de 2n +\sigma_n$.

\begin{lem}\label{cnln0}
The pair $(C_{[n]},L_{[n]})$ converges toward $\big((1-s)_{0\le s \le 1}, (0)_{0\le s \le 1}\big)$ in distribution in the space $(\K,d_\K)^2$.
\end{lem}

\begin{pre}
The first step consists in showing the convergence of the first component
$$C_{[n]} \to (1-s)_{0\le s \le 1}.$$
At first, we will consider bridges instead of first-passage bridges.

\paragraph{Step 1.} Let $(S_i)_{i\ge 0}$ be a simple random walk started at~$0$, and, for all $p\in [0,1]$, let $(S^{(p)}_i)_{i\ge 0}$ be a random walk started at~$0$ with steps having the distribution $p \, \delta_{2(1-p)} + (1-p) \,\delta_{-2p}$. It is a simple computation to see that, for any measurable function~$f$ and any~$k$,
$$\E{f\big( (S_i)_{0\le i \le k} \big)}=\E{\big(4p(1-p)\big)^{-k/2} \lp\frac{p}{1-p}\rp^{-(S^{(p)}_k+k(2p-1))/2} f\big( (S^{(p)}_i+i(2p-1))_{0\le i \le k} \big)}.$$
(Note that $(S^{(p)}_i+i(2p-1))_{i\ge 0}$ is a random walk whose steps have the distribution $p \, \delta_{1} + (1-p) \,\delta_{-1}$.) Let us fix $n\in \N$ and $\eps >0$. Applying the latter equality, we obtain that
\begin{equation}\label{pbskn0}
\Pb\lp\sup_{0\le i \le k_n} \lt S_i + i\frac{\sigma_n}{k_n} \rt > \eps \sigma_n \ \bigg| \ S_{k_n}=-\sigma_n\rp=\Pb\lp\sup_{0\le i \le k_n} \lt S^{(p_n)}_i \rt > \eps \sigma_n \ \bigg| \ S^{(p_n)}_{k_n}=0\rp,
\end{equation}
if we choose $p_n\de 1/2-\sigma_n/2k_n$.

For $m\in \Z$, it should be clear that, under $\Pb( \cdot \ | \ S^{(p_n)}_{k_n}=2m)$, the path $(S^{(p_n)}_i)_{0\le i \le k_n}$ is uniformly distributed among the paths going from~$0$ to~$2m$ and having steps with value $2(1-p_n)$ or $-2p_n$. Then, changing uniformly a $-2p_n$-step into a $2(1-p_n)$-step, we obtain a path with law $\Pb( \cdot \ | \ S^{(p_n)}_{k_n}=2(m+1))$ that always lies above the previous one. This observation shows the stochastic domination 
$$\Pb\lp \cdot \ \big| \ S^{(p_n)}_{k_n}=2m\rp \preceq \Pb\lp \cdot \ \big| \ S^{(p_n)}_{k_n}=2(m+1)\rp,$$
from which we obtain that
\begin{align*}
\Pb\lp S^{(p_n)}_{k_n}\ge 0 \rp\ \Pb\lp\sup_{0\le i \le k_n} S^{(p_n)}_i > \eps \sigma_n \ \bigg| \ S^{(p_n)}_{k_n}=0\rp\hspace{-20mm}\\
		&=   \sum_{m=0}^{\8} \Pb\lp S^{(p_n)}_{k_n}=2m\rp \ \Pb\lp\sup_{0\le i \le k_n} S^{(p_n)}_i > \eps \sigma_n \ \bigg| \ S^{(p_n)}_{k_n}=0\rp \\
		&\le \sum_{m=0}^{\8} \Pb\lp S^{(p_n)}_{k_n}=2m\rp \ \Pb\lp\sup_{0\le i \le k_n} S^{(p_n)}_i > \eps \sigma_n \ \bigg| \ S^{(p_n)}_{k_n}=2m\rp \\
		&\le \Pb\lp\sup_{0\le i \le k_n} S^{(p_n)}_i > \eps \sigma_n\rp.
\end{align*}
The term $\Pb( S^{(p_n)}_{k_n}\ge 0 )$ is equal to $\Pb(\cB(k_n,p_n) \ge k_n p_n)$, where $\cB(k_n,p_n)\de S^{(p_n)}_{k_n}/2 + k_n p_n$ has a binomial distribution with parameters~$k_n$ and~$p_n$. By \cite[Theorem~2]{hamza95}, this quantity is larger than $1/2$. Adding to this the fact that $(S^{(p_n)}_i)_{i\ge 0}$ is a martingale, we obtain, by applying Doob's inequality, that
$$\Pb\lp\sup_{0\le i \le k_n} S^{(p_n)}_i > \eps \sigma_n \ \bigg| \ S^{(p_n)}_{k_n}=0\rp \le \frac{2}{\eps^2 \sigma_n^2}\ \E{\lp S^{(p_n)}_{k_n}\rp^2}=\frac{8p_n(1-p_n)k_n}{\eps^2 \sigma_n^2}\le\frac{2k_n}{\eps^2 \sigma_n^2}.$$
Using a similar argument to bound $\Pb(\inf_{0\le i \le k_n} S^{(p_n)}_i < -\eps \sigma_n \ | \ S^{(p_n)}_{k_n}=0)$, we see that the quantity \eqref{pbskn0} is smaller than $4k_n/\eps^2 \sigma_n^2$.

\bigskip

Finally, the construction of discrete first-passage bridges from discrete bridges provided in \cite[Theorem~1]{bertoin03ptf} yields
\begin{align*}
\Pb\lp \sup_{0\le s\le 1} \lt C_{[n]}(s)-(1-s) \rt > \eps \rp
&=\Pb\lp \sup_{0\le i \le k_n} \lt C_{n}(i)-\sigma_n+ i\frac{\sigma_n}{k_n} \rt > \eps \sigma_n\rp\\
&\le \Pb\lp\sup_{0\le i \le k_n} \lt S_i + i\frac{\sigma_n}{k_n} \rt > \frac{\eps}{2}\, \sigma_n \ \bigg| \ S_{k_n}=-\sigma_n\rp\\
&\le \frac{16k_n}{\eps^2 \sigma_n^2} = \frac{16}{\eps^2} \lp\frac{2n}{\sigma_n^2}+\frac{1}{\sigma_n}\rp \ton 0.
\end{align*}

\paragraph{Step 2.} Now that we have the convergence of the first component, let us prove the convergence of the pair $(C_{[n]},L_{[n]})$. As explained in the proof of \cite[Proposition~15]{bettinelli10slr}, it is sufficient to show that, for every $q \ge 2$, there exists a constant $K_q$ satisfying, for all~$n$ and all $0 \le s \le t \le 1/2$ for which $k_ns$ and $k_n t$ are integers,
$$\E{\lt S_{k_n t} - S_{k_n s} \rt^{q}\ \big| \ S_{k_n}=-\sigma_n} \le K_q \,\sigma_n^q \lt t -s \rt^{q/ 2}.$$
Using the same method as above (with the same value of~$p_n$), we see that the left-hand side is equal to
$$\E{\lt S^{(p_n)}_{k_n t} - S^{(p_n)}_{k_n s} - \sigma_n (t-s) \rt^{q}\ \big| \ S^{(p_n)}_{k_n}=0}.$$
We need to bound the quantity
\begin{equation}\label{eqquo}
\E{\lt S^{(p_n)}_{k_n t} - S^{(p_n)}_{k_n s} \rt^{q}\ \big| \ S^{(p_n)}_{k_n}=0} = \E{\lt S^{(p_n)}_{k_n t} - S^{(p_n)}_{k_n s} \rt^{q} \ \frac{Q^{S^{(p_n)}}_{k_n(1-t)}\lp S^{(p_n)}_{k_n t} \rp}{Q^{S^{(p_n)}}_{k_n}\lp 0 \rp}},
\end{equation}
where we used the notation $Q_a^{S^{(p_n)}}(b) \de \Pb(S_a^{(p_n)}=b)$ and the Markov property at time $k_n t$. Using the simple fact\footnote{Observe that, when $p\in (0,1)$, $\Pb(\cB(m,p)=r)\ge \Pb(\cB(m,p)=r-1)$ if and only if $r\le (m+1)p$.} that for a binomial variable $\cB(m,p)$ with parameters $m\in \N$ and $p\in [0,1)$, we have
$$\sup_{r\ge 0} \Pb\big(\cB(m,p)=r\big) = \Pb\big(\cB(m,p)=\lf (m+1)p \rf\big),$$
we see that the quotient in the right-hand side of \eqref{eqquo} is smaller than
$$\frac{\Pb\big( \cB(k_n(1-t),p_n)=\lf (k_n(1-t)+1)p_n \rf \big)}{\Pb\big( \cB(k_n,p_n)=k_n p_n \big)}\mathop{\sim}\limits_{n\to \8} \frac{1}{\sqrt{1-t}}\le \sqrt 2,$$
so that it is uniformly bounded in~$n$ by some finite constant~$K$. Finally, we conclude thanks to Rosenthal's Inequality \cite[Theorem 2.9 and 2.10]{petrov95ltp} that there exists a constant (depending on~$q$) $K'_q$ such that \eqref{eqquo} is smaller than
$$K\ \E{\lt S^{(p_n)}_{k_n t} - S^{(p_n)}_{k_n s} \rt^{q}} \le K'_q\ \E{\lt S^{(p_n)}_1 \rt^{q}} \, k_n^{q/2} \lt t -s \rt^{q/ 2} \le (K_q-1) \,\sigma_n^q \lt t -s \rt^{q/ 2},$$
with $K_q\de K'_q\, 2^q \sup_n (k_n/\sigma_n^2)^{q/2} +1< \8$. This completes the proof.
\end{pre}

\section{Proofs using the Brownian snake}\label{secsnake}

In this section, we prove Lemmas~\ref{minb}, \ref{pcb}, \ref{legalllemb}, and complete the proof of Theorem~\ref{thmdimh}. To this end, we will need some notions about the Brownian snake. We refer the reader to~\cite{legall99sbp} for a complete description of this object. Recall that we denoted by~$\K$ the space of continuous real-valued functions on $\R_+$ killed at some time, and that we wrote $\zeta(w)$ the lifetime of an element $w \in \K$. We also use the notation $\widehat w \de w(\zeta(w))$ for the final value of a path~$w\in \K$. From now on, we will work on the space $\Omega' \de \C(\R_+, \K)$ of continuous functions from~$\R_+$ into~$\K$, equipped with the topology of uniform convergence on every compact subset of~$\R_+$. We write $W_s \de \omega(s)$ the canonical process on~$\Omega'$, and denote by $\zeta_s\de \zeta(W_s)$ its lifetime.

For $w\in \K$, we denote the law of the Brownian snake started from~$w$ by~$\Pb_w$. This means that, under~$\Pb_w$, the process $(\zeta_s)_{s\ge 0}$ has the law of a reflected Brownian motion on~$\R_+$ started from~$\zeta(w)$, and that the conditional distribution of $(W_s)_{s \ge 0}$ knowing $(\zeta_s)_{s\ge 0}$, denoted by~$\Theta_w^\zeta$, is characterized by
\begin{itemize}
 \item $W_0 = w$, $\Theta_w^\zeta$ a.s.
 \item the process $(W_s)_{s\ge 0}$ is time-inhomogeneous Markov under~$\Theta_w^\zeta$ and, for $0\le s \le s'$,
	\begin{itemize}
	 \item $W_{s'}(t) = W_s(t)$ for all $0\le t \le \zeta_r$, $\Theta_w^\zeta$ a.s., where $\zeta_r \de \inf_{r'\in[s,s']} \zeta_{r'}$,
	 \item under  $\Theta_w^\zeta$, the process $\lp W_{s'}(\zeta_r + t) \rp_{0\le t \le \zeta_{s'} - \zeta_r}$ is independent of $W_s$ and distributed as a real Brownian motion started from $W_s(\zeta_r)$ and stopped at time $\zeta_{s'} - \zeta_r$.
	\end{itemize}
\end{itemize}

We suppose here that $\zeta(w)>0$. Let us set $I_a \de \inf\{ s\, : \, \zeta_s = a \}$ and let us define the probability measure on~$\Omega'$
$$\Pb_w^0 \de \Pb_w(\, \cdot \, |\, I_0 = 1).$$
This conditioning may be properly defined by saying that, under~$\Pb_w^0$, the law of $(\zeta_s)_{0 \le s \le 1}$ is the law of a first-passage Brownian bridge on $[0,1]$ from $\zeta(w)$ to~$0$, the law of $(\zeta_s)_{s \ge 1}$ is the law of a reflected Brownian motion on $[1, +\infty)$ started from~$0$, and the conditional distribution of $(W_s)_{s \ge 0}$ knowing $(\zeta_s)_{s\ge 0}$ is~$\Theta_w^\zeta$.

We denote by $\oo_\sigma\in \K$ the function $s \in [0,\sigma] \mapsto 0$. Under~$\Pb_{\oo_\sigma}^0$, the process $((\zeta_s)_{0 \le s \le 1},(\widehat W_s)_{0 \le s \le 1})$ has the same law as the process $\big( F_{[0,1]}^{\sigma \to 0}, Z_{[0,1]}\big)$ defined during Section~\ref{secbbfpbbs}. If we denote by~$\bbB$ the law on~$\K$ of a Brownian bridge on $[0,\sigma]$ from~$0$ to~$0$, multiplied by the factor~$\sqrt 3$, we then obtain that, under
$$\int_\K \bbB (dw) \, \Pb_w^0(d\omega),$$
the process $((\zeta_s)_{0 \le s \le 1},(\widehat W_s)_{0 \le s \le 1})$ has the same law as~$(C_\8,\Li)$ (under the common probability measure~$\Pb$).

\bigskip

We denote by~$n(de)$ the It\^o measure of positive Brownian excursions, whose normalization is given by the relation $n(\sup e > \eps) = 1/2\eps$, and we set
$$\N_x \de \int_{\C(\R_+,\R_+)} n(de)\, \Theta_{\bar x}^e$$
the excursion measure of the Brownian snake away from the path $\bar x :0 \mapsto x$. Under~$\Pb_w$, let us denote by $(\alpha_i, \beta_i)$, $i\in I$, the excursion intervals of $s\in [0, I_0] \mapsto \zeta_s-\ulz_{s}$, that is, the connected components of the open set $[0,I_0]\cap \{s \, : \, \zeta_s > \ulz_{s}\}$. For $i\in I$, we define $W^{(i)} \in \C(\R_+, \K)$ by setting, for $s\ge 0$,
$$W^{(i)}_s(t) = W_{(\alpha_i+s) \wedge \beta_i}(\zeta_{\alpha_i}+t), \qquad 0 \le t \le \zeta_s^{(i)}\de \zeta_{(\alpha_i+s) \wedge \beta_i} - \zeta_{\alpha_i}.$$

One of the main ingredients to our proofs is the following lemma.

\begin{lem}[{\cite[Lemma~V.5]{legall99sbp}}]\label{poisexc}
The point measure
\begin{equation*}\label{ppm}
\sum_{i\in I} \delta_{(\zeta_{\alpha_i},W^{(i)})}(dt\, d\omega)
\end{equation*}
is under~$\Pb_w$ a Poisson point measure on $\R_+\times \C(\R_+, \K)$ with intensity
$$\cN_w(dt \, d\omega) \de 2\, \un{[0,\zeta(w)]}(t) dt \, \N_{w(t)}(d\omega).$$
\end{lem}

We will also need the explicit ``law'' of the minimum of the Brownian snake's head under~$\N_x$.

\begin{lem}[{\cite[Lemma~2.1]{legallweill06cbt}}]\label{nxlaw}
For all $x$, $y \in \R$ with $y < x$,
$$\N_x \lp \min_{s\ge 0} \widehat W_s < y \rp = \frac{3}{2 (x- y)^2}.$$
\end{lem}

\bigskip

With this setting, we have two singular conditionings: one being $I_0=1$, and the second one being the fact that~$w$ is under $\bbB(dw)$ a bridge, instead of a Brownian motion. The first step in our proofs will generally be to dispose of the first of these conditionings (and sometimes the second as well), making us work under~$\Pb_w$ instead of~$\Pb_w^0$. This will usually be done by a simple absolute continuity argument, at least for almost sure properties. Another difficulty will arise from the factor~$\sqrt 3$, and we will sometimes need to take extra care because of it.

\subsection{Proof of Lemma~\ref{minb}}\label{secminb}

Thanks to Lemma~\ref{poisexc}, we will derive Lemma~\ref{minb} from the following similar result under~$\N_x$, which is due to Le~Gall and Weill~\cite{legallweill06cbt}.

\begin{prop}[{\cite[Proposition~2.5]{legallweill06cbt}}]\label{minw}
There exists~$\N_x$ a.e.\ a unique instant where $(\widehat W_s)_{s \ge 0}$ reaches its minimum.
\end{prop}

\begin{pre}[Proof of Lemma~\ref{minb}]
From a simple absolute continuity argument, it is sufficient to show that, for every $a\in [0,\zeta(w)]$, the process $(\widehat W_s)_{0 \le s \le I_a}$ reaches its minimum only once~$\Pb_w$ a.s., for every~$w$ belonging to a subset of~$\K$ of full $\bbB$-measure.
Without any assumption on~$w$, Lemmas~\ref{poisexc} and~\ref{nxlaw} imply that, $\Pb_w$ a.s.,
the process $(\widehat W_s)_{0 \le s \le I_a}$ does not reach its minimum on two different intervals of the form $[\alpha_i, \beta_i]$, $i\in I$. Moreover, the probability that it reaches its minimum more than once on some such interval is smaller than
\begin{multline*}
\Pb_w\lp \exists i \in I\, :\, \exists \alpha_i \le s < t \le \beta_i \, : \, \widehat W_s = \widehat W_t = \min_{s \in [\alpha_i,\beta_i]}\widehat W_s \rp\\
	= 1 - \exp \lp- 2 \int_0^{\zeta(w)} \! dt\ \N_{w(t)}\lp \exists s < t \, : \, \widehat W_s = \widehat W_t = \min_{s\ge 0}\widehat W_s \rp \rp=0,
\end{multline*}
by Proposition~\ref{minw}.

We will now see see that $(\widehat W_s)_{0 \le s \le I_a}$ does not reach its minimum on $[0,I_a]\bs \bigcup_{i\in I} [\alpha_i,\beta_i]$, which will complete the proof. It is at this time that we make extra assumptions on~$w$. The so-called snake property shows that
$$\lb\widehat W_s\, : \, s\in [0,I_a]\bs \bigcup_{i\in I} (\alpha_i,\beta_i)\rb = \big\{ w(t)\, : \, a \le t \le \zeta(w)\big\},$$
so that it will be enough to see that, $\Pb_w$ a.s., $\min_{0\le s \le I_a} \widehat W_s < \min_{[a, \zeta(w)]} w$. Using Lemma~\ref{poisexc} then Lemma~\ref{nxlaw}, we obtain
\begin{align*}
\Pb_w\lp \min_{0\le s \le I_a} \widehat W_s < \min_{[a, \zeta(w)]} w \rp
							&=1-\exp \lp- 2 \int_a^{\zeta(w)} \!dt\ \N_{w(t)}\lp \min_{s\ge 0}\widehat W_s < \min_{[a, \zeta(w)]} w \rp \rp\\
							&=1-\exp \lp- 3 \int_a^{\zeta(w)} \!dt\ \Big( w(t)-\min_{[a, \zeta(w)]} w \Big)^{-2}\rp.
\end{align*}
An easy application of L\'evy's modulus of continuity (see for example \cite[Theorem~I.2.7]{revuz99cma}) shows that, $\bbB(dw)$ a.s., this quantity equals~$1$.
\end{pre}

\subsection{Proof of Lemma~\ref{pcb}}\label{secippre}

For a continuous function $f:[0,\ell] \to \R$, we write $\IPl(f)$ (resp.\ $\IPr(f)$) the set of its left-increase points (resp.\ right-increase points). Remember that $s\in (0,\ell]$ is a left-increase point of~$f$ if there exists $t\in [0,s)$ satisfying $f(r) \ge f(s)$ for all $t\le r \le s$, and that a right-increase point is defined in a symmetrical way. We also denote by $\IP(f) = \IPl(f) \cup \IPr(f)$ the set of all its increase points. Due to the fact that the points $I_a$, $a\in [0,\zeta(w)]$ are left-increase points of~$\zeta$ and do not always lie in $\cup_{i\in I} (\alpha_i, \beta_i]$, we cannot directly apply the same strategy as in the previous section and derive Lemma~\ref{pcb} from a similar statement under~$\N_x$. Instead, we use a technique of covering intervals inspired from~\cite{bertoin91ilp} and a theorem of Shepp~\cite{shepp72clr}. In~\cite{bertoin91ilp}, Bertoin is interested in a similar problem: he characterizes the L\'evy processes~$X$ for which the set $\IPr(X)\cap \IPl(-X)$ is almost surely empty. Our method gives, in particular, another proof to \cite[Lemma~3.2]{legall08slb}, which states that the set
$$\IP\big((\zeta_s)_{0 \le s\le \ell}\big) \cap \IP\lp (\widehat W_s)_{0 \le s\le \ell} \rp$$
is~$\N_x$ a.e.\ empty. (Recall that we write $\ell\de \sup\{s\ge 0\, : \, \zeta_s >0\}$.) This comes very roughly from the fact that, if~$\zeta$ and~$\widehat W$ do not share any increase points on $[0,I_0]$, in particular, they do not share any increase points on any $(\alpha_i,\beta_i)$ either, and, by Lemma~\ref{poisexc}, the process restricted to $(\alpha_i,\beta_i)$ is then ``distributed'' under~$\N_x$.

For $y \in \R$, we set $\bmT_{y} \de \inf\{s \ge 0\, : \, w(s)=y\}$, where~$w$ is the canonical process on~$\K$, and, for $y < a$ and $\kappa >0$, we denote by $P_{a,\kappa}^{(y,\infty)}$ the law on~$\K$ of a standard Brownian motion multiplied by~$\kappa$, started from~$a$ and stopped at time~$\bmT_{y}$. For $x>0$, we also denote by $P_{a,\kappa}^x$ the law of a standard Brownian motion multiplied by~$\kappa$, started from~$a$ and stopped at time~$x$. When we omit the value of~$\kappa$, it will be assumed to be~$1$.

Although quite long to properly write in full detail, our strategy is pretty simple. One of the main difficulty comes from the two levels of randomness of the Brownian snake. In contrast to the previous proof where we worked under~$\Pb_w$ for a fixed~$w\in \K$, we will need here to work under $\bbB(dw)\Pb_w^0(d\omega)$ and see~$w$ as random. As a consequence, we will need to consider the timescale of~$\zeta$ and~$W$, as well as the timescale of~$w$. Juggling from one to the other may also cause confusion.

\begin{figure}[ht]
		\psfrag{x}[r][r]{$x$}
		\psfrag{s}[][]{$s$}
		\psfrag{z}[r][r]{$\zeta_s$}
		\psfrag{a}[][]{$w(\zeta_s)$}
		\psfrag{i}[][]{$I_0$}
		\psfrag{e}[][]{\textcolor{red}{$\zeta$}}
		\psfrag{w}[][]{\textcolor{blue}{$w$}}
	\centering\includegraphics[width=14cm]{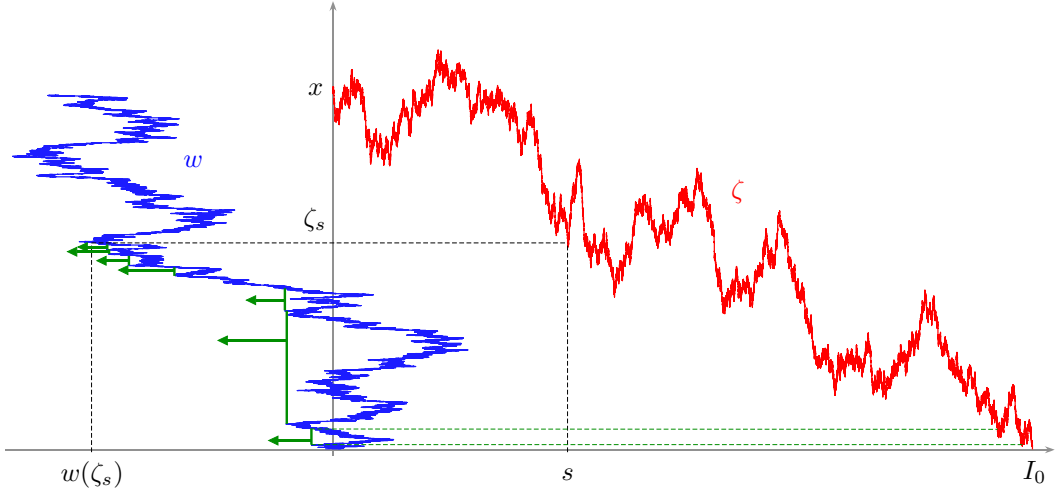}
	\caption[Increase points of $(\zeta,\widehat W)$.]{Visual aid for the proof of Lemma~\ref{pcb}. On this picture, the timescale of~$\zeta$ and~$W$ is horizontal, whereas the timescale of~$w$ is vertical. The point~$s$ satisfies $\zeta_s=\ulz_s=\inf\{t\, : w(t)= w(\zeta_s)\}$. We represented by (green) solid lines the seven longest excursions of $w-\ul w$ before time~$\zeta_s$, and the (green) arrows represent the minimum of~$\widehat W$ on the intervals that correspond in the timescale of~$W$ to these excursions.}
	\label{figip}
\end{figure}

In order to facilitate the reading of our proof, we outline it now. By absolute continuity arguments, we get rid of the conditionings and work under $P_{0,\sqrt 3}^x(dw) \, \Pb_{w}(d\omega)$ instead of $\bbB(dw)\Pb_w^0(d\omega)$. Using symmetry and time-reversal, we mainly need to focus on right-increase points of~$\widehat W$ that are also left-increase points of~$\zeta$. It should not be too hard to convince oneself\footnote{To be more accurate, if~$s$ does not satisfy this hypothesis, we apply the Markov property at some (rational) time~$a$ close enough before~$s$ so that $\zeta_s=\inf_{[a,s]} \zeta$. When doing so, we work under $P_0^x(dw) \, \Pb_{w}(d\omega)$ instead of $P_{0,\sqrt 3}^x(dw) \, \Pb_{w}(d\omega)$.} that it suffices to look at points $s\in [0,I_0]$ such that $\zeta_s=\ulz_s$. If~$s$ is such a point and also a right-increase point of~$\widehat W$, we will first see that~$s$ is not the starting point of an excursion of $\zeta-\ulz$, that is, $s\notin\{\alpha_i\ : \ i\in I\}$. As a result, as moreover~$s$ is a right-increase point of~$\widehat W$ and~$\ulz$ is non-increasing, we obtain that~$\ulz_s$ is a left-increase point of~$w$. By an argument similar as before, we may restrict our attention to points~$s$ satisfying $\zeta_s=\ulz_s= \inf\{t\, : w(t)= w(\zeta_s)\}$. See Figure~\ref{figip}.

Then, we will consider the excursions of $w-\ul w$ and look at the minimum of~$\widehat W$ on the intervals corresponding in the timescale of~$\widehat W$ to these excursions. Using~\cite{shepp72clr}, we will see that, as close as we want before~$\zeta_s$, we can find an excursion of $w-\ul w$ where the corresponding minimum of~$\widehat W$ is smaller than~$w(\zeta_s)$, prohibiting~$\zeta_s$ from being a left-increase point of~$w$.

We start with a lemma stating that the extremities of any excursion interval $(\alpha_i,\beta_i)$ are not increase points of the process~$\widehat W$ restricted to this interval $(\alpha_i,\beta_i)$.

\begin{lem}\label{lemipnx}
Let $w\in \K$. Then, $\Pb_w(d\omega)$ a.s., for all $i\in I$,
$$\alpha_i \notin \IPr\lp (\widehat W_s)_{0 \le s\le {I_0}} \rp \sand \beta_i \notin \IPl\lp (\widehat W_s)_{0 \le s\le {I_0}} \rp.$$
\end{lem}

\begin{pre}
It is enough to show that $\N_x$ a.e.\ $0 \notin \IPr\lp (\widehat W_s)_{0 \le s\le \ell} \rp$. Indeed, this entails by Lemma~\ref{poisexc} that
\begin{multline*}
\Pb_w \lp \exists i \in I \, : \, \alpha_i \in \IPr\lp (\widehat W_s)_{0 \le s\le {I_0}} \rp \rp\\
	= 1 - \exp \lp -2 \int_0^{\zeta(w)} \! dt\ \N_{w(t)}\lp 0 \in \IPr\lp (\widehat W_s)_{0 \le s\le \ell} \rp \rp \rp =0.$$
\end{multline*}
Then, by the time-reversal property under~$\N_x$ (the process $(\zeta_{\ell-s}, W_{\ell-s})_{0\le s\le \ell}$ has under $\N_x$ the same distribution as $(\zeta_{s}, W_{s})_{0\le s\le \ell}$), we see that $\N_x$ a.e.\ $\ell \notin \IPl\lp (\widehat W_s)_{0 \le s\le \ell} \rp$, and we conclude in the same way.

Let~$e$ be an excursion. By definition, under $\Theta_{\bar x}^e$, the process
$$\lp \widehat W (\sup \{s \le \ell/2 \, :\, e_s=y\}) \rp_{0\le y \le e_{\ell/2}}$$
has the law $P_x^{e_{\ell/2}}$. The desired result follows since $n(de)$ a.e.\ $\sup\{s \le \ell/2 \, :\, e_s=y\}\to 0$ as $y\to 0$.
\end{pre}

The following lemma will only be used for $\kappa=1$ or $\kappa=\sqrt 3$ in what follows, but the proof works for any $\kappa \le \sqrt 3$, so that we consider all these values.

\begin{lem}\label{pioneer}
Let $\kappa \le \sqrt 3$ and $x>0$. The sets
$$\A \de \lb s \in \IPr\lp (\widehat W_s)_{0 \le s\le {I_0}} \rp \, : \, \zeta_s = \ulz_s= \inf\{t\, : w(t)= w(\zeta_s)\} \rb$$
and
$$\cB \de \lb s \in \IPl\lp (\widehat W_s)_{0 \le s\le {I_0}} \rp \, : \, \zeta_s = \ulz_s= \sup\{t\, : w(t)= w(\zeta_s)\} \rb$$
are $P_{0,\kappa}^x(dw) \, \Pb_{w}(d\omega)$ a.s.\ empty.
\end{lem}

\begin{pre}
For technical reasons, it will be easier to work with the measure $P_{0,\kappa}^{(-y,\infty)}(dw)$ instead of $P_{0,\kappa}^x(dw)$. Let $y >0$. We are going to show that the set~$\A$ is $P_{0,\kappa}^{(-y,\infty)}(dw) \, \Pb_{w}(d\omega)$ a.s.\ empty. This will entail in particular that the set
$$\lb s \in \IPr\lp (\widehat W_s)_{I_x \le s\le {I_0}} \rp \, : \, \zeta_s = \ulz_s= \inf\{t\, : w(t)= w(\zeta_s)\} \rb$$
is $P_{0,\kappa}^{(-y,\infty)}(dw\ | \ \bmT_{-y} \ge x)\, \Pb_{w}(d\omega)$ a.s.\ empty. By the Markov property of the Brownian snake at time~$I_x$, under the latter measure, the distribution of
$$\big( (w(s))_{0\le s \le x}, (\zeta_{I_x+s})_{0 \le s \le I_0-I_x}, (W_{I_x+s})_{0 \le s \le I_0-I_x} \big)$$
is precisely $P_{0,\kappa}^x(dw\ | \ \ul w_x \ge -y) \, \Pb_{w}(d\omega)$. Letting $y \to \infty$ yields that~$\A$ is  $P_{0,\kappa}^x(dw) \, \Pb_{w}(d\omega)$ a.s.\ empty.

\paragraph{Step 1.} Let us denote by $(u_j, v_j)$, $j \in J$, the excursion intervals of $w - \ul w$, and
$$w^{(j)}(s) \de w((u_j+s) \wedge v_j) - w(u_j)\qquad j \in J.$$
We will need to find the distribution under $P_{0,\kappa}^{(-y,\infty)}(dw) \, \Pb_{w}(d\omega)$ of the point measure
$$\PP \de \sum_{j\in J} \delta_{(-w(u_j),\, m^{(j)})} \qquad \text{ where }m^{(j)}\de w(u_j)-\min_{[I_{v_j},I_{u_j}]} \widehat W.$$
To this end, we adapt a computation of Miermont \cite[Lemma~31]{miermont11bms}. By It\^o's excursion theory, the point measure
$$\sum_{j\in J} \delta_{(-w(u_j),\, w^{(j)})}$$
is under $P_{0,\kappa}^{(-y,\infty)}(dw)$ a Poisson point measure on $\R_+ \times \K$ with intensity $\kappa^{-1}\un{[0,y]}(t) dt \, 2\kappa_* n(de)$, where $\kappa_* n(de)$ denotes the pushforward of $n(de)$ by $e\mapsto \kappa \, e$. Using Lemma~\ref{poisexc}, we may see the $m^{(j)}$, $j\in J$, as independent marks on $w^{(j)}$, $j\in J$, with law $\Pb_{w^{(j)}} ( - \min \widehat W \in dz )$. The marking theorem of Poisson point measures \cite[Marking Theorem]{kingman93pp} shows that, under $P_{0,\kappa}^{(-y,\infty)}(dw) \, \Pb_{w}(d\omega)$, $\PP$ is also a Poisson point measure on $\R_+ \times \R_+$ with intensity
$$\kappa^{-1} \un{[0,y]}(t) dt \int_\K 2n(de) \, \Pb_{\kappa e} \lp - \min \widehat W \in dz \rp.$$
To compute explicitly this intensity, we use Lemmas~\ref{poisexc} and~\ref{nxlaw}, and then Bismut's description of~$n$ \cite[Theorem~XII.4.7]{revuz99cma}:
\begin{align*}
\int_\K 2n(de) \, \Pb_{\kappa e} \lp - \min \widehat W \ge z \rp 
				&= \int_\K 2n(de) \lp 1 - \exp \lp - \int_0^\ell \frac{3\, ds}{(\kappa e_s + z)^2} \rp \rp\\
				&= \int_\K 2n(de) \int_0^\ell \frac{3\, dt}{(\kappa e_t + z)^2} \, \exp \lp - \int_t^\ell \frac{3\, ds}{(\kappa e_s + z)^2} \rp \\
				&= 6 \int_0^\infty \frac{da}{(\kappa a + z)^2} \, E_a^{(0,\infty)}\lbr \exp \lp - \int_0^{\bmT_0} \frac{3\, ds}{(\kappa w(s) + z)^2} \rp \rbr\\
				&= 6 \int_0^\infty \frac{da}{(\kappa a + z)^2} \, E_{a+z/\kappa}^{(0,\infty)}\lbr \exp \lp - \frac{3}{\kappa^2} \int_0^{\bmT_{z/\kappa}} \frac{ds}{(w(s))^2} \rp \rbr.
\end{align*}
Using the absolute continuity relations between Bessel processes with different indices, which is due to Yor \cite[Exercise~XI.1.22]{revuz99cma}, and the fact that reflected Brownian motion is a $1$-dimensional Bessel process, we see that
\begin{align*}
E_{a+z/\kappa}^{(0,\infty)}\lbr \exp \lp - \frac{3}{\kappa^2} \int_0^{\bmT_{z/\kappa}} \frac{ds}{(w(s))^2} \rp \rbr
		&= \lim_{t\to \infty} E_{a+z/\kappa}^{(0,\infty)}\lbr \exp \lp - \frac{3}{\kappa^2} \int_0^{\bmT_{z/\kappa}\wedge t} \frac{ds}{(w(s))^2} \rp \rbr\\
		&= \lim_{t\to \infty} E_{a+z/\kappa}^{\langle 2+2\nu \rangle}\lbr \lp \frac{w(\bmT_{z/\kappa}\wedge t)}{a+z/\kappa} \rp^{-\nu-1/2} \rbr\\
		&=\lp \frac{z}{\kappa a + z}\rp^{-\nu-1/2} P_{a+z/\kappa}^{\langle 2+2\nu \rangle}\lbr \bmT_{z/\kappa} < \infty \rbr\\
		&=\lp \frac{z}{\kappa a + z}\rp^{\nu -1/2}.
\end{align*}
where $\nu \de \sqrt{24+\kappa^2}\, /2\kappa$ and $P_a^{\langle 2+2\nu \rangle}$ denotes the distribution of a Bessel process of dimension $2+2\nu$. In the last line, we used the fact that, for $b > c$, $P_{b}^{\langle 2+2\nu \rangle}\lbr \bmT_{c} < \infty \rbr= (c/b)^{2\nu}$ (see \cite[Chapter~XI]{revuz99cma}). Putting all this together and differentiating with respect to~$z$, we obtain that the intensity of~$\PP$ is
$$\un{[0,y]}(t) dt \, \frac{\lambda}{z^2} \, dz, \qquad \text{ where } \lambda \de \frac{12\kappa^{-1}}{\sqrt{24+\kappa^2} + \kappa} \ge 1.$$

\paragraph{Step 2.}  Let $\sum_{k\in K} \delta_{(t_k,z_k)}$ be a Poisson random measure with intensity $dt \, \lambda z^{-2} dz$. Then, by the restriction property of Poisson random measures, for all $\eps >0$, $\sum_{k\in K} \delta_{(t_k,z_k)} \un{\{z_k \le \eps\}}$ is a Poisson random measure with intensity $dt \, \lambda z^{-2} \un{\{z\le \eps\}} \, dz$. By a theorem of Shepp~\cite{shepp72clr}, we obtain that the random set
\begin{equation}\label{shepp}
\bigcup_{k\in K\, : \, z_k \le \eps} (t_k,t_k+z_k)
\end{equation}
is a.s.\ equal to $\R$. We used here the fact that $\lambda \ge 1$, ensuring that~$\R$ is covered with ``small'' intervals. See in particular the remark on high frequency coverings in \cite[Section~5]{shepp72clr}. As a result, the set
$$\bigcup_{z_k \le \eps,\, 0\le t_k \le y} (t_k,t_k+z_k)$$
a.s.\ contains $[\eps,y]$. Because the point measure $\sum_{k\in K} \delta_{(t_k,z_k)}\un{\{0\le t_k \le y\}}$ has the same law as~$\PP$ under $P_{0,\kappa}^{(-y,\infty)}(dw) \, \Pb_{w}(d\omega)$, we find that, $P_{0,\kappa}^{(-y,\infty)}(dw) \, \Pb_{w}(d\omega)$ a.s., for all rational $\eps>0$, the set $[-y,-\eps]$ is contained in
$$\operatorname{Cov}_\eps \de \bigcup_{j\in J\, : \, m^{(j)} \le \eps} \lp w(u_j)-m^{(j)}, w(u_j)\rp.$$

Now, let us assume that~$\A$ is not empty, and let us take $s\in \A$. By Lemma~\ref{lemipnx}, $s \notin \{\alpha_i, \, i\in I\}$. As a result, there exists $\eta >0$ such that $\widehat W_r \ge \widehat W_s$ for all $r \in [s, I_{\zeta_s-\eta}]$. In particular, for all $j\in J$ such that $[u_j,v_j] \subseteq [\zeta_s-\eta, \zeta_s]$, we have $w(u_j)-m^{(j)} \ge \widehat W_s=w(\zeta_s)$. As $\zeta_s = \inf\{t\, : w(t)= w(\zeta_s)\}$, we can find a rational $\eps >0$ satisfying $w(\zeta_s) \le \ul w(\zeta_s-\eta)-\eps \le -\eps$.

Let $j\in J$ be such that $m^{(j)} \le \eps$. If $u_j \le \zeta_s-\eta$, then $w(u_j)-m^{(j)} \ge \ul w(\zeta_s-\eta)-\eps \ge w(\zeta_s)$. If $u_j \in [\zeta_s-\eta, \zeta_s]$, we already observed that $w(u_j)-m^{(j)} \ge w(\zeta_s)$. Finally, if $u_j \ge \zeta_s$, then $w(u_j) \le w(\zeta_s)$. In all cases, $w(\zeta_s) \notin (w(u_j)-m^{(j)}, w(u_j))$. We found a point $w(\zeta_s) \in [-y,-\eps]$ that does not belong to $\operatorname{Cov}_\eps$. This can only happen with probability~$0$.

\bigskip

Similar arguments show that the set~$\cB$ is $\ol P_{0,\kappa}^{(-y,\infty)}(dw)\, \Pb_{w}(d\omega)$ a.s.\ empty, where we write $\ol P_{0,\kappa}^{(-y,\infty)}$ the pushforward of $P_{0,\kappa}^{(-y,\infty)}$ under $w\mapsto \ol w \de (w(\zeta(w)-s))_{0\le s \le \zeta(w)}$. This entails that~$\cB$ is also $\ol P_{0,\kappa}^{x}(dw)\, \Pb_{w}(d\omega)$ a.s.\ empty, and, by time-reversal, $\ol w - w(x)$ has under $P_{0,\kappa}^{x}(dw)$ the same distribution as~$w$, so that the result also holds $P_{0,\kappa}^{x}(dw) \, \Pb_{w}(d\omega)$ a.s. We leave the details to the reader.
\end{pre}

\begin{nota}\label{notroot3}
We can see from this proof that the value $\sqrt 3$ is critical. The theorem of Shepp actually entails that the set~\eqref{shepp} is a.s.\ equal to~$\R$ if and only if $\lambda \ge 1$, that is, $\kappa \le \sqrt 3$. In the case $\kappa > \sqrt 3$, we can thus find an $\eps>0$ and a point~$z$ in $[-y,-\eps]$ that does not belong to $\operatorname{Cov}_\eps$. It is then not very hard to see that $I_{\inf\{t\ : \ w(t)=z\}}$ is a right-increase point of $\widehat W$. As a result, the set $\IP\big((\zeta_s)_{0 \le s\le \ell}\big) \cap \IP\lp (\widehat W_s)_{0 \le s\le \ell} \rp$ is not empty.
\end{nota}

We may now proceed to the proof of Lemma~\ref{pcb}. We define
$$I_b^{(a)} \de \inf\{ s \ge a \, : \, \zeta_s = b \}.$$

\begin{pre}[Proof of Lemma~\ref{pcb}]
As above, we start by working under $P_{0,\kappa}^{x}(dw) \, \Pb_{w}(d\omega)$ with $\kappa \le \sqrt 3$.

\paragraph{Step 1.} The first step consists in treating the left-increase points of~$\zeta$. To do so, we will use the Markov property of the Brownian snake and ``insert'' rational numbers in order to be able to apply the previous lemmas. Let $b \in [0,x]$. Because the process
$$\lp \big( w(b+r)-w(b) \big)_{0 \le r \le x-b}, \big( (W_r(b+t))_{0 \le t \le \zeta_r-b} \big)_{0 \le r \le I_b} \rp$$
has under $P_{0,\kappa}^{x}(dw) \, \Pb_{w}(d\omega)$ the law $P_{0,\kappa}^{x-b}(dw) \, \Pb_{w}(d\omega)$, we see by Lemma~\ref{pioneer} that the set
$$\A^b \de \lb s \in \IPr\lp (\widehat W_s)_{0 \le s\le {I_b}} \rp \, : \, \zeta_s = \ulz_s= \inf\{t \ge b\, : w(t)= w(\zeta_s)\} \rb$$
is $P_{0,\kappa}^{x}(dw) \, \Pb_{w}(d\omega)$ a.s.\ empty. Similarly, for $b \in [0,x]$, $c > b$, and $a$, the Markov property shows that the process
$$\lp \big( W_a(b+r)-W_a(b) \big)_{0 \le r \le c-b}, \big( (W_{I_c^{(a)}+r}(b+t))_{0 \le t \le \zeta_{I_c^{(a)}+r}-b} \big)_{0 \le r \le I_b^{(a)}-I_c^{(a)}} \rp$$
has under $P_{0,\kappa}^{x}(dw) \, \Pb_{w}(d\omega\ | \ a \le I_0,\, \ulz_a <b, \, c < \zeta_a )$ the law $P_0^{c-b}(dw)\, \Pb_{w}(d\omega)$. (Beware that here the factor $\kappa$ does not appear.) As a result, Lemmas~\ref{lemipnx} and~\ref{pioneer} successively show that, on the event $\{a \le I_0, \, \ulz_a <b, \, c < \zeta_a\}$, the sets
$$\C_a^{b,c} \de \lb s \in \IPr\lp (\widehat W_s)_{I_c^{(a)} \le s\le {I_b^{(a)}}} \rp \cap \IPr(\zeta)  \, : \, \zeta_s = \inf_{[a,s]} \zeta\rb$$
and
$$\A_a^{b,c} \de \lb s \in \IPr\lp (\widehat W_s)_{I_c^{(a)} \le s\le {I_b^{(a)}}} \rp \, : \, \zeta_s = \inf_{[a,s]} \zeta= \inf\{t \ge b\, : W_a(t)= W_a(\zeta_s)\} \rb$$
are $P_{0,\kappa}^{x}(dw) \, \Pb_{w}(d\omega)$ a.s.\ empty. As a result, we obtain that $P_{0,\kappa}^{x}(dw) \, \Pb_{w}(d\omega)$ a.s., for all rational values of~$a$, $b$, and~$c$, these sets are empty.

\bigskip

Now, if the set $\IPl\big((\zeta_s)_{0 \le s \le I_0}\big) \cap \IPr \big((\widehat W_s)_{0 \le s \le I_0} \big)$ is not empty, let~$s$ be a point lying in it. Let us first suppose that $\zeta_s = \ulz_s$. By Lemma~\ref{lemipnx}, we know that $s \notin \{\alpha_i, i\in I\}$. This implies that $\zeta_s \in \IPl(w)$. As local minimums of Brownian motion are distinct, we can find a rational $b \in [0,\zeta_s]$ such that $\zeta_s = \inf\{t \ge b\, : w(t)= w(\zeta_s)\}$, and $s \in \A_b$. Otherwise, $\zeta_s > \ulz_s$. As $s \in \IPl(\zeta)$, we can find a rational $a\in [0,s)$ such that $\zeta_a > \zeta_s$ and $\zeta_s = \inf_{[a,s]} \zeta$. If $s \in \IPr(\zeta)$, we can find rationals $b \in(\ulz_a, \zeta_s)$ and $c \in (\zeta_s, \zeta_a)$ so that $s\in \C_a^{b,c}$. If $s \notin \IPr(\zeta)$, then $\zeta_s \in \IPl(W_a)$. We can then find rationals~$b$ and~$c$ such that $s\in \A_a^{b,c}$.

Summing up, we obtain that $\IPl\big((\zeta_s)_{0 \le s \le I_0}\big) \cap \IPr \big((\widehat W_s)_{0 \le s \le I_0} \big)$ is $P_{0,\kappa}^{x}(dw) \, \Pb_{w}(d\omega)$ a.s.\ empty. By a similar argument, we show that the set $\IPl\big((\zeta_s)_{0 \le s \le I_0}\big) \cap \IPl \big((\widehat W_s)_{0 \le s \le I_0} \big)$ is also $P_{0,\kappa}^{x}(dw) \, \Pb_{w}(d\omega)$ a.s.\ empty.

\paragraph{Step 2.} We now use a time-reversal argument under~$\N_y$ to treat the right-increase points of~$\zeta$. By translation, the quantity
$$\Delta \de \N_y\lp \IPl\big((\zeta_s)_{0 \le s\le \ell}\big) \cap \IP\big((\widehat W_s)_{0 \le s\le \ell} \big) \neq \varnothing \rp$$
does not depend on~$y$. Using the Poissonian decomposition of the excursions of $\zeta-\ulz$, we see that, $P_{0,\kappa}^{x}(dw)$ a.s., 
$$0 = \Pb_{w}\lp \IPl\big((\zeta_s)_{0 \le s \le I_0}\big) \cap \IP \big((\widehat W_s)_{0 \le s \le I_0} \big) \neq \varnothing \rp \ge 1-e^{-2x \Delta},$$
so that $\Delta =0$. Note that a priori we only have an inequality, because some left-increase points of $(\zeta_s)_{0 \le s \le I_0}$ may well lie outside of the set $\cup_{i \in I} (\alpha_i, \beta_i]$. Using time-reversal under~$\N_y$, we find that $\IPr\big((\zeta_s)_{0 \le s \le \ell}\big) \cap \IP \big((\widehat W_s)_{0 \le s \le \ell} \big)$ is also $\N_y$ a.e.\ empty. As announced at the beginning of this section, we re-obtained here \cite[Lemma~3.2]{legall08slb}. It is then easier to deal with right-increase points of $(\zeta_s)_{0 \le s \le I_0}$, because they all lie in $\cup_{i \in I} [\alpha_i, \beta_i)$: using once again the Poissonian decomposition of the excursions of $\zeta-\ulz$, we find (for any $w\in \K$)
$$\Pb_{w}\lp \IPr\big((\zeta_s)_{0 \le s \le I_0}\big) \cap \IP \big((\widehat W_s)_{0 \le s \le I_0} \big) \neq \varnothing \rp = 1-e^{-2\zeta(w) \Delta}=0.$$

Putting it all together, we showed that $\IP\big((\zeta_s)_{0 \le s \le I_0}\big) \cap \IP \big((\widehat W_s)_{0 \le s \le I_0}\big)$ is $P_{0,\kappa}^{x}(dw) \, \Pb_{w}(d\omega)$ a.s.\ empty. Using the fact that the distribution of $(w(s))_{0 \le s \le \sigma-\eps}$ under~$\bbB$ is absolutely continuous with respect to the distribution of $(w(s))_{0 \le s \le \sigma-\eps}$ under $P_{0,\sqrt 3}^{\sigma-\eps}(dw)$, we obtain that, $\bbB(dw)\, \Pb_{w}(d\omega)$ a.s., for all rational $\eps \in (0,\sigma)$, $\IP\big((\zeta_s)_{I_{\sigma-\eps} \le s \le I_0}\big) \cap \IP \big((\widehat W_s)_{I_{\sigma-\eps} \le s \le I_0}\big)=\varnothing$. Standard properties of Brownian motion show that, $\Pb_w$ a.s.\ $0 \notin \IP(\zeta)$ and $\inf_{\eps \in (0,\sigma)\cap\bbQ} I_{\sigma-\eps} =0$, so that $\bbB(dw)\, \Pb_{w}(d\omega)$ a.s. $\IP\big((\zeta_s)_{0 \le s \le I_0}\big) \cap \IP \big((\widehat W_s)_{0 \le s \le I_0}\big)=\varnothing$. Using another absolute continuity argument and the fact that $1\notin \IPl(\widehat W)$ (because $0\notin \IP(w)$), we conclude that $\int_\K \bbB (dw) \, \Pb_w^0(d\omega)$ a.s.\ $\IP\big((\zeta_s)_{0 \le s \le 1}\big) \cap \IP \big((\widehat W_s)_{0 \le s \le 1}\big)=\varnothing$. This completes the proof.
\end{pre}

\subsection{Proof of Lemma~\ref{legalllemb}}

As explained in the proof of \cite[Lemma~11]{bettinelli10tsl}, Lemma~\ref{legalllemb} is a consequence of the following two lemmas, which we state here directly in terms of the Brownian snake. We will use a strategy similar to that of Section~\ref{secminb} and derive these lemmas from similar statements under~$\N_x$, namely \cite[Lemma~5.3]{legall07tss} and \cite[Lemma~6.1]{legall08glp}. We denote the Lebesgue measure on~$\R$ by~$\Leb$.

\begin{lem}\label{lemlg1}
Let $w\in \K$. $\Pb_w^0$ a.s., for every $\eta > 0$, for all $x\in [0,1]$ and all $l<r$ such that,
\begin{itemize}
	\item either $0 < l < r < x$ and $\zeta_l=\zeta_r=\inf_{[l,x]}\zeta$,
	\item or $x < l < r < 1$ and $\zeta_l=\zeta_r=\inf_{[x,r]}\zeta$,
\end{itemize}
the condition $\inf_{[l,r]} \widehat W < \widehat W_l - \eta$ implies that
$$\liminf_{\eps \to 0} \eps^{-2} \Leb \Big( \Big\{ s \in [l,r] \, : \, \widehat W_s < \widehat W_l - \eta + \eps\,;\ \forall y \in [\zeta_l, \zeta_s],
		\widehat W_{\sup\{t\le s\,:\, \zeta_t=y\}} > \widehat W_l-\eta+\frac\eps 8 \Big\} \Big) >0.$$
\end{lem}

\begin{lem}\label{lemlg2}
For every $p\ge 1$ and every $\delta \in (0,1]$, there exists a constant $c_{p,\delta} < \infty$ such that, for every $\eps > 0$,
$$\int_\K \bbB (dw) \, \mathbb E_w^0 \lbr \lp\int_0^1 
\un{\{\widehat W_s \le \min_{0 \le r \le 1} \widehat W_r + \eps\}}\, ds\rp^p \rbr \le c_{p,\delta} \, \eps^{4p-\delta}.$$
\end{lem}

\begin{pre}[Proof of Lemma~\ref{lemlg1}]
By an absolute continuity argument, it is sufficient to show the result under~$\Pb_w$ (while replacing~$1$ with~$I_0$). By \cite[Lemma~5.3]{legall07tss} and the preceding claim in~\cite{legall07tss}, the result holds under~$\N_0$ (while replacing~$1$ with~$\ell$) and the same result also holds if we replace the hypothesis on~$l$ and~$r$ by $l=0$ and $r=\ell$.

Now observe that, under~$\Pb_w$, if the numbers~$l$, $r$, $x$, and $\eta$ satisfy the hypothesis but not the conclusion of our statement, then there exists an $i\in I$ such that, either $(l,r)=(\alpha_i,\beta_i)$, in which case $0$, $\beta_i-\alpha_i$, and~$\eta$ also satisfy the hypothesis but not the conclusion for the process $(\zeta^{(i)},W^{(i)})$, or $l-\alpha_i$, $r-\alpha_i$, $x-\alpha_i$, and $\eta$ also satisfy the hypothesis but not the conclusion for the process $(\zeta^{(i)},W^{(i)})$. It is then an easy application of Lemma~\ref{poisexc} to show that the probability that there exists such numbers is equal to~$0$.
\end{pre}

\begin{pre}[Proof of Lemma~\ref{lemlg2}]
We fix $p\ge 1$ and $\delta \in (0,1]$. 

\paragraph{Step 1.}
We use \cite[Lemma~6.1]{legall08glp} and Bismut's description of~$n$ \cite[Theorem~XII.4.7]{revuz99cma} in order to derive a result similar to the one we seek but without the conditionings.

Let $0 < x < y$. Using Bismut's description of~$n$, we obtain that
\begin{equation}\label{biseq}
\N_0\lp \int_0^\ell dt \, \un{\{x\le \zeta_t \le y\}} \lp\int_0^{\ell} \un{\{\widehat W_s \le \min \widehat W + \eps\}}\, ds\rp^p \rp = \int_x^y da \, \Phi_\eps(a)
\end{equation}
where
$$\Phi_\eps(a) \de  \int_\K P_0^{a} (dw) \int_{\Omega'\times\Omega'} \Pb_{w}(dW^1)\Pb_{w}(dW^2) \lp\sum_{k=1}^2\int_0^{I_0^k} \un{\{\widehat W_s^k \le m + \eps\}}\, ds\rp^p,$$
$m \de \min \widehat W^1 \wedge \min \widehat W^2$, $I_0^1 \de \inf\{ s\, : \, \zeta(W_s^1) = 0 \}$ and $I_0^2 \de \inf\{ s\, : \, \zeta(W_s^2) = 0 \}$.
We now give an expression for $\Phi_\eps(a)$ that involves only one Brownian snake instead of two. As the Lebesgue measure of $\bigcup_{i\in I} [\alpha_i,\beta_i]$ is~$\Pb_w$ a.s.\ equal to~$I_0$, we see that
\begin{align*}
\sum_{k=1}^2\int_0^{I_0^k} \un{\{\widehat W_s^k \le m + \eps\}}\, ds &= \sum_{k=1}^2 \sum_{i\in I^k}\int_0^{\beta_i^k-\alpha_i^k} \un{\{\widehat W_s^{k,(i)} \le m + \eps\}}\, ds
\end{align*}
where we use an obvious notation for quantities defined in terms of~$W^k$. Using the fact proven during the proof of Lemma~\ref{minb} that, under~$\Pb_w$, $(\widehat W_s)_{0 \le s \le I_0}$ reaches its minimum on $\bigcup_{i\in I} [\alpha_i,\beta_i]$, we see that
$$m = \min \big\{\min X\ : \ X\in \big\{\widehat W^{1,(j)},\, j\in I^1\big\}\cup\big\{\widehat W^{2,(j)},\, j\in I^2\big\}\big\}.$$
Moreover, Lemma~\ref{poisexc}, together with a small computation, shows that the point process
$$\sum_{i\in I^1} \delta_{(\zeta^1_{\alpha_i^1},W^{1,(i)})}(dt\, d\omega)+\sum_{i\in I^2} \delta_{(\zeta^2_{\alpha_i^2},W^{2,(i)})}(dt\, d\omega)$$
under $\Pb_{w}(dW^1)\Pb_{w}(dW^2)$ has the same distribution as
$$\sum_{i\in I} \delta_{(\zeta_{\alpha_i}/2,W^{(i)})}(dt\, d\omega)$$
under~$\Pb_{\tilde w}(dW)$, where $\tilde w : s \in [0,2\zeta(w)] \mapsto w(s/2)$.
As a result,
\begin{equation*}
\int_{\Omega'\times\Omega'} \Pb_{w}(dW^1)\Pb_{w}(dW^2) \lp\sum_{k=1}^2\int_0^{I_0^k} \un{\{\widehat W_s^k \le m + \eps\}}\, ds\rp^p
=\mathbb E_{\tilde w} \lbr \lp\int_0^{I_0} \un{\{\widehat W_s \le \min \widehat W + \eps\}}\, ds\rp^p \rbr.
\end{equation*}
As~$\tilde w$ has under $P_0^{a} (dw)$ the same distribution as~$w$ under $P_{0,1/\sqrt 2}^{2a} (dw)$, we obtain that
\begin{equation}\label{phia}
\Phi_\eps(a) = \int_\K P_{0,1/\sqrt 2}^{2a} (dw)\, \mathbb E_{w} \lbr \lp\int_0^{I_0} \un{\{\widehat W_s \le \min \widehat W + \eps\}}\, ds\rp^p \rbr.
\end{equation}

\bigskip
Using equation~\eqref{biseq} and H\"older's inequality, we obtain that for any $k\ge 2$ and $A>0$,
$$\int_x^y da \,\Phi_\eps(a) \le \N_0\lp\lp \int_0^\ell dt \, \un{\{0\le \zeta_t \le A\}}\rp^k\rp^{1/k} \ 
\N_0\lp \lp\int_0^{\ell} \un{\{\widehat W_s \le \min \widehat W + \eps\}}\, ds\rp^{pk/(k-1)} \rp^{1-1/k}.$$
Noticing that the first term of the right-hand side is finite (this integral may actually easily be computed using \cite[Proposition~III.3]{legall99sbp}) and bounding the second term using \cite[Lemma~6.1]{legall08glp}, we see that there exists a finite constant $C_{A,p,\delta}$ such that, for any $0<x<y<A$,
$$\int_x^y da \,\Phi_\eps(a) \le C_{A,p,\delta} \, \eps^{4p-\delta}.$$

Let us now dispose of the integration over~$a$. A simple application of the snake's Markov property shows that, for any $a\in [-A/2,A/2]$, we have $\Phi_{\eps}(A) \le 2^p(\Phi_{\eps}(A/2+a) + \Phi_{\eps}(A/2-a))$. Integrating over~$a$ and applying the latter inequality shows that there exists a finite constant $C'_{A,p,\delta}$ such that
\begin{equation}\label{phiab}
\Phi_{\eps}(A) \le C'_{A,p,\delta} \, \eps^{4p-\delta}.
\end{equation}

\paragraph{Step 2.} Combining~\eqref{phia} and~\eqref{phiab} almost gives the desired result. There are still two problems we need to address. First, the diffusion factor in~\eqref{phia} is $1/\sqrt 2$ and it will be~$\sqrt 3$ in the process we are interested in. We now turn to this difficulty. Let us consider a sequence $(X_i)_{i\in I}$ of i.i.d.\ random Bernoulli variables with mean~$1/6$, independent from any other variables. Then, the marking theorem of Poisson point measures entails that the process $\sum_{i\in I} \delta_{(6\zeta_{\alpha_i},W^{(i)})}$ has, under $\int_\K P_{0,{\sqrt 3}}^{a} (dw) \, \Pb_{w}(d\omega)$, the same distribution as the process
$$\sum_{i\in I} \delta_{(\zeta_{\alpha_i},W^{(i)})}\un{\{X_i=1\}},\quad\text{ under }\quad\int_\K P_{0,1/\sqrt 2}^{6a} (dw) \, \Pb_{w}(d\omega).$$
As a result, writing $I'\de \{i\in I\, : \, X_i=1\}$, we obtain that
\begin{align*}
\int_\K P_{0,{\sqrt 3}}^{a} (dw) \, \mathbb E_{w} \lbr \lp\int_0^{I_0} \un{\{\widehat W_s \le \sulW_{I_0} + \eps\}}\, ds\rp^p \rbr\hspace{-30mm}\\
	&=\int_\K P_{0,{\sqrt 3}}^{a} (dw) \, \mathbb E_{w} \lbr \lp\sum_{i\in I}\int_0^{\beta_i-\alpha_i} \un{\{\widehat W_s^{(i)} \le \min \{\min \widehat W^{(j)},\, j\in I\} + \eps\}}\, ds\rp^p \rbr\\
	&=\int_\K P_{0,1/\sqrt 2}^{6a} (dw) \, \mathbb E_{w} \lbr \lp\sum_{i\in I'}\int_0^{\beta_i-\alpha_i} \un{\{\widehat W_s^{(i)} \le \min \{\min \widehat W^{(j)},\, j\in I'\} + \eps\}}\, ds\rp^p \rbr\\
	&\le 6 \int_\K P_{0,1/\sqrt 2}^{6a} (dw) \, \mathbb E_{w} \lbr \lp\int_0^{I_0} \un{\{\widehat W_s \le \sulW_{I_0} + \eps\}}\, ds\rp^p \rbr = 6\,\Phi_\eps(3a).
\end{align*}

\paragraph{Step 3.} Finally, it remains to treat the two conditionings. We notice that
\begin{align*}
\int_0^1 \un{\{\widehat W_s \le \sulW_1 + \eps\}}\, ds 
			&= \int_0^{I_{\sigma/2}} \un{\{\widehat W_s \le \sulW_{1} + \eps\}}\, ds + \int_{I_{\sigma/2}}^1 \un{\{\widehat W_s \le \sulW_1 + \eps\}}\, ds\\
			&\le \int_0^{I_{\sigma/2}} \un{\{\widehat W_s \le \sulW_{I_{\sigma/2}} + \eps\}}\, ds 
									+ \int_{I_{\sigma/2}}^1 \un{\{\widehat W_s \le \min_{I_{\sigma/2} \le r \le 1} \widehat W_r + \eps\}}\, ds.
\end{align*}
As both terms in the previous line have the same law under $\int_\K \bbB (dw) \, \Pb_w^0(d\omega)$ (notice that, under~$\bbB(dw)$, the processes $(w(t))_{0\le t \le \sigma/2}$ and $(w(\sigma/2+t)-w(\sigma/2))_{0\le t \le \sigma/2}$ have the same law), we obtain
$$\int_\K \bbB (dw) \, \mathbb E_w^0 \lbr \lp\int_0^1 \un{\{\widehat W_s \le \sulW_{1} + \eps\}}\, ds\rp^p \rbr 
		\le 2^{p+1} \int_\K \bbB (dw) \, \mathbb E_w^0 \lbr \lp\int_0^{I_{\sigma/2}} \un{\{\widehat W_s \le \sulW_{I_{\sigma/2}} + \eps\}}\, ds\rp^p \rbr.$$

A consequence of \cite[Equation~(19)]{bettinelli10slr} is that
\begin{align*}
\mathbb E_w^0 \lbr \lp\int_0^{I_{\sigma/2}} \un{\{\widehat W_s \le \sulW_{I_{\sigma/2}} + \eps\}}\, ds\rp^p \rbr
		&=\mathbb E_w \lbr \lp\int_0^{I_{\sigma/2}} \un{\{\widehat W_s \le \sulW_{I_{\sigma/2}} + \eps\}}\, ds\rp^p 
				\ \frac {q_{1 - {I_{\sigma/2}}}(\sigma/2)}{q_1(\sigma)}\, \un{\lb {I_{\sigma/2}}<1\rb}\rbr,
\end{align*}
where
$$q_a(x) \de - \frac x {\sqrt{2\pi a^{3}}}\, \exp\lp -\frac{x^2}{2a}\rp.$$
As a result, we obtain
\begin{align*}
\mathbb E_w^0 \lbr \lp\int_0^{I_{\sigma/2}} \un{\{\widehat W_s \le \sulW_{I_{\sigma/2}} + \eps\}}\, ds\rp^p \rbr
		&\le c_\sigma\ \mathbb E_w \lbr \lp\int_0^{I_{\sigma/2}} \un{\{\widehat W_s \le \sulW_{I_{\sigma/2}} + \eps\}}\, ds\rp^p \rbr,
\end{align*}
where
$$c_\sigma \de \sup_{a > 0} \frac {q_a(\sigma/2)}{q_1(\sigma)} < \8.$$

We then dispose of the second conditioning. By using \cite[Equation~(18)]{bettinelli10slr}, it is not hard to see that
\begin{multline*}
\int_\K \bbB (dw) \, \mathbb E_w \lbr \lp\int_0^{I_{\sigma/2}} \un{\{\widehat W_s \le \sulW_{I_{\sigma/2}} + \eps\}}\, ds\rp^p \rbr\\
		\le \sqrt 2 \int_\K P_{0,\sqrt 3}^{\sigma/2} (dw) \, \mathbb E_{w} \lbr \lp\int_0^{I_0} \un{\{\widehat W_s \le \sulW_{I_0} + \eps\}}\, ds\rp^p \rbr \le 6\sqrt 2\, C'_{3\sigma/2,p,\delta} \, \eps^{4p-\delta}.
\end{multline*}
Setting $c_{p,\delta} \de 2^{p+1} c_\sigma 6\sqrt 2\, C'_{3\sigma/2,p,\delta}$ concludes the proof.
\end{pre}

\subsection{Upper bound for the Hausdorff dimension of \hr{$\partial\qis$}{dq}}\label{secboundaryub}

We may now end the proof of Theorem~\ref{thmdimh}.

\begin{pre}[Proof of Theorem~\ref{thmdimh} (upper bound)]
Under~$\Pb_w$, we set $T_x \de \inf \{ r \ge 0 \, :\, \zeta_r= \zeta(w) - x \}$. For $s \in [0,I_0]$, we set $s^+ \de \sup\{t\, : \, \ulz_t = \ulz_s \}$. Using the same kind of reasoning as in the previous sections, we see that it is enough to show that, for any pseudo-metric~$d$ on $[0,I_0]$ such that
$$d(s,t) \le \widehat W_s + \widehat W_t -2\min_{r\in[s^+, t]} \widehat W_r, \qquad 0 \le s \le  t \le I_0,$$
we have
$$\Pb_w(\dH(\{T_x,\ 0 \le x \le \zeta(w)\},d)\le 2)=1,\qquad \bbB(dw) \text{ a.s.}$$ 

\bigskip

We fix $\eta \in (0,1)$. We will cover $\{T_x,\ 0 \le x \le \zeta(w)\}$ by small open balls. Let us first bound the distance between two points in $\{T_x,\ 0 \le x \le \zeta(w)\}$. Using \cite[Lemma~5.1]{legall07tss} and the fact that, $\bbB(dw)$ a.s., $w$ is $1/(2+\eta)$-H\"older continuous, it is not hard to see that there exists a (random) constant $c < \8$ such that, $\bbB(dw)$ a.s.\ $\Pb_w$ a.s.,
\begin{equation}\label{dzeta}
\lt \widehat W_s - \widehat W_t \rt \le c \, \lp d_\zeta(s,t) \rp^{\frac{1}{2+\eta}}, \qquad \text{ for all }s,t \in [0,I_0],
\end{equation}
where
$$d_\zeta(s,t) \de \zeta_s + \zeta_t -2\min_{r\in[s \wedge t, s \vee t]} \zeta_r, \qquad s,t\in [0,I_0].$$
Let $0 \le x \le y \le \zeta(w)$, and $m(x,y) \in [T_x^+,T_y]$ be such that $\widehat W_{m(x,y)}=\min_{s\in[T_x^+,T_y]}\widehat W_s$. When~\eqref{dzeta} holds, we have
\begin{align}
d(T_x,T_y) &\le \widehat W_{T_x} + \widehat W_{T_y} -2\widehat W_{m(x,y)} \notag\\
		 &\le c \, \lp \big( d_\zeta(T_x,m(x,y)) \big)^{\frac{1}{2+\eta}} + \big( d_\zeta(T_y,m(x,y)) \big)^{\frac{1}{2+\eta}} \rp\notag\\
		 &\le 2c \, \big( y - x +2\big(\zeta_{m(x,y)}- \ulz_{m(x,y)} \big) \big)^{\frac{1}{2+\eta}}\label{flbound}.
\end{align}

In order to control the term $(\zeta_{m(x,y)}- \ulz_{m(x,y)})$ in the above inequality, we sort out the excursions going ``too high.'' Namely, we fix $\eps>0$, and set
$$I^{(\eps)} \de \lb i \in I \, : \, \sup_{s\ge 0} \zeta^{(i)}_s > \eps \rb.$$
By Lemma~\ref{poisexc}, the cardinality of~$I^{(\eps)}$ is under~$\Pb_w$ a Poisson random variable with mean
$$2 \int_0^{\zeta(w)} dt\, \N_{w(t)}\lp \sup_{s\ge 0} \zeta_s > \eps \rp= 2\zeta(w)\ n\lp \sup_{s\ge 0} e_s > \eps \rp = \frac{\zeta(w)}{\eps}.$$
In particular, $|I^{(\eps)}| < \8$, $\Pb_w$ a.s. We denote by $B(s,r)\subseteq [0,I_0]$ the open ball of radius~$r$ centered at~$s$, for the pseudo-metric~$d$, and, for $i \in I$, we set $x_i \de \zeta(w) - \zeta_{\alpha_i}$. If $\delta \de 2c\, (3\eps)^{1/(2+\eta)}$, we claim that the set
$$\big\{ B\lp T_{x_i},\delta \rp, \ i \in I^{(\eps)} \big\} \cup \lb B\lp T_{k\eps},\delta \rp, \ 0 \le k \le \frac{\zeta(w)}{\eps} \rb$$
is a covering of $\{T_x,\ 0 \le x \le \zeta(w)\}$. To see this, let us take a point $y \in [0, \zeta(w)]$, and let us consider $x\de \max\{ s \in \{0\}\cup\{x_i,\ i \in I^{(\eps)}\}\, : \, s \le y \}$. Observe that, by construction, all the excursions of $\zeta-\ulz$ with support included in the interval $[T_x^+, T_y]$ have a height smaller than~$\eps$. Because $T_x^+ \le m(x,y) \le T_y$, we see that $\zeta_{m(x,y)}- \ulz_{m(x,y)} \le \eps$. Then, if $y-x < \eps$, by~\eqref{flbound}, we have $T_y \in B(T_x,\delta)$. If $y-x \ge \eps$, then $y-\lf y/\eps \rf \eps < \eps$, and $\lf y/\eps \rf \eps \ge x$, so that $\zeta_{m(\lf y/\eps \rf \eps,y)}- \ulz_{m(\lf y/\eps \rf \eps,y)} \le \eps$. This yields that $T_y \in B(T_{\lf y/\eps \rf \eps},\delta)$, by~\eqref{flbound}.

The $(2+\eta)(1+\eta)$-value of this covering is less than 
$$\lp |I^{(\eps)}| +  \frac{\zeta(w)}{\eps} + 1 \rp \, (2\delta)^{(2+\eta)(1+\eta)} \le c'\, \lp |I^{(\eps)}| +  \frac{\zeta(w)}{\eps} + 1 \rp \, \eps^{1+\eta},$$
for some constant~$c'$, independent of~$\eps$. By Chebyshev's inequality, we see that with $\Pb_w$-probability at least $1-\eps/\zeta(w)$, we have $|I^{(\eps)}| \le 2\zeta(w)/\eps$. We conclude that, $\bbB(dw)$ a.s., the $(2+\eta)(1+\eta)$-Hausdorff content of $\{T_x,\ 0 \le x \le \zeta(w)\}$ is $\Pb_w$ a.s.\ equal to~$0$, so that $\dH(\{T_x,\ 0 \le x \le \zeta(w)\},d)\le (2+\eta)(1+\eta)$. Finally, letting $\eta \to 0$ yields the result.
\end{pre}

\section{Developments and open questions}\label{secdev}

\subsection{Quadrangulations with a simple boundary}

We considered in this work quadrangulations with a boundary that is not necessarily simple. It is natural to ask ourselves what happens if we require the boundary to be simple. This translates into a conditioning of the coding forest from which some technical difficulties arise and our results may not straightforwardly be adapted to this case. We expect, however, to find the same limit, up to some factor modifying the length of the boundary. 

Very roughly, the intuition is that the 0-regularity of the boundary implies that a large quadrangulation with a boundary should typically consist in one large quadrangulation with a simple boundary on which small quadrangulations with a boundary are grafted. As a result, if we remove these small components on the boundary, the first quadrangulation should not be too far from a quadrangulation with a simple boundary having roughly the same number of faces but a significantly smaller boundary.

We expect that such results may be rigorously derived from our analysis with a little more work.

\subsection{Application to self-avoiding walks}

\begin{figure}[ht]
		\psfrag{a}[][]{\textit{no cycles}}
		\psfrag{1}[][]{($\times 1$)}
		\psfrag{s}[][]{($\times (\sigma-1)$)}
		\psfrag{n}[][]{($\times n$)}
	\centering\includegraphics{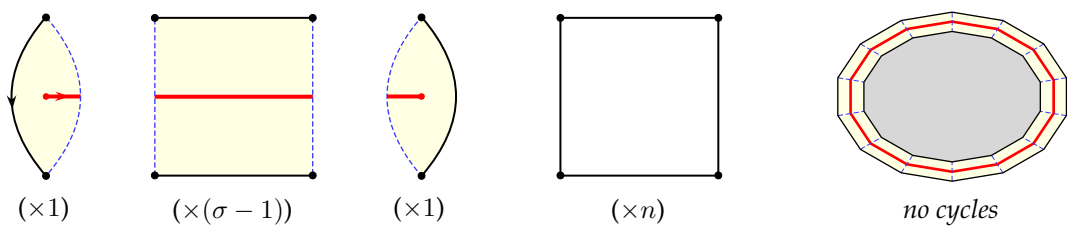}
	\caption[Borot, Bouttier and Guitter's model of self-avoiding walks.]{Borot, Bouttier and Guitter's model of self-avoiding walks. The distinguished edges are the thin dashed (blue) lines, whereas the other edges are the thicker black solid lines. The path is drawn in an even thicker solid (red) line. It starts in the root half-step tile and ends in the other half-step tile. We added an arrowhead on it in the root half-step tile to symbolize this fact.}
	\label{saw_model}
\end{figure}

We present here a model of self-avoiding walks on random quadrangulations, which is adapted from~\cite{bouttier11on}. In the latter reference, Borot, Bouttier and Guitter study a model of loops on quadrangulations. We can easily adapt their model to the case of self-avoiding walks, and we see that it is directly related to quadrangulations with a boundary. We call \textbf{step tile} a quadrangle in which two opposite half-edges incident to the quadrangle are distinguished, and \textbf{half-step tile} a face of degree~$2$ in which one incident half-edge is distinguished. On the figures, we draw a (red) line linking the two distinguished edges in the step tiles, as well as a line linking the distinguished edge to the center of the face in the half-step tiles (see Figure~\ref{saw_model}). These lines will constitute the self-avoiding walk of the model. 

Let $n \ge 0$ and $\sigma \ge 1$ be two integers. A map whose faces consist in two half-step tiles, $\sigma-1$ step tiles, and~$n$ quadrangles is called an \textbf{($n$,$\sigma$)-configuration} if it satisfies the following:
\begin{itemize}
	\item the reverse of every distinguished half-edge is also a distinguished half-edge,
	\item there are no cyclic chains of step tiles,
	\item the root of the map is the half-edge that is not distinguished in one of the two half-step tiles.
\end{itemize}

\begin{figure}[ht]
	\centering\includegraphics{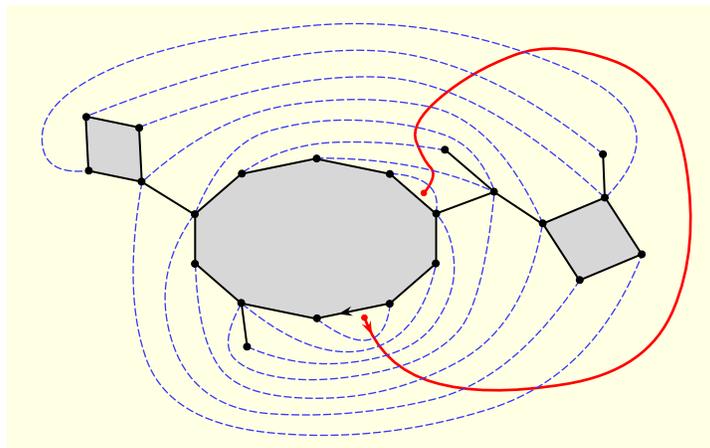}
	\caption[A quadrangulation with a boundary and the corresponding configuration.]{A quadrangulation with a boundary and the corresponding configuration. On this figure, $\sigma=15$.}
	\label{saw_bij}
\end{figure}

We claim that the ($n$,$\sigma$)-configurations are in one-to-one correspondence with the quadrangulations with a boundary having~$n$ internal faces and~$2\sigma$ half-edges on the boundary. Indeed, let us take an ($n$,$\sigma$)-configuration. It has~$2\sigma$ distinguished half-edges forming~$\sigma$ different edges. When removing these~$\sigma$ edges, we obtain a map having~$n$ quadrangles and one other face of degree $1+2(\sigma-1)+1=2\sigma$, this face being incident to the root. This is thus a quadrangulation of~$\Qns$. Conversely, let us consider a quadrangulation of~$\Qns$, and let $\e_1=\e_*$, $\e_2$, \dots, $\e_{2\sigma}$ be the half-edges incident to its external face, read in the clockwise order. We add extra edges (that do not cross each other) linking~$\e_{2\sigma-i}^+$ to~$\e_{i+1}^+$ for $0 \le i \le \sigma-1$. We thus create two faces of degree~$2$ (one incident to~$\e_*$ and one incident to~$\e_{\sigma+1}$) as well as $\sigma-1$ faces of degree~$4$. These faces are half-step tiles and step tiles when we distinguish the half-edges composing the extra edges we added. It is then easy to see that this map is an ($n$,$\sigma$)-configuration (see Figure~\ref{saw_bij}). Moreover, the composition of the two operations we described here (in an order or the other) is clearly the identity, so that our claim follows.

We believe that some results may be derived from our work and this bijection. For example, we think possible to show that, when $\sigma_n/\sqrt{2 n} \to \sigma\in (0,\8)$, a uniform ($n$,$\sigma_n$)-configuration converges (up to extraction) toward a metric space with a marked path in some sense. Moreover, the limiting space should be homeomorphic to the $2$-dimensional sphere and have dimension~$4$ a.s. The marked path should also have dimension~$2$.

Here again, however, some technical difficulties arise. The main problem is that the gluing operation tremendously modifies the metric. We may easily define this gluing operation in the continuous setting by considering the quotient of~$\qis$ by the coarsest equivalence relation for which $\fl(x)\sim \fl(\sigma-x)$, $x\in [0,\sigma]$ (with the notation of the end of Section~\ref{secsurfb}), but some care is required when dealing with this quotient. In particular, it is not clear that the points identified in this quotient are solely the one identified by the equivalence relation. Then, it also remains to show that the convergence still holds after the gluing operation.

\bigskip

Understanding the scaling limit of quadrangulations with a simple boundary and this gluing operation could also lead to some interesting results on random self-avoiding walks on random quadrangulations, as quadrangulations with a marked self-avoiding walk are in one-to-one correspondence with quadrangulations with a simple boundary by a bijection similar to the one described above.

We end this section by mentioning that the model we presented corresponds to the particular case Borot, Bouttier and Guitter called \textit{rigid} in~\cite{bouttier11on}. We may also consider the more general case in which the two distinguished half-edges of a step tile are not required to be opposite. In this model, we also expect to find the same limits.

\nocite{*}
\bibliographystyle{alpha}
\bibliography{_slrpqb_bib}
\end{document}